\documentclass[oneside,english]{amsart}
\usepackage[T1]{fontenc}
\usepackage[latin9]{inputenc}
\usepackage{amsthm}
\usepackage{amssymb}
\usepackage[pdftex]{graphicx}
\usepackage{bm}

\makeatletter
\numberwithin{equation}{section}
\numberwithin{figure}{section}

\makeatother

\usepackage{babel}
\begin{document}

\title[\tiny{Dirichlet forms on self-similar sets with overlaps}]{Dirichlet forms on self-similar sets with overlaps}
\author{Shiping Cao}
\address{Department of Mathematics, Cornell University, Ithaca, NY 14850, U.S.A.}
\curraddr{} \email{sc2873@cornell.edu}
\thanks{}

\author{Hua Qiu}
\address{Department of Mathematics, Nanjing University, Nanjing, 210093, P. R. China.}
\curraddr{} \email{huaqiu@nju.edu.cn}
\thanks{The research of the second author was supported by the Nature Science Foundation of China, Grant 11471157.}



\date{}

\dedicatory{}
\begin{abstract} We study Dirichlet forms and Laplacians on self-similar sets with overlaps. A notion of   ``finitely ramified of finite type($f.r.f.t.$) nested structure'' for self-similar sets is introduced. It allows us to reconstruct a class of self-similar sets in a graph-directed manner by a modified setup of Mauldin and Williams, which satisfies the property of finite ramification. This makes it  possible to extend the technique developed by Kigami for analysis on $p.c.f.$ self-similar sets to this more general framework. Some basic properties related to $f.r.f.t.$ nested structures are investigated. Several non-trivial examples and their Dirichlet forms are provided.

\end{abstract}
\maketitle

\section{Introduction }
Analysis, especially the theory of Laplacians, on fractals has been extensively developed on certain self-similar fractals, see [K1-K7, S1-S2] and the references therein. To define the Dirichlet forms and Laplacians directly and constructively, it typically requires that the fractals have the property of finite ramification, which means that any connected subset of the fractals can be disconnected by removing finitely many appropriate points. The class of $p.c.f.$(\textit{post-critically finite}) self-similar sets introduced by Kigami [K2] satisfies perfectly this requirement.
Let $K$ be a \textit{p.c.f.} self-similar set, satisfying the self-similar identity $K=\bigcup_{i=1}^N F_iK$ with $N\geq 2$ and $\{F_i\}_{i=1}^N$ being an iterated function system($i.f.s.$ for short) of contractive similitudes. The union of  intersections of  cells of level $1$,
\[
C_K=\bigcup_{1\leq i<j\leq N}F_iK\cap F_jK,
\]
is a finite set and disconnects the fractal $K$ into small pieces, i.e.,
\begin{equation}
K\setminus C_K=\bigsqcup_{i=1}^N (F_iK\setminus C_K),
\end{equation}
where we use ``$\bigsqcup$'' to denote the disjoint union. 
In addition, there is a finite subset $V_0$ in $K$, called the \textit{boundary} of $K$, satisfying
\begin{equation}
V_0\cup C_K\subset \bigcup_{i=1}^N F_i V_0.
\end{equation}

The above characteristics of $p.c.f.$ self-similar sets are essential in the construction of self-similar Dirichlet forms. Similar descriptions  can be applied to some other finitely ramified fractals, for example, finitely ramified graph-directed fractals [HMT, HN, M1, MG] and some Julia sets [ADS, FS, RT, SST]. 

\begin{figure}[h]
	\centering
	\includegraphics[height=3cm]{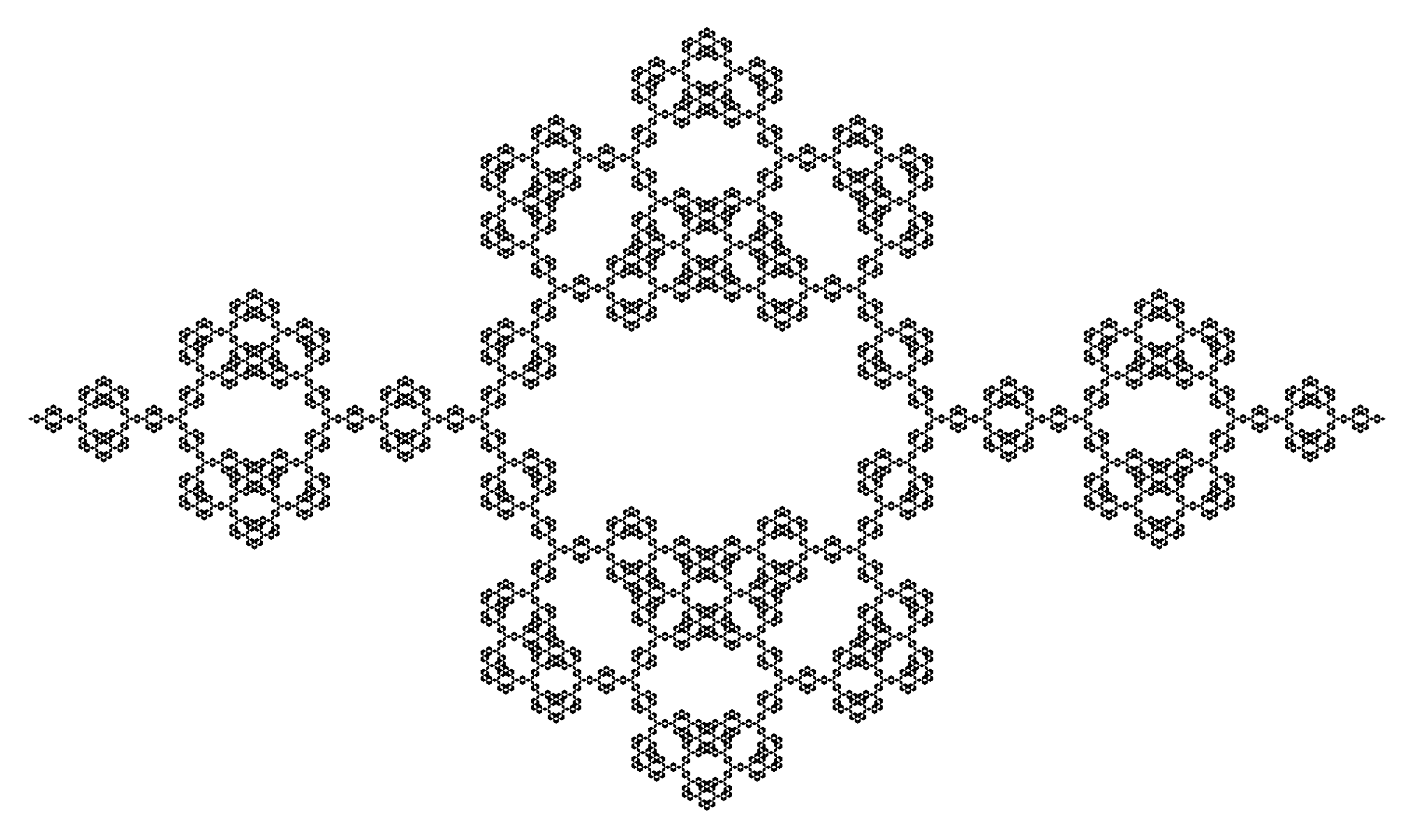}
	\begin{center}
		\caption{The diamond fractal.}
	\end{center}
\end{figure}

It is desirable to enlarge the class of self-similar sets on which the sprit of Kigami's construction of Dirichlet forms works. In this paper, we will focus our interest on the analysis on self-similar sets allowing overlaps(it may happen that $\# F_iK\cap F_jK=\infty$ for some $i\neq j$) and satisfying the finitely ramified requirement. One such example that is not $p.c.f.$ and has been well studied is the \textit{diamond fractal} [KSW, M2], see Figure 1.1. 
However, there are few other examples solved. Our main objective is to set up a condition for self-similar sets which is as general as possible and will include $p.c.f.$ self-similar sets, so that we can construct Dirichlet forms through simple extensions of the techniques developed previously. The main idea is based on the same device as that for the diamond fractal, i.e., a finitely ramified graph-directed reconstruction of the fractal, which is a modified setup of Mauldin and Williams [MW]. We will show that this strategy can be adapted to a general class of self-similar sets. Note that we will work in $\mathbb{R}^n$, but it would be possible to extend analogous results to abstract metric spaces. We will give several non-trivial examples with solutions of Dirichlet forms, see Figure 1.2, 1.3. 

\begin{figure}[h]
	\centering
	\includegraphics[height=3.7cm]{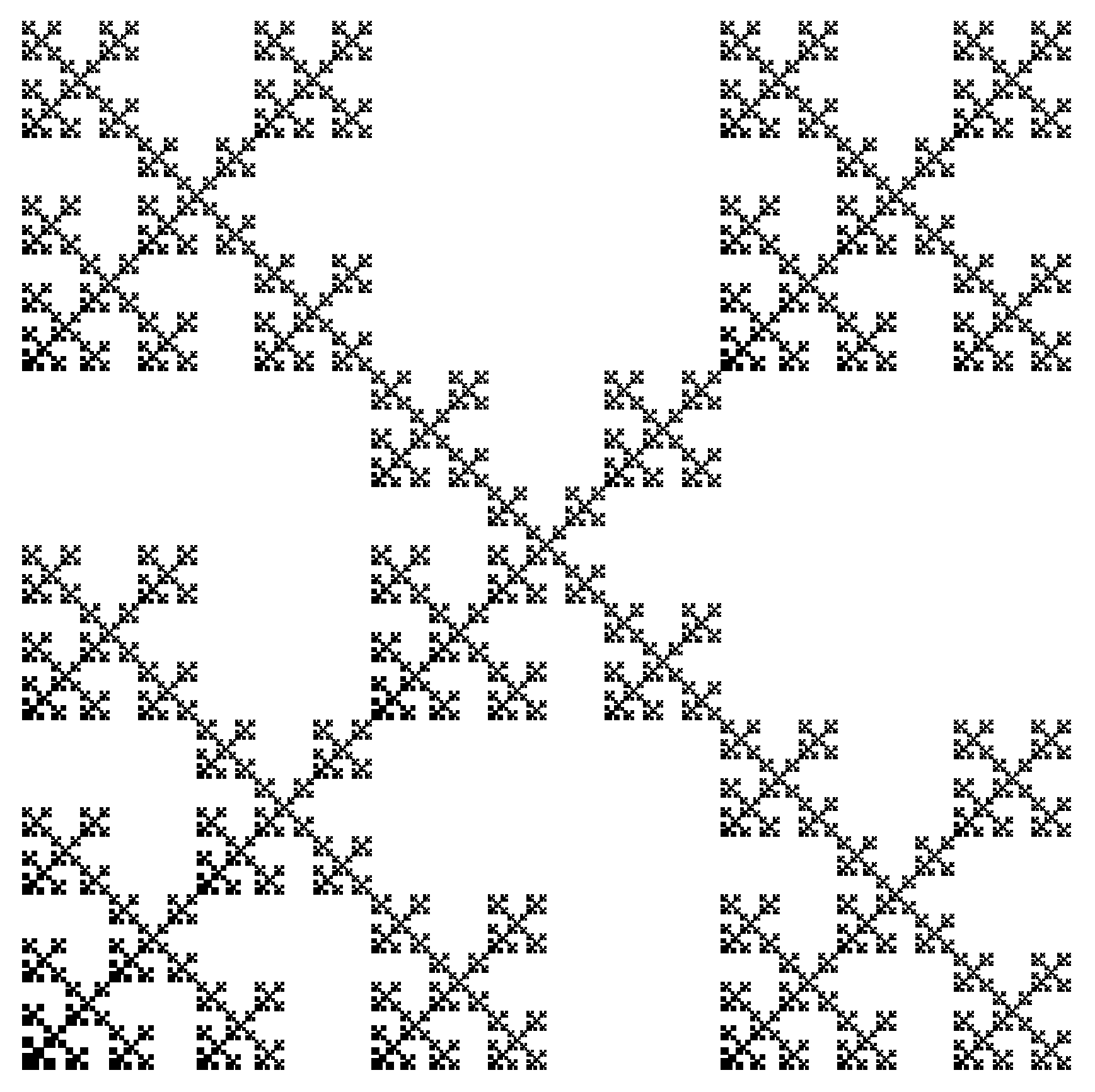}\qquad
	\includegraphics[height=3.7cm]{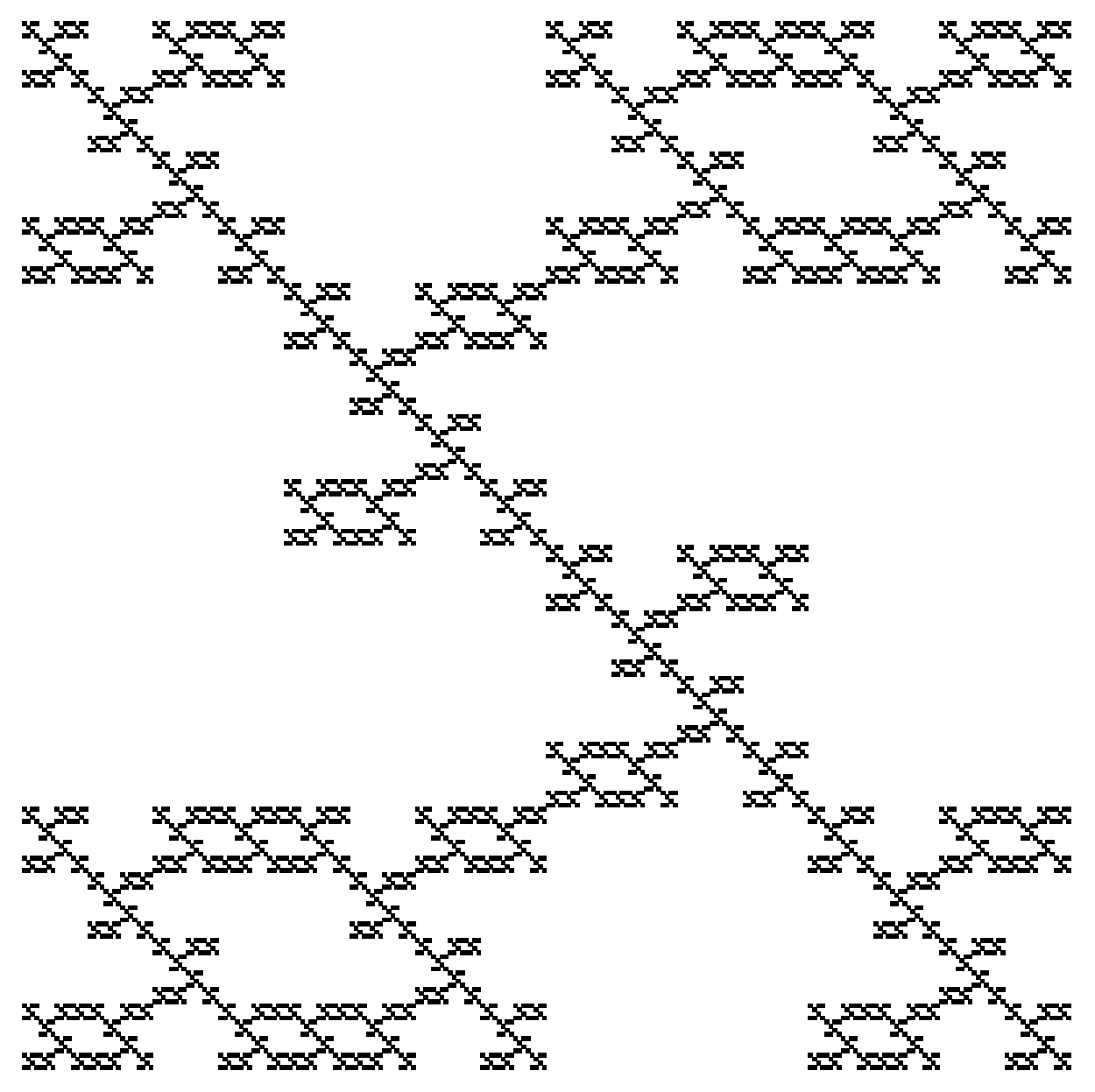}
	\begin{center}
		\caption{Vicsek type fractals with overlaps.}
	\end{center}
\end{figure}

\begin{figure}[h]
	\centering
	\includegraphics[height=3.7cm]{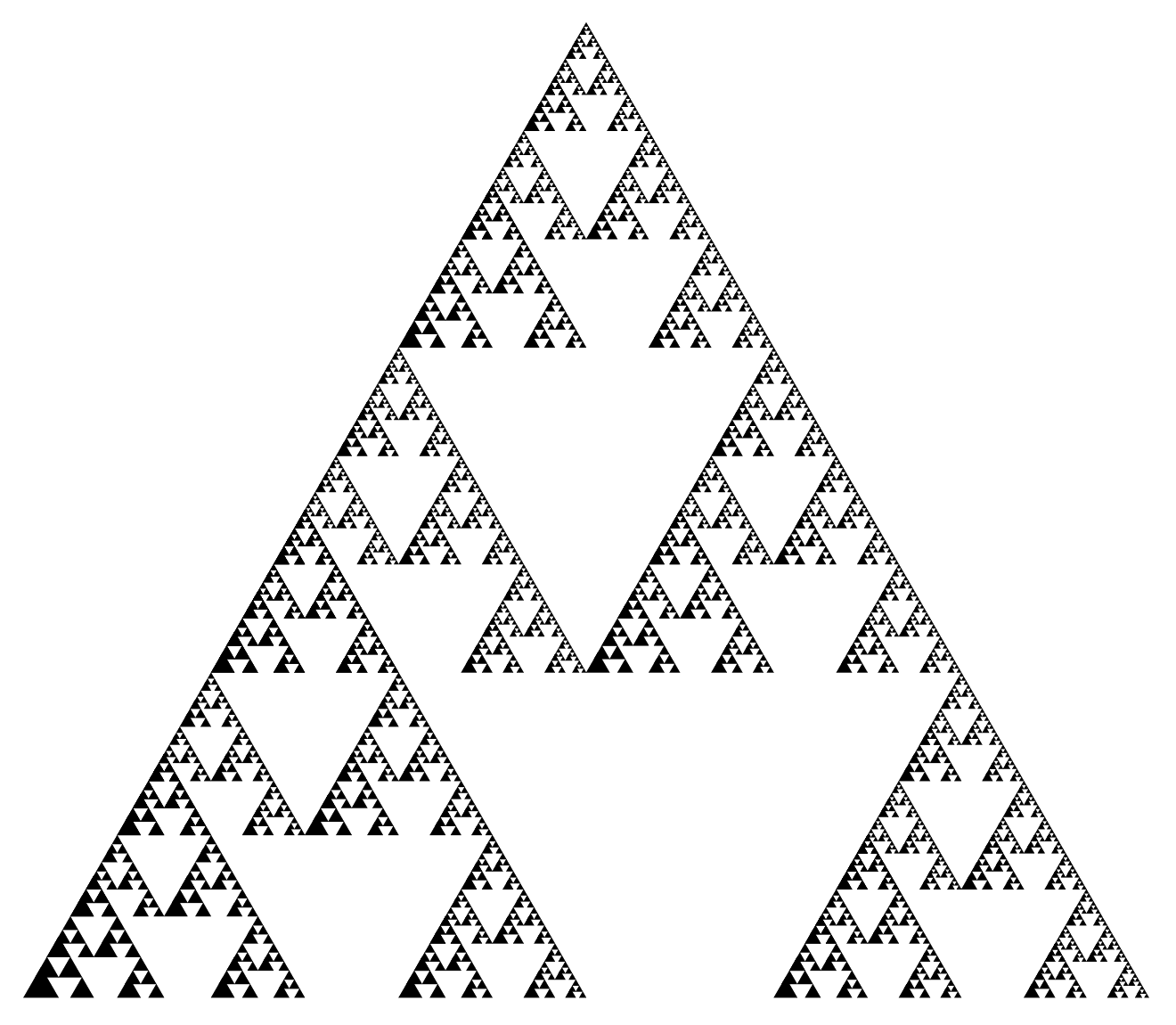}\qquad
	\includegraphics[height=3.7cm]{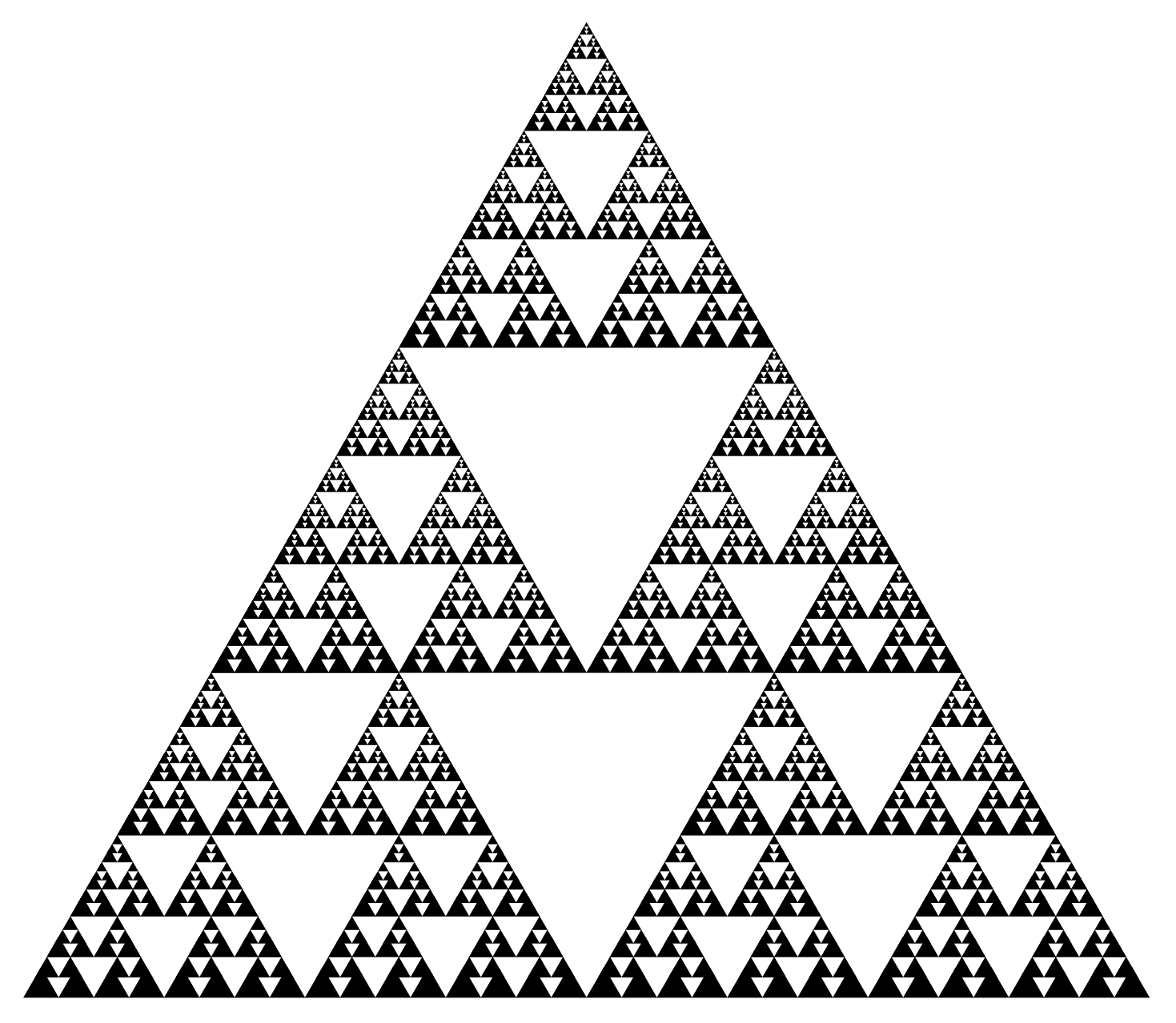}
	\begin{center}
		\caption{Sierpinski gasket type fractals with overlaps.}
	\end{center}
\end{figure}

To be specific, we will introduce a class of self-similar sets allowing overlaps, called \textit{finitely ramified of finite type}($f.r.f.t.$ for short) self-similar sets. For a self-similar set $K$, we always assume that it is connected. Call a connected compact subset $A$ in $K$ an \textit{island}  if $\#A\cap \overline{K\setminus A}<\infty$. Roughly, we need a class of islands(including $K$ itself)  in $K$ satisfying the following two requirements. The first one is each island $A$ can be split into a finite union of at least two islands, called the \textit{children} of $A$. Say two islands are \textit{equivalent} if they are similar, and we always choose the same way to split equivalent islands.  Then the second requirement is that there are finite types of equivalent islands involved. We will introduce an index set $\Lambda$  to help us to describe a structure of nested islands, called $f.r.f.t.$ nested structure, on $K$ to characterize the two requirements rigidly. The class of $p.c.f.$ self-similar sets is included in $f.r.f.t.$ self-similar sets, where the structure of islands consists of all similar copies of $K$ generated by the $i.f.s.$ and the number of types of such islands is $1$. Essentially, an $f.r.f.t.$ nested structure of $K$ will provide us a graph-directed construction of $K$ having the property of finite ramification. This framework allows us to apply Kigami's technique dealing with Laplacians on $p.c.f.$ self-similar sets to $f.r.f.t.$ self-similar sets through simple extensions. Then the problem of finding Dirichlet form and Laplacian on $K$ is to investigate the existence and uniqueness of a fixed point for a certain multiple dimensional renormalization mapping. This requires a detailed analysis of the overlaps when the similitudes are iterated. 

Hambly and Nyberg [HN] have used the same idea to consider the so-called \textit{finitely ramified graph-directed($f.r.g.d.$ for short) fractal families}  and study the spectral asymptotics of eigenvalue counting functions of the associated Laplacians by establishing a multidimensional renewal theorem. It should be pointed out that the concept of $f.r.f.t.$ self-similar sets is in fact an intermediate one between $p.c.f.$ self-similar sets and $f.r.g.d.$ fractals.  

There are some basic questions arisen naturally related to $f.r.f.t.$ self-similar sets. 

The first one is when a self-similar set $K$ possesses an $f.r.f.t.$ nested structure. The $f.r.f.t.$ requirement may allow overlaps when the similitudes are iterated, but it seems that the overlapping types among distinct comparable(with respect to diameters) similar copies of $K$ should be finite. On the other hand, it is natural to formulate each island by a finite union of comparable similar copies of $K$, and if so, it is reasonable to believe that the number of copies of $K$ that make up of each island  should have a uniform control. Basing on these observations, we propose two conditions for $K$. One is the \textit{finite neighboring type property}, the other is the \textit{finite chain length property}. They can not imply each other. We prove that $K$ possesses an $f.r.f.t.$ nested structure if both these two conditions hold. These two conditions are quite general, at least they hold for all known examples. But we do not know whether they are necessary for the $f.r.f.t.$ requirement.

Once the self-similar set $K$ is $f.r.f.t.$, it will have infinitely many $f.r.f.t.$ nested structures.  Let $\mathcal{S}$ and $\mathcal{S}'$ be two $f.r.f.t.$ nested structures on $K$. There are two possibilities. We say $\mathcal{S}$ is \textit{derived} from $\mathcal{S}'$ if each island in $\mathcal{S}$ is an iteration of islands in $\mathcal{S}'$ under the graph-directed construction of $\mathcal{S}'$. For example, look at a $p.c.f.$ self-similar set $K$ associated with an $i.f.s.$ $\{F_i\}_{i=1}^N$. Then $\{F_i\circ F_j\}_{i,j=1}^N$ is also an $i.f.s.$ of $K$ which consists of $N^2$ possible similitudes by iterating $F_i$'s twice. These two $i.f.s.$'s will naturally induce two $f.r.f.t.$ nested structures. Obviously, the latter one is derived from the former one. It may also happen that $\mathcal{S}$ and $\mathcal{S}'$ are not derived from each other, see Example 1 in Section 4 for example. It is naturally to consider the relationship between Dirichlet forms associated with distinct $f.r.f.t.$ nested structures. We will restrict to consider those Dirichlet forms with a homogenous property, i.e., the resulting Laplacians are locally translation invariant, and prove that under a mild condition, different $f.r.f.t.$ nested structures will give rise to same homogenous Dirichlet forms.

The ``finite type'' assumption is quite useful for calculating the Hausdorff dimension of a self-similar set $K$ with overlaps. By reconstructing $K$ in a graph-directed manner, one can determine the Hausdorff dimension of $K$ in terms of the spectral radius of certain weighted incidence matrix. See [JY, L, LN, RW] and the references therein. In [LN], Lau and Ngai formulated a so-called \textit{ generalized finite type condition} for self-similar sets, which extends  the well-known \textit{open set condition}.  The  $p.c.f.$ condition and open set condition are two distinct separation conditions for different purposes of investigation. We remark that the $f.r.f.t.$ condition and generalized finite type condition are  natural extensions of the $p.c.f.$ condition and the open set condition  in the sense of Mauldin and Williams, respectively.

The remainder of this paper is organized as follows.

Firstly, the topology of $f.r.f.t.$ nested structures will be discussed from Section 2 to Section 4. In Section 2, we will introduce the definition of the $f.r.f.t.$ nested structures, and discuss some basic topological properties of these structures. In Section 3, we will explore various conditions for a self-similar set to have an $f.r.f.t.$ nested structure. We will investigate the relationship between $f.r.f.t.$ self-similar sets and other well established finitely ramified fractals including $p.c.f.$ self-similar sets, and $f.r.g.d.$ fractals.  At the end of this section, we will introduce the finite neighboring type property and finite chain length property  for a self-similar set $K$, and prove that they are sufficient for $K$ to possess an $f.r.f.t.$ nested structure. Then in Section 4, we will describe the $f.r.f.t.$ nested structures of the four examples shown in Figure 1.2. and Figure 1.3.

Secondly, we will deal with the Dirichlet forms associated with the $f.r.f.t.$ nested structures in Section 5 and Section 6. In Section 5, we will introduce the concept of harmonic structures for $f.r.f.t.$ nested structures analogous to that for $p.c.f.$ self-similar sets [K2] and $f.r.g.d.$ fractals [HN], which are solutions of canonical fixed-point problems. We will construct Dirichlet forms on the four examples shown in Section 4. In Section 6, for an $f.r.f.t.$ self-similar set $K$, we will restrict to consider the homogeneous harmonic structures, which together with the proper normalized Hausdorff measure will induce homogenous Dirichlet forms on $K$. We prove that under some assumption, different $f.r.f.t.$ nested structures on $K$ will induce same homogeneous Dirichlet forms.

Thirdly, in Section 7, we will briefly go through the spectral asymptotics  of the constructed Laplacians  by virtue of the normalized Hausdorff measure on a $f.r.f.t.$ self-similar set $K$. The result in this section is a direct application of the result developed in [HN]. 
 
Lastly, in Appendix we will present two interesting examples, which are finitely ramified  but not $f.r.f.t.$ self-similar sets, to show that the finite neighboring type property and finite chain length property, proposed in Section 3, can not imply each other. The analysis on these fractals are not clear.

\section{$f.r.f.t.$ nested structures}
\textbf{Definition 2.1.} \textit{Let $K$ be a connected self-similar set and $\{K_\alpha\}_{\alpha\in\Lambda}$ be a countable collection of distinct compact connected subsets in $K$, containing no singleton, satisfying that}

\textit{1. there is an index $ \vartheta\in \Lambda$, called the root of $\Lambda$, such that $K=K_\vartheta$;}

\textit{2. for any $\alpha\in\Lambda$, there is a finite set $\Lambda_\alpha\subset \Lambda$ with $\#\Lambda_\alpha\geq 2$, such that $K_\alpha=\bigcup_{\beta\in\Lambda_\alpha}K_\beta$, call $\alpha$ the parent of $\beta$, and $\beta$ the child of $\alpha$; }

\textit{3.  any  $\alpha\in\Lambda\setminus\{\vartheta\}$ is an offspring of $\vartheta$. }

\textit{Call $\{K_\alpha,\Lambda_\alpha\}_{\alpha\in \Lambda}$ a} nested structure \textit{of $K$, and $\Lambda$ its} index set.

Note that for any index $\alpha\in\Lambda\setminus\{\vartheta\}$, it has at least one parent. We denote by $P^1(\alpha)$\big(or simply $P(\alpha)$\big) the set of parents of $\alpha$. Trivially, write $P(\vartheta)=\emptyset$. For $k\geq 2$, write $P^k(\alpha)=\bigcup_{\beta\in P^{k-1}(\alpha)} P(\beta)$ inductively.  Let $\Lambda^{(0)}_\vartheta=\{\vartheta\}$ and for $k\geq 1$, let $\Lambda^{(k)}_\vartheta=\bigcup_{\alpha\in \Lambda^{(k-1)}_\vartheta} \Lambda_\alpha$ inductively. It is easy to see that $\Lambda=\bigcup_{k=0}^\infty \Lambda^{(k)}_\vartheta$.

\textbf{Remark.} \textit{Suppose we have a structure $\{K_\alpha,\Lambda_\alpha\}_{\alpha\in \Lambda}$ which only satisfies condition $1,2$ in Definition 2.1, then there exists a subset $\Lambda'\subset \Lambda$ such that $\{K_\alpha,\Lambda_\alpha\cap \Lambda'\}_{\alpha\in \Lambda'}$ is a nested structure of $K$.}

\textit{Proof.} For any $\alpha\in\Lambda$, we still denote $P(\alpha)$ the set of parents of $\alpha$. Notice that $P(\alpha)$ may be an empty set. 

We will show that there is a minimal subset $\Lambda'\subset \Lambda$ such that $\{K_\alpha,\Lambda_\alpha\cap \Lambda'\}_{\alpha\in \Lambda'}$ satisfies condition $1,2$, with $\Lambda_\alpha$ and  $\Lambda$ replaced by $ \Lambda_\alpha\cap\Lambda'$ and $\Lambda'$ respectively, and the minimality of $\Lambda'$ will imply condition $3$. 

First, the existence of a minimal subset $\Lambda'$ is ensured by the Zorn's lemma. In fact, if $\{\Lambda^{(i)}\}_{i\in I}$ is a chain of decreasing nested subsets of $\Lambda$ such that $\{K_\alpha,\Lambda_\alpha\cap\Lambda^{(i)}\}_{\alpha\in \Lambda^{(i)}}$ satisfies condition $1,2$, then  the structure $\{K_\alpha,\Lambda_\alpha\cap\Lambda^{(\infty)}\}_{\alpha\in\Lambda^{(\infty)}}$  with $\Lambda^{(\infty)}=\bigcap_{i\in I} \Lambda^{(i)}$ also satisfies condition $1,2$, since for each $\alpha\in \Lambda^{(\infty)}$, we have $ \Lambda_\alpha\cap\Lambda^{(\infty)}=\Lambda_\alpha\cap \Lambda^{(i)}$ for all $i\geq i_0$, for some $i_0\in I$.

Next, suppose the minimal subset $\Lambda'$ does not satisfy condition $3$, then there exists $\gamma\in \Lambda'$ such that $\gamma\notin\bigcup_{k=0}^\infty (\Lambda')^{(k)}_\vartheta$. It is easy to see $P(\gamma)\cap \bigcup_{k=0}^\infty (\Lambda')^{(k)}_\vartheta=\emptyset$, and inductively for any $m\geq 1$,  $P^m(\gamma)\cap \bigcup_{k=0}^\infty (\Lambda')^{(k)}_\vartheta=\emptyset$. So $\Lambda'\setminus\bigcup_{m=0}^\infty P^m(\gamma)$ is a proper subset in $\Lambda'$ satisfying condition $1,2$, which is a contradiction to the minimality of $\Lambda'$. \hfill$\square$\\

For a connected self-similar set $K$(allowing overlaps), let $N\geq 2$ and  $\{F_i\}_{i=1}^{N}$ be its $i.f.s.$. For $n\geq 1$, denote $W_n=\{1,2,\cdots,N\}^n$ the set of words of length $n$, together with $W_0=\{\vartheta\}$, and $W_*=\bigcup_{n\geq 0}W_n$. For $w=w_1w_2\cdots w_n\in W_*$, write $F_w=F_{w_1}\circ F_{w_2}\circ\cdots\circ F_{w_n}$ for short. Note that it is possible that $F_w=F_{w'}$ with $w\neq w'$. By removing all but the smallest in the lexicographical order (or any fixed order) words, we obtain an index set $W_\#\subset W_*$ such that $\{F_wK\}_{w\in W_\#}$ consists of distinct copies of $K$ and $\{F_wK\}_{w\in W_\#}=\{F_wK\}_{w\in W_*}$. Then it is easy to see that $\{F_w K, \Lambda_w\}_{w\in W_\#}$ is a  \textit{canonical} nested structure of $K$,  where 
$
\Lambda_w=\{u\in W_\#\cap W_{n+1}: F_u K\subset F_wK\}
$
 for $w\in W_n$.

We list three requirements that need to be assumed on the structure  $\{K_\alpha,\Lambda_\alpha\}_{\alpha\in\Lambda}$.

For $\alpha,\beta\in\Lambda$, write $\alpha\sim\beta$ if there exists a similitude $\phi_{\alpha,\beta}$ such that  $\phi_{\alpha,\beta}(K_\alpha)=K_\beta.$ Fix the similitude $\phi_{\alpha,\beta}$ so that $\phi_{\alpha,\alpha}=id$ and $\phi_{\gamma,\beta}\circ\phi_{\alpha,\gamma}=\phi_{\alpha,\beta}, \forall\alpha\sim\beta, \beta\sim\gamma$. Denote $\Lambda/\sim$ the collection of equivalent classes in $\Lambda$ with respect to ``$\sim$''. Let $C_\alpha=\bigcup_{\beta,\beta'\in\Lambda_\alpha} K_\beta\cap K_{\beta'}$ and define
\[
V_\alpha=\bigcup_{\beta\sim\alpha,\beta\neq\vartheta}\phi_{\beta,\alpha}\Big(K_\beta\cap \big(\bigcup_{n\geq 1}\bigcup_{\gamma\in P^n(\beta)}C_{\gamma}\big)\Big).
\]

\textbf{A1.} \textit{
	Assume $\#(\Lambda/\sim)<\infty$. }

\textbf{A2.} \textit{For $\alpha\sim\alpha'$, there is a one to one correspondence between $\Lambda_\alpha$ and $\Lambda_{\alpha'}$ such that  $\forall \beta\in\Lambda_\alpha$, there exists a unique $\beta'\in \Lambda_{\alpha'}$ satisfying $\beta\sim\beta'$ and $\phi_{\beta,\beta'}=\phi_{\alpha,\alpha'}$.}

\textbf{A3.} \textit{ 
Assume $\# V_\alpha<\infty$ for all $\alpha\in \Lambda$.}

\textbf{Definition 2.2.} \textit{We say a  nested structure $\{K_\alpha,\Lambda_\alpha\}_{\alpha\in \Lambda}$ is} finitely ramified of finite type\textit{($f.r.f.t.$ for short) if A1, A2 and A3 are satisfied, and call $K$ an} $f.r.f.t.$ self-similar set.

We denote by $M=\# (\Lambda/\sim)$ for simplicity, and called it the \textit{number of types} of  $\{K_\alpha,\Lambda_\alpha\}_{\alpha\in \Lambda}$.

 It is easy to prove that for a $p.c.f.$ self-similar set, $W_*=W_\#$ and the  canonical  structure $\{F_wK,\Lambda_w\}_{w\in W_*}$ is $f.r.f.t.$ with $M=1$. The proof will be given in Section 3.  

In general, an $f.r.f.t.$  nested structure $\{K_\alpha,\Lambda_\alpha\}_{\alpha\in \Lambda}$ will have analogous properties as (1.1) and (1.2).

\textbf{Proposition 2.3.} \textit{Let  $\{K_\alpha,\Lambda_\alpha\}_{\alpha\in \Lambda}$ be an $f.r.f.t.$ nested  structure}. Then for any $\alpha\in\Lambda$, we have

\textit{(a).  $\# C_\alpha<\infty$, and $K_\alpha\setminus C_\alpha=\bigsqcup_{\beta\in \Lambda_\alpha} K_\beta\setminus C_\alpha$;}

\textit{(b). $K_{\beta}\cap K_{\beta'}=V_{\beta}\cap V_{\beta'}$, $\forall \beta, \beta'\in \Lambda_\alpha$;}

\textit{(c). For $\alpha\in \Lambda\setminus\{\vartheta\}$, $\#P(\alpha)=1$, i.e. the parent of $\alpha$ is unique;}

\textit{(d). $V_\alpha\cup C_\alpha\subset V_{\alpha,1}:=\bigcup_{\beta\in \Lambda_\alpha} V_\beta$.}

\textit{Proof.} (a). Clearly, for any $\beta\in \Lambda_\alpha$, $K_\beta\cap C_\alpha\subset V_\beta$ by the definition of $V_\beta$. Thus, $C_\alpha$ is a finite set as $C_\alpha=\bigcup_{\beta\in\Lambda_\alpha} C_\alpha\cap K_\beta\subset \bigcup_{\beta\in\Lambda_\alpha} V_\beta$. The equaility is obvious as $K_\beta\cap K_{\beta'}\subset C_\alpha, \forall \beta, \beta'\in \Lambda_\alpha$.

(b). Obviously, $K_\beta\cap K_{\beta'}\subset C_\alpha$, thus
\[K_\beta\cap K_{\beta'}=(K_\beta\cap C_\alpha)\cap (K_{\beta'}\cap C_\alpha)\subset V_\beta\cap V_{\beta'}\subset K_\beta\cap K_{\beta'}.\]
The equality follows immediately.

(c).  First we will prove by induction that $\#(K_\gamma\cap K_{\gamma'})<\infty$ for any $n\geq 1$ and any $\gamma\neq\gamma'\in \Lambda^{(n)}_\vartheta$. For $n=1$, the conclusion holds by (b). Now assume the conclusion holds for $n-1$. Let $\gamma, \gamma'$ be two distinct index in $\Lambda_\vartheta^{(n)}$. If $\gamma,\gamma'$ has a same parent, then the conclusion holds by (b). Otherwise, choose any $\zeta\in P(\gamma)\cap \Lambda^{(n-1)}_\vartheta,\zeta'\in P(\gamma')\cap \Lambda^{(n-1)}_\vartheta$, then $K_\gamma\cap K_{\gamma'}\subset K_{\zeta}\cap K_{\zeta'}$, which still implies that $\#(K_\gamma\cap K_{\gamma'})<\infty$ by the inductive assumption.

Suppose $\alpha\neq\vartheta$ and $\gamma,\gamma'$ are two distinct parents of $\alpha$. Assume $\gamma\in \Lambda^{(n)}_{\vartheta},\gamma'\in \Lambda^{(n')}_{\vartheta}$ with $n\geq n'$. Then for any $\zeta\in P^{n-n'}(\gamma)\cap\Lambda^{(n')}_{\vartheta}$($P^0(\gamma)=\{\gamma\}$), we have either  $\#(K_\zeta\cap K_{\gamma'})<\infty$ or $\zeta=\gamma'$. The former case is obviously impossible, and the latter case implies $K_\alpha=K_\alpha\cap K_\gamma\subset C_{\gamma'}$, a contradiction to (a). Thus $\#P(\alpha)=1$.

(d). We have already seen that $C_\alpha\subset V_{\alpha,1}:= \bigcup_{\beta\in \Lambda_\alpha} V_\beta$. Thus, we only need to show $V_\alpha\subset V_{\alpha,1}$. For any point $x\in V_\alpha$,  there exists an index $\beta\in \Lambda_\alpha$ such that $x\in K_\beta$. By the definition of $V_\alpha$, there exists $\alpha'\sim \alpha$ such that $x\in K_{\alpha}\cap \phi_{\alpha',\alpha}(C_{\gamma'})$ for some $\gamma'\in \bigcup_{n\geq 1}P^n(\alpha')$. Moreover, by A2, there exists $\beta'\in \Lambda_{\alpha'}$ such that $\beta'\sim \beta$ with $\phi_{\beta',\beta}=\phi_{\alpha',\alpha}$. Thus, 
\[x\in  K_{\beta}\cap \phi_{\alpha',\alpha}(C_{\gamma'})=\phi_{\beta',\beta}(K_{\beta'}\cap C_{\gamma'})\subset V_\beta.\]
So we get $V_\alpha\subset V_{\alpha,1}$, and thus $V_\alpha\cup C_\alpha\subset V_{\alpha,1}$. 
\hfill$\square$\\
 
\textbf{Remark on notations. }\textit{From Proposition 2.3(c), for a given $f.r.f.t.$  nested  structure $\{K_\alpha,\Lambda_\alpha\}_{{\alpha\in\Lambda}}$, we can simplify some notations. Since $\forall \alpha\in \Lambda\setminus\{\vartheta\}$, $P(\alpha)$ contains exactly one index, we can view $P(\cdot)$ as a mapping from $\Lambda\setminus\{\vartheta\}$ onto $\Lambda$. So the notations $\gamma=P(\alpha),K_{P(\alpha)},C_{P(\alpha)}$ make sense. In particular, the set $V_\alpha$ can be written as}
 \[V_\alpha=\bigcup_{\beta\sim\alpha,\beta\neq\vartheta}\phi_{\beta,\alpha}\Big(K_\beta\cap \big(\bigcup_{n\geq 1}C_{P^n(\beta)}\big)\Big).\]
 \textit{We will use $P^{-n}(\vartheta)$ instead of $\Lambda^{(n)}_{\vartheta}$, and the new notation can be applied to general indices, i.e., $P^{-n}(\alpha)$.}\\

In the following, we call $K_\alpha$ an \textit{island} with index $\alpha$, and $V_\alpha$ the \textit{boundary} of $K_\alpha$(sometimes we denote it by $\partial K_\alpha$). Write $V_0=V_\vartheta$, the boundary of $K$.   For $n\geq 1$, denote $V_n=\bigcup_{\alpha\in P^{-n}(\vartheta)}V_\alpha$, $V_*=\bigcup_{n=0}^{\infty} V_n$, and call points in $V_*$ \textit{vertices} in $K$. Note that $V_0$ can be  empty  if there is no $\alpha\sim \vartheta$. This can only happen for $V_\vartheta$ as $K_\alpha\cap C_{P(\alpha)}$ is not empty, $\forall\alpha\neq \vartheta$.
Obviously, the above notations agree with the ones for $K$ if $K$ is a $p.c.f.$ self-similar set and  $\{F_w K,\Lambda_w\}_{w\in {W_*}}$ is the canonical  nested structure.

For a general $f.r.f.t.$  nested structure $\{K_\alpha,\Lambda_\alpha\}_{\alpha\in \Lambda}$, we have

\textbf{Proposition 2.4.} \textit{  For any $n\geq 0$, $\# V_n<\infty$, and 
	 $V_n\subset V_{n+1}$}. Moreover, $V_*$ is dense in $K$.

\textit{Proof.}  For $n\geq 0$, obviously we have $\# V_n<\infty$. $V_n\subset V_{n+1}$ is an easy consequence of Proposition 2.3(d). Hence it suffices to show that 
\[
\lim_{n\to\infty} \sup\{diam(K_\alpha)|\alpha\in P^{-n}(\vartheta)\}=0,
\]
where $diam(A)=\sup_{x,y\in A}d(x,y)$ is the \textit{diameter} of the set $A$. Denote $\lambda=\sup\{diam(K_\beta)/ diam(K_\alpha): \alpha\neq\beta\in \Lambda, \text{ with }\beta\sim\alpha \text{ and }\alpha=P^m(\beta)\text{ for some }1\leq m\leq M\}$. Obviously $\lambda<1$ as the supremum is taken over finite cases, since there are only finite number of equivalent classes. It is  easy to see that for any $n\geq M$, and any $\alpha\in P^{-n}(\vartheta)$,
there always exist $\beta,\beta'\in\{\alpha,P(\alpha),\cdots,P^M(\alpha)\}$ such that $\beta\sim \beta'$ and $\beta'=P^m(\beta)$ for some $1\leq m\leq M$, which results that
\[diam(K_\alpha)/diam(K_{P^{M}(\alpha)})\leq diam(K_\beta)/diam(K_{\beta'})\leq \lambda.\]
Thus, if $n$ is sufficient large, then for any $\alpha\in P^{-n}(\vartheta)$,
\[diam(K_\alpha)\leq \lambda\cdot diam(K_{P^M(\alpha)})\leq \cdots\leq \lambda^{[\frac{n}{M}]}\cdot diam(K),\]
where the right side of the inequailty goes to $0$ obviously, as $n$ goes to infinity. This gives that $V_*$ is dense in $K$. \hfill$\square$

\textbf{Remark.} \textit{ One may weaken the assumptions in Definition 2.2 by requiring  $\phi_{\alpha,\beta}$ to be just a homeomorphism instead of similitude.   Then Proposition 2.3 still holds and the notations $V_n, V_*$ still make sense. However, the weakened assumptions could not ensure that $V_*$ is dense in $K$. The following is a counterexample. }

Consider the line segment $I=[0,1]$, which can be viewed as a self-similar set with the $i.f.s.$, $F_1: x\rightarrow \frac{1}{2}x$, $F_2: x\rightarrow\frac{1}{2}x+\frac 12.$ 
We introduce another two piecewise linear mappings
\[\tilde{F}_1(x)=\begin{cases}x&\text{  if }x \in [0,\frac{1}{4}],\\\frac{1}{3}x+\frac{1}{6}&\text{  if }x \in [\frac{1}{4},1],\end{cases}\text{ and }\quad
\tilde{F}_2(x)=\begin{cases}x+\frac{1}{2}&\text{  if }x \in [0,\frac{1}{4}],\\\frac{1}{3}x+\frac{2}{3}&\text{  if }x \in [\frac{1}{4},1].\end{cases}
\]
Choose $\Lambda$ to be $W_*=\bigcup_{n\geq 0}\{1,2\}^n$. Denote $I_\vartheta=I$ and for $w=w_1w_2\cdots w_n\in W_*$ with $n\geq 1$, denote
\[I_w=\tilde{F}_{w}I \text{ with }\tilde{F}_{w}=\tilde{F}_{w_1}\circ\tilde{F}_{w_2}\circ\cdots \circ\tilde{F}_{w_n}.\]
 Obviously, $\{I_w,\Lambda_w\}_{w\in W_*}$ is a nested structure of $1$ type,  with $\Lambda_w$ being the same notation as that for $p.c.f.$ self-similar sets, and $\phi_{w,w'}=\tilde{F}_{w'}\tilde{F}^{-1}_w$. Noticing that $(0,\frac{1}{4})\subset \tilde{F}^n_{1}I$ for any $n\geq 1$, $V_*$ is not dense in $I$ since $V_*\cap (0,\frac{1}{4})=\emptyset$.

Throughout the following text, for an $f.r.f.t.$  nested structure $\{K_\alpha,\Lambda_\alpha\}_{\alpha\in{\Lambda}}$, we use $\{\mathcal{T}_1,\mathcal{T}_2,\cdots,\mathcal{T}_M\}$ to denote the equivalent classes in $\Lambda$ according to ``$\sim$'', called the \textit{types} of islands. For $\alpha\in\Lambda$, we denote $t(\alpha)$ the integer in $\{1,2,\cdots,M\}$ such that $\alpha\in \mathcal{T}_{t(\alpha)}$, call $\mathcal{T}_{t(\alpha)}$ the type of $\alpha$.  Without loss of generality, we always assume that the type of $\vartheta$ is $\mathcal{T}_1$. Islands with the same type only differ by a similitude. 

\section{Relationship with $p.c.f.$ structures and $f.r.g.d.$ structures}
In this section, we will discuss the relationship of $f.r.f.t.$ self-similar sets with $p.c.f.$ self-similar sets and $f.r.g.d.$ fractals. We will show that the concept of $f.r.f.t.$ structures is an intermediate one, i.e., a connected $p.c.f.$ self-similar set has a canonical $f.r.f.t.$  nested structure, and an $f.r.f.t.$  nested structure implies a finitely ramified graph-directed construction.  We will also provide some other sufficient conditions for a connected self-similar set to possess a $f.r.f.t.$ nested structure.

First, let's look at the relationship with $p.c.f.$ self-similar sets. 

We briefly recall the definition of $p.c.f.$ self-similar structures. Let $K$ be a self-similar set associated with an $i.f.s.$ $\{F_i\}_{i=1}^N$. Denote $W_n$, $W_*$ as before, and write $\Sigma=\{1,2,\cdots,N\}^{\mathbb{N}}$ the one-sided \textit{shift space} with the \textit{shift map} $\sigma: \Sigma\rightarrow\Sigma$ defined as $\sigma(\omega_1\omega_2\cdots)=\omega_2\omega_3\cdots$.  Let $\pi:\Sigma\to K$ be the natural projection defined as
\[\pi(\omega)=\bigcap_{n=1}^\infty F_{[\omega]_n}(K),\quad\omega\in \Sigma,\]
where $[\omega]_n=\omega_1\omega_2\cdots\omega_n\in W_n$ for $\omega=\omega_1\omega_2\cdots$. Let $C_K=\bigcup_{1\leq i<j\leq N}F_iK\cap F_jK$ be the  union of intersections of cells of level $1$. Define the \textit{critical set} to be $\mathcal{C}=\pi^{-1}(C_K)$ and the \textit{post-critical set} to be $\mathcal{P}=\bigcup_{n\geq 1}\sigma^n(\mathcal{C})$. Call the pair $(K,\{F_i\}_{i=1}^N)$ a post-critically finite($p.c.f.$ for short) structure if $\mathcal{P}$ is a finite set. In addition, call such $K$ a $p.c.f.$ self-similar set.

\textbf{Theorem 3.1.} \textit{Let $K$ be a connected p.c.f. self-similar set, then  the canonical nested structure $\{F_wK,\Lambda_w\}_{w\in W_*}$ is $f.r.f.t.$ with $M=1$.}

\textit{Proof.} Obviously,  $W_*=W_\#$, where $W_\#$ is the same notation as introduced in Section 2.
It suffices to verify A1,A2,A3 in Definition 2.2. For the sake of uniformity, in the remaining proof, we write $\Lambda=W_*$, and  for any $\alpha\in \Lambda$, denote $K_\alpha=F_\alpha K$, $\Lambda_\alpha=\{\alpha i| i\in\{1,2,\cdots,N\}\}$ and $C_\alpha=\bigcup_{i\neq j}K_{\alpha i}\cap K_{\alpha j}$.

Firstly, for any $\alpha,\beta\in \Lambda$, choose the natural similitude $\phi_{\alpha,\beta}=F_\beta\circ F_\alpha^{-1}$ such that $K_\beta=\phi_{\alpha,\beta}(K_\alpha)$. This induces a trivial equivalent relation ``$\sim$'' in $\Lambda$ that every two indices are equivalent. So A1 is true, with $\Lambda/\sim=\{\mathcal{T}_1\}$. 

Secondly, for any $\alpha,\beta\in\Lambda$ and $1\leq i\leq N$, it is direct to check that $\phi_{\alpha i,\beta i}=F_\beta\circ F_i \circ(F_\alpha\circ F_i)^{-1}=F_\beta\circ F_\alpha^{-1}=\phi_{\alpha,\beta}$, which yields a one to one correspondence between $\Lambda_\alpha$ and $\Lambda_\beta$. This gives A2.

Lastly, notice that for any $\alpha\in\Lambda$ we have $C_\alpha=F_\alpha C_K$, thus
\[\begin{aligned}
V_\vartheta&=\bigcup_{\alpha\in \Lambda\setminus\{\vartheta\}}\phi_{\alpha,\vartheta}\big(K_\alpha\cap(\bigcup_{n\geq 1}\bigcup_{\gamma\in P^n(\alpha)}C_\gamma)\big)=\bigcup_{\alpha\in W_*\setminus\{\vartheta\}} F_\alpha^{-1}(K_\alpha\cap \bigcup_{0\leq n<|\alpha|}C_{[\alpha]_n})\\
&=\bigcup_{\alpha\in W_*\setminus\{\vartheta\}}K\cap F_\alpha^{-1}(\bigcup_{0\leq n<|\alpha|} C_{[\alpha]_n})=\bigcup_{\alpha\in W_*\setminus\{\vartheta\}}\bigcup_{0\leq n<|\alpha|} K\cap F_\alpha^{-1}(C_{[\alpha]_n})\\
&=\bigcup_{\beta\in W_*}\bigcup_{\gamma\in W_*\setminus \{\vartheta\}} K\cap F_{\beta\gamma}^{-1}(C_\beta)=\bigcup_{\beta\in W_*}\bigcup_{\gamma\in W_*\setminus \{\vartheta\}}K\cap F_\gamma^{-1}(C_K)\\
&=K\cap \big(\bigcup_{\gamma\in W_*\setminus \{\vartheta\}}F_\gamma^{-1}(C_K)\big)=\pi (\mathcal{P}).
\end{aligned}\]
So $V_\vartheta$ is finite as $\mathcal{P}$ is finite, A3 is satisfied.\hfill$\square$\\

Next, we discuss the relationship between $f.r.f.t.$ self-similar sets and $f.r.g.d.$ fractals. Meanwhile, we show how to reconstruct a self-similar set in a  finitely ramified graph-directed manner by virtue of its $f.r.f.t.$ nested structure. Some notations will be frequently used in the following sections.

Recall the concepts of graph-directed construction and $f.r.g.d.$ fractals, which can be found in detail in [HN]. Let $G=(S,E)$ be a \textit{directed graph}, where $S$ is the set of \textit{states}(call vertices in this graph states to avoid confusion) and $E$ is the set of \textit{edges} of the graph. Note that multiple edges and loops are allowed. For an edge $e\in E$, denote by $i(e)$ the \textit{initial state} of $e$, and $f(e)$ the \textit{final state} of $e$. 

\textbf{Definition 3.2.} \textit{ Let $G=(S,E)$ be a directed graph. Assign each $e\in E$ a similitude $\psi_e$ with similarity ratio $l_e$, and each $s\in S$ a compact connected set $J_s$. Call $\mathcal{G}=(S,E,\{\psi_e\}_{e\in E})$ a} graph-directed construction  \textit{if the following conditions are satisfied, }

\textit{1. $\forall s\in S$, there is at least one edge $e\in E$, such that $s=i(e)$;}

\textit{2. $\forall s\in S$, $\bigcup_{i(e)=s}\psi_eJ_{f(e)}\subset J_s$;}

\textit{3. For a cycle $e_1e_2\cdots e_n$, where cycle means $f(e_k)=i(e_{k+1}), \forall k=1,2,\cdots, n-1$ and $f(e_n)=i(e_1)$, we have $\prod_{k=1}^n l_{e_k}<1$. }

It is well-known that there is a unique vector of compact sets $\mathcal{K}=\{K_s\}_{s\in S}$ such that $K_s=\bigcup_{i(e)=s}\psi_e K_{f(e)}, \forall s\in S$. We call them the \textit{invariant sets} of the graph-directed construction $\mathcal{G}$.

We define a \textit{shift space} associated with $\mathcal{G}$ to address points in $K_s,s\in S$. A finite sequence of edges in $G$, denoted by $\bm{e}=e_1e_2\cdots e_n$, is called a \textit{walk} if $f(e_k)=i(e_{k+1}),\forall 1\leq k\leq n-1$. We write $|\bm{e}|=n$ for the length of the walk. An infinite sequence of edges is called an \textit{infinite walk}, denoted by $\bm{\epsilon}=\epsilon_1\epsilon_2\cdots$, if for any $n\geq 1$, the first $n$ steps $[\bm{\epsilon}]_n=\epsilon_1\epsilon_2\cdots\epsilon_n$ is a walk of length $n$. Denote by $E_*$ the collection of finite walks in $G$, and $E^\infty$ the space of infinite walks. For convenience, let $i(\bm{e})=i(e_1)$ or $i(\bm{\epsilon})=i(\epsilon_1)$ the initial state of a walk, and let $f(\bm{e})=f(\bm{e}_{|e|})$ the final state of a walk. Then we define a projection $\pi:E^\infty\to \bigcup_{s\in S}{K_s}$ by 
\[\pi(\bm{\epsilon})=\bigcap_{n=1}^\infty \psi_{[\bm{\epsilon}]_n}K_{f(\bm{[\epsilon]_n})},\]
where we use the notation $\psi_{\bm{e}}=\psi_{e_1}\circ\psi_{e_2}\circ\cdots\circ\psi_{e_n}$.

Analogous to  $p.c.f.$ self-similar sets, for each $s\in S$, we introduce the set of level $1$ intersection $C_s=\bigcup_{e\neq e'\in i^{-1}(s)}\psi_eK_{f(e)}\cap \psi_{e'}K_{f(e')}$, where $e\in i^{-1}(s)$ means $i(e)=s$, the \textit{critical set} $\mathcal{C}_{\mathcal{G}}=\bigcup_{s\in S}\pi^{-1}(C_s)$, and the \textit{post-critical set} $\mathcal{P}_{\mathcal{G}}=\bigcup_{n=1}^\infty \sigma^n (\mathcal{C}_{\mathcal{G}})$, where $\sigma$ is the \textit{shift map} on $E^\infty$, i.e., $\sigma(\epsilon_1\epsilon_2\cdots)=\epsilon_2\epsilon_3\cdots$. Then for each $s\in S$, we write $V_s=\{\pi(\bm{\epsilon}):i(\bm{\epsilon})=s,\bm{\epsilon}\in\mathcal{P}_{\mathcal{G}}\}$.

\textbf{Definition 3.3.} \textit{A family  $\mathcal{K}=\{K_s\}_{s\in S}$ constructed by the graph-directed construction $\mathcal{G}=(S,E,\{\psi_e\}_{e\in E})$ is called a} finitely ramified graph-directed\textit{($f.r.g.d.$  for short)} fractal family \textit{ if $V_s$ is finite  for each $s\in S$. Each member $K_s\in \mathcal{K}$ is called an }$f.r.g.d.$ fractal. \\

Now we proceed to clarify the relationship between $f.r.f.t.$ self-similar sets and $f.r.g.d.$ fractals.

Let us start with an $f.r.f.t.$ nested structure $\{K_\alpha,\Lambda_\alpha\}_{\alpha\in \Lambda}$, with $\Lambda/\sim=\{\mathcal{T}_1,\mathcal{T}_2,\cdots,\mathcal{T}_M\}$.  For $1\leq i\leq M$, choose an element $\alpha_i$ in  $\mathcal{T}_i$(for convenience, we always require $K_{\alpha_i}$ has the largest diameter in islands of type $\mathcal{T}_i$, obviously $\alpha_1=\vartheta$), 
then we have
\begin{equation}\label{eq23}
\begin{aligned}
K_{\alpha_i}=\bigcup_{\beta\in\Lambda_{\alpha_i}} K_\beta=\bigcup_{\beta\in \Lambda_{\alpha_i}}\phi_{\alpha_{t(\beta)},\beta}(K_{\alpha_{t(\beta)}}),\quad \forall 1\leq i\leq M.
\end{aligned}
\end{equation}
Thus we can construct a directed graph $G=(S,E)$ with the \textit{state} set $S$   and the \textit{edge} set $E$  defined as
\begin{equation}
S=\{\mathcal{T}_i\}_{i=1}^M,\quad E=\{(\mathcal{T}_i, \mathcal{T}_{t(\beta)}): 1\leq i\leq M,\beta\in \Lambda_{\alpha_i}\}.
\end{equation}
Obviously, $\mathcal{G}=(S,E,\{\phi_{\alpha_{t(\beta)},\beta}\})$ is a graph-directed construction, and $K$ is an $f.r.g.d.$ fractal as a member of the $f.r.g.d.$ family $\{K_{\alpha_i}\}_{i=1}^M$. 

\textbf{Theorem 3.4.} \textit{ A connected self-similar set possesses an $f.r.f.t.$ nested structure if and only if it is an $f.r.g.d.$ fractal. }

\textit{Proof.} The ``only if'' part is true as discussed before. 

Let's look at the ``if'' part. Suppose $K$ is a connected $f.r.g.d.$ fractal, then it must be a member of an $f.r.g.d.$ fractal family $\mathcal{K}=\{K_s\}_{s\in S}$ with an associated graph-directed construction $\mathcal{G}=(S,E,\{\psi_e\}_{e\in E})$. Let $s^*$ be the associated state of $K$, and $S'$ be the collection of all states that appear in the walks emanating from $s^*$. Without loss of generality, assume $S'=S$, otherwise we only need to consider the subgraph whose states set is $S'$ instead. 

Define $\Lambda=\{\bm{e}\in E_*|i(\bm{e})=s^*\}\cup\{\vartheta\}$, write $K_\vartheta=K$ and  for any $\bm{e}\in \Lambda\setminus\{\vartheta\}$, write $K_{\bm{e}}=\psi_{\bm{e}}K_{f(\bm{e})}$. We claim that  $\{K_\alpha,\Lambda_\alpha\}_{\alpha\in\Lambda}$ is an $f.r.f.t.$ nested structure of $K$. 

In fact, it is easy to see that  $\{K_\alpha,\Lambda_\alpha\}_{\alpha\in\Lambda}$ is a nested  structure with $\Lambda_{\bm{e}}=\{\bm{e} e:e\in E,i(e)=f(\bm{e})\}$ for any $\bm{e}\in\Lambda\setminus\{\vartheta\}$, and trivially $\Lambda_{\vartheta}=\{e\in E:i(e)=s^*\}$.

We only need to verify the assumptions A1,A2,A3. We regard $\vartheta$ as an empty walk and let $f(\vartheta)=s^*$, $\psi_\vartheta$ be the identity map for consistency.  For A1, there are $M=\# S$ types of islands, and $\bm{e}\sim\bm{e}'$ if and only if $f(\bm{e})=f(\bm{e}')$ with $\phi_{\bm{e},\bm{e}'}=\psi_{\bm{e}'}\circ\psi_{\bm{e}}^{-1}$. For A2, it is easy to check that $\phi_{\bm{e}e,\bm{e}'e}=\phi_{\bm{e},\bm{e}'}$ for any $e\in E$ with $i(e)=f(\bm{e})=f(\bm{e}')$, and this provides a one to one correspondence between $\Lambda_{\bm{e}}$ and $\Lambda_{\bm{e}'}$. As for A3, for any $\bm{e}\in\Lambda$ we have
\[\begin{aligned}
V_{\bm{e}}:&=\bigcup_{\bm{e}'\sim\bm{e},\bm{e}'\neq \vartheta} \phi_{\bm{e}',\bm{e}}\big(K_{\bm{e}'}\cap (\bigcup_{n\geq 1}\bigcup_{\bm{e}''\in P^n(\bm{e}')}\psi_{\bm{e}''}C_{f(\bm{e}'')})\big)\\
&=\bigcup_{\bm{e}'\sim\bm{e},\bm{e}'\neq \vartheta}K_{\bm{e}}\cap \phi_{\bm{e}',\bm{e}}(\bigcup_{0\leq n< |\bm{e}'|} \psi_{[\bm{e}']_n}C_{f([\bm{e}']_n)})\\&=\psi_{\bm{e}}\Big(K_{f(\bm{e})}\cap \big(\bigcup_{\bm{e}'\sim\bm{e},\bm{e}'\neq\vartheta}\psi_{\bm{e}'}^{-1}(\bigcup_{0\leq n<|\bm{e}'|} \psi_{[\bm{e}']_n}{C}_{f([\bm{e}']_n)})\big)\Big)\\
&=\psi_{\bm{e}} (K_{f(\bm{e})}\cap \bigcup_{\bm{e}'\sim \bm{e},\bm{e}'\neq\vartheta} \psi_{\bm{e}'}^{-1}C_{i(\bm{e}')})=\psi_{\bm{e}}(V_{f(\bm{e})}),
\end{aligned}\] 
which is finite since $V_{f(\bm{e})}$ is finite.

Thus we have proved that   $\{K_\alpha,\Lambda_\alpha\}_{\alpha\in\Lambda}$ is an $f.r.f.t.$ nested structure of $K$. 
\hfill$\square$

In the remaining part of this section, we will show some other sufficient conditions for the existence of an $f.r.f.t.$ nested structure of a given connected self-similar set.  For a self-similar set $K$, we write $\{F_i\}_{i=1}^N$, $W_*, W_\#$ as before. We denote $c_i$ the \textit{similarity ratio} of $F_i$, $1\leq i\leq N$, and   $c_*=\min\{c_1,c_2,\cdots,c_N\}$. For $w=w_1 w_2\cdots w_n\in W_*$,  write $c_w=c_{w_1} c_{w_2}\cdots c_{w_n}$ the similarity ratio of $F_w$ for short. We will introduce some conditions concerning overlaps.

 \textbf{Definition 3.5.} \textit{A self-similar set $K$ is said to satisfy the} finite neighboring type property \textit{ if there are only finitely many similitudes $h=F_w^{-1}F_u$ with $w,u\in W_*$ and $F_wK\cap F_uK\neq \emptyset$, and with similarity ratio $c_h\in(c_*, 1/c_*)$.}

 This condition, formulated in algebraic terms, was introduced in [B] and [BR] by Bandt and Rao to describe algorithms to verify the open set condition. It is also related with the finite type concept in  [LN, NW] by Lau and Ngai, Ngai and Wang to determine the Hausdorff dimension of certain self-similar sets with overlaps. In fact, it is a more restrictive but simpler version than that in [LN, NW].

 We call a finite sequence of distinct copies $\gamma=(F_{w^{(1)}}K, F_{w^{(2)}}K,\cdots, F_{w^{(n)}}K)$ with $w^{(i)}\in W_*$, $\forall 1\leq i\leq n$, an \textit{overlapping chain} if $\#\big(F_{w^{(i+1)}}K\cap (\bigcup_{1\leq j\leq i}F_{w^{(j)}}K)\big)=\infty, \forall 1\leq i\leq n-1$, and $F_{w^{(i)}}K\nsubseteq F_{w^{(j)}}K$ for distinct $1\leq i,j\leq n$, and call $n$ the \textit{length} of $\gamma$. Moreover, for $0<\delta<1$, we call the chain $\gamma$ a \textit{$\delta$-overlapping chain} if $\delta\leq c_{w^{(i)}}c_{w^{(j)}}^{-1}\leq 1/\delta$ for any $F_{w^{(i)}}K$ and $F_{w^{(j)}}K$ contained in $\gamma$. Denote $\mathcal{L}(K)$ the supremum of the lengths of overlapping chains in $K$, and  $\mathcal{L}_\delta(K)$ the supremum of the lengths of  $\delta$-overlapping chains in $K$, $\forall 0<\delta <1$.

 \textbf{Definition 3.6.} \textit{A self-similar set $K$ is said to satisfy the} finite chain length property \textit{ if $\mathcal{L}(K)<\infty$, and  the} finite $\delta$-chain length property \textit{ if  $\mathcal{L}_\delta(K)<\infty$ for some $0< \delta <1$.}

 Obviously, $\mathcal{L}(K)<\infty$ implies $\mathcal{L}_\delta(K)<\infty$ for any $0<\delta<1$. And for $0<\delta_1\leq \delta_2<1$, $\mathcal{L}_{\delta_1}(K)<\infty$ implies $\mathcal{L}_{\delta_2}(K)<\infty$. Nevertheless, we have

 \textbf{Proposition 3.7.} \textit{Let $K$ be a self-similar set. Then the following two conditions are equivalent:}

 \textit{1. there exists  $0<\delta\leq c_*$ such that $\mathcal{L}_\delta(K)<\infty$;}

 \textit{2. for any $0<\delta\leq c_*$, $\mathcal{L}_\delta(K)<\infty$.}

 \textit{Proof.} We only need to prove ``$1\Rightarrow 2$''.

 Fix a $0<\delta\leq c_*$ such that $\mathcal{L}_\delta(K)<\infty$. We just need to show $\mathcal{L}_{\delta'}(K)<\infty$ for any $0<\delta'\leq \delta$. Let  $\gamma$ be a $\delta'$-overlapping chain with  $\gamma=(F_{w^{(1)}}K, F_{w^{(2)}}K,\cdots, F_{w^{(n)}}K)$. We denote by $\lambda=\max\{c_{w^{(i)}}: 1\leq i\leq n\}$, and write $W_\lambda=\{w\in W_\#: \lambda c_*< c_w\leq \lambda\}$. Then for each $1\leq i\leq n$, we could choose $\widetilde{w}^{(i)}$ in $W_\lambda$ such that $F_{w^{(i)}}K\subset F_{\widetilde{w}^{(i)}}K.$ It is clear that, after deleting the repeated ones if necessary, $(F_{\widetilde{w}^{(1)}}K, F_{\widetilde{w}^{(2)}}K,\cdots, F_{\widetilde{w}^{(n)}}K)$ forms a $c_*$-overlapping chain, we denote it by $\widetilde{\gamma}$. Since $\delta\leq c_*$, $\widetilde{\gamma}$ is also a $\delta$-overlapping chain, and hence the length of $\widetilde{\gamma}$ is no more than $\mathcal{L}_\delta(K)$. Notice that the similarity ratio of elements in $\gamma$ have a lower bound $\lambda \delta'$, and the similarity ratio of elements in $\widetilde{\gamma}$ have an upper bound $\lambda$. Then an easy calculation shows that the length of $\gamma$ is no more than $\# W_{\delta'}\mathcal{L}_\delta(K)$. From the arbitrariness of $\gamma$, we have proved that $\mathcal{L}_{\delta'}(K)<\infty$. \hfill$\square$

 Then we have

 \textbf{Theorem 3.8.} \textit{Let $K$ be a connected self-similar set, satisfying the finite neighboring type property and the finite $\delta$-chain length property for some $\delta\leq c_*$. Then $K$ possesses an f.r.f.t. nested structure.} 

 \textit{Proof.} By Proposition 3.7, without loss of generality, we just take $\delta=c_*$. For $\lambda>0$, as in the proof of Proposition 3.7, we write $W_\lambda=\{w\in W_\#: \lambda c_*< c_w\leq \lambda\}$. Notice that 
$
K=\bigcup_{w\in W_\lambda} F_wK.
$

  By choosing $\lambda_0$ sufficiently small, we could assume $\# W_{\lambda_0} >\mathcal{L}_{c_*}(K)$. 
Then rewrite the decomposition $K=\bigcup_{w\in W_{\lambda_0}} F_wK$ into the following finite union of connected compact subsets,
\begin{equation}\label{eq311}
K=\bigcup\big\{K_\gamma: \gamma \text{ is a }  c_*\text{-overlapping chain with elements in} \{F_wK: w\in W_{\lambda_0}\}\big \},
\end{equation}
with $K_\gamma:=\bigcup_{w\in \gamma} F_wK$, and with $\#(K_\gamma\cap K_{\gamma'})<\infty$ for distinct $\gamma,\gamma'$. Since $\# W_{\lambda_0} >\mathcal{L}_{c_*}(K)$, the union in (\ref{eq311}) contains at least two distinct $K_\gamma$'s. 

 Next, by choosing $\lambda_1$ sufficiently small, for each $K_\gamma$ in (\ref{eq311}), we could assume that $K_\gamma=\bigcup\{F_wK: w\in W_{\lambda_1}, F_wK\subset K_\gamma\}$ with $\# \{w: w\in W_{\lambda_1}, F_wK\subset K_\gamma\} >\mathcal{L}_{c_*}(K)$.
Similarly as above, we then redecompose  each $K_\gamma$ into a finite (at least two) union of connected compact subsets which intersect each other only at finite points,
\begin{equation}\label{eq312}
K_\gamma=\bigcup\big\{K_{\gamma'}: \gamma' \text{ is a }  c_*\text{-overlapping chain with elements in} \{F_wK: w\in W_{\lambda_1}\}\big \}.
\end{equation}

 By using the finite neighboring type property, and noticing that $\mathcal{L}_{c_*}(K)<\infty$, we obtain that there are only finitely many types of $c_*$-overlapping chains  involved in (\ref{eq311}) and (\ref{eq312}) up to similitudes. So (\ref{eq311}) and (\ref{eq312}) provide an $f.r.g.d.$ graph-directed reconstruction of the self-similar set $K$. Then by using Theorem 3.4, $K$ possesses an $f.r.f.t.$ nested structure. Thus we have proved the theorem. \hfill$\square$

 Obviously, we have

 \textbf{Corollary 3.9.} \textit{Let $K$ be a connected self-similar set, satisfying both the finite neighboring type property and the finite chain length property. Then $K$ possesses an f.r.f.t. nested structure.} 

 \textbf{Remark.} The finite neighboring type property and the finite chain length property can not imply each other. See counterexamples in Appendix.

\section{Examples of  $f.r.f.t.$ nested  strutures}
In this section, we will show some examples possessing $f.r.f.t.$  nested structures.

\textbf{Example 1.(Vicsek set with overlaps)} 
Let $\{q_i\}_{i=1}^4$ be the four vertices of a square, and $q_5$ be the center of the square. The \textit{Vicsek set with overlaps}, denoted by $\mathcal{V}^o$, is the invariant set of the $i.f.s.$ $\{F_i\}_{i=1}^5$,
\[F_1: x\rightarrow\frac{1}{2} x+\frac{1}{2}q_1,\quad F_i: x\rightarrow\frac{1}{3} x+\frac{2}{3}q_i,i=2,3,4,5,\]
see Figure 1.2(left) for $\mathcal{V}^o$. There are two types of connected compact subsets in $\mathcal{V}^o$, including $\mathcal{V}^o$, among which each can be split into similar copies by removing finite number of  appropriate points. See Figure 4.1 for an illustration, where the left one is $\mathcal{V}^o$ and the right one is $\bigcup_{i=1,2,4,5}F_i\mathcal{V}^o$. We use different colors to indicate the types of the pieces, and dot the points of overlaps. 
 \begin{figure}[h]
	\centering
	\includegraphics[height=3cm]{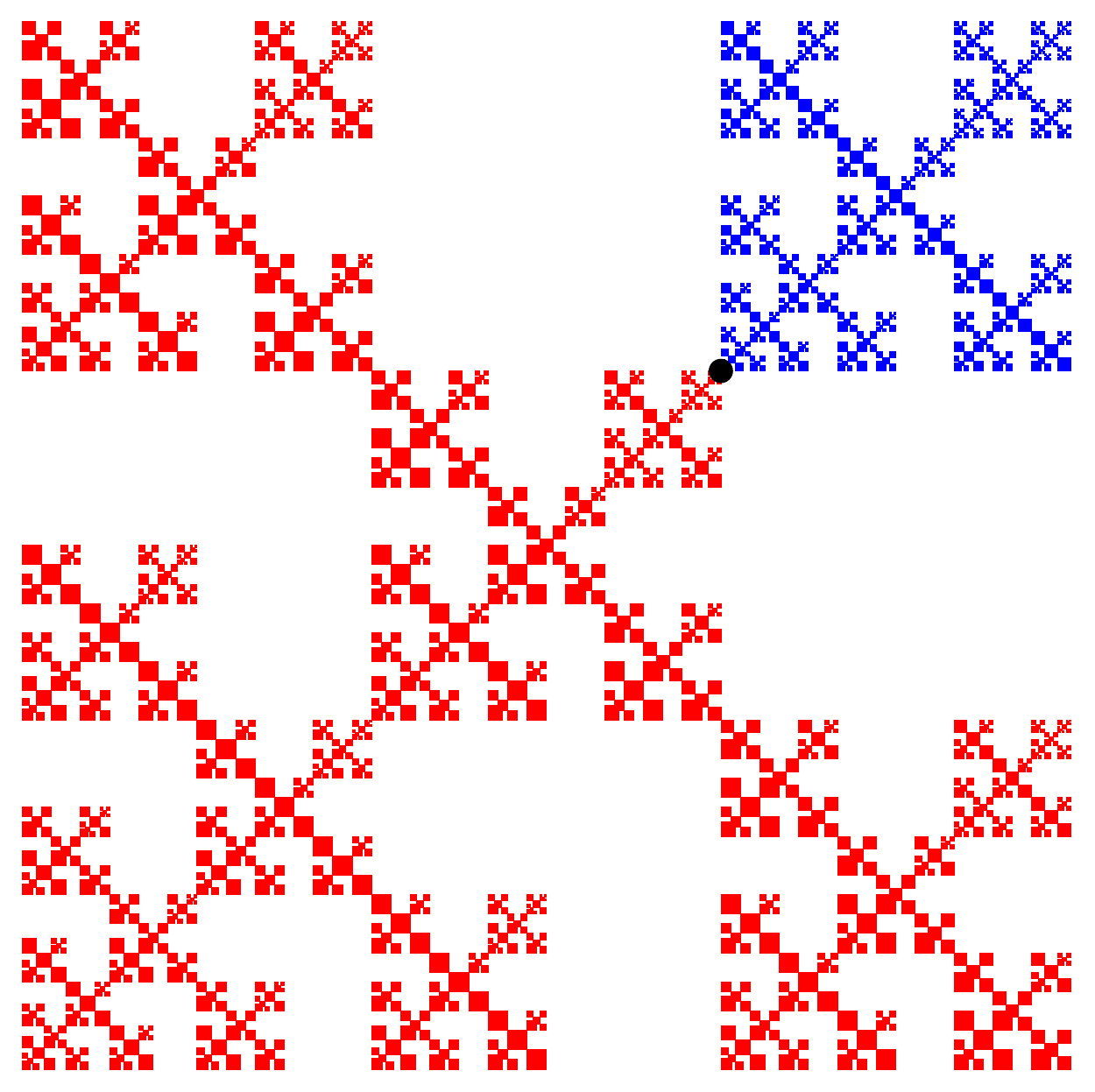}\qquad
	\includegraphics[height=3cm]{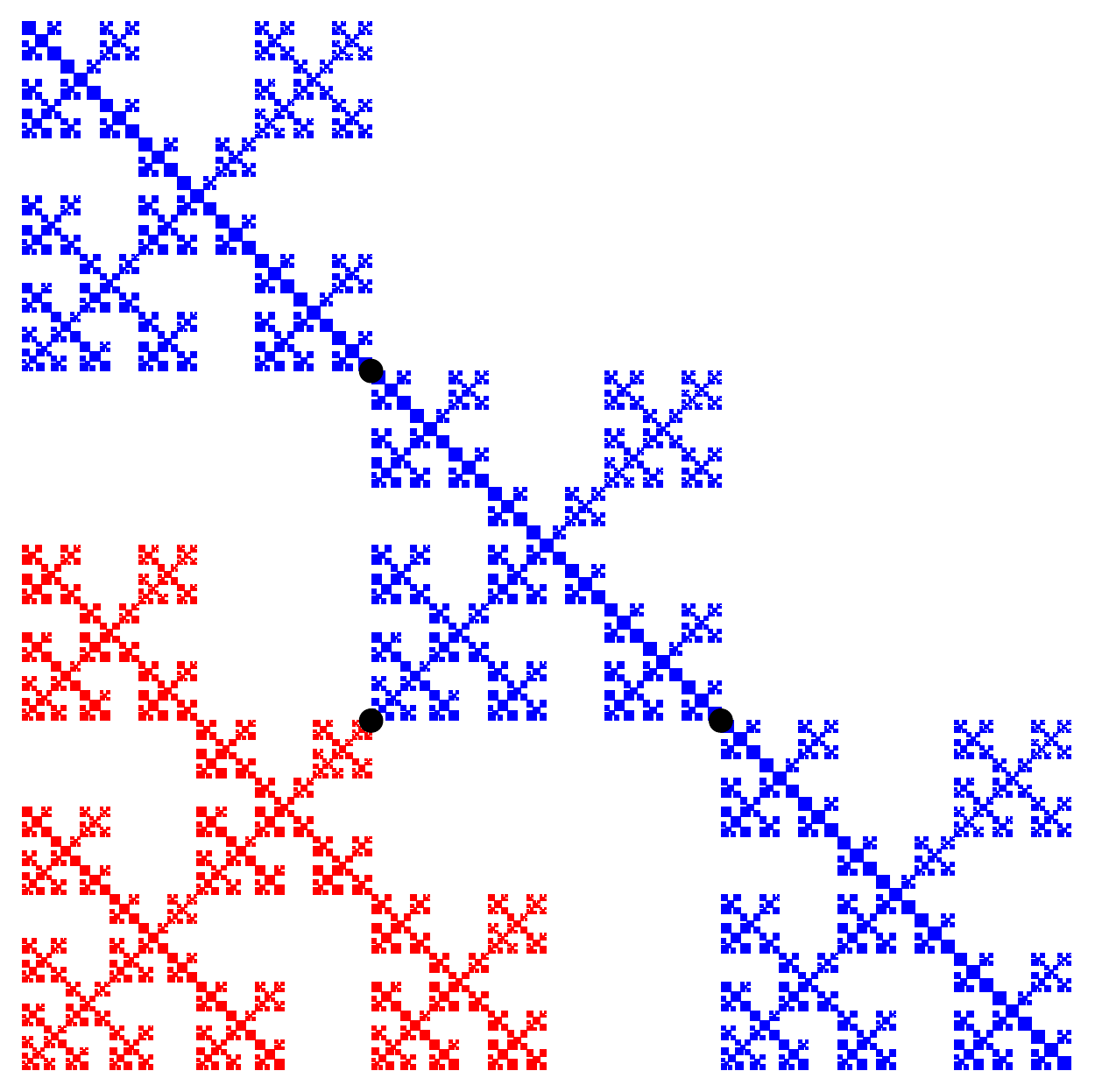}
		\begin{picture}(0,0) \thicklines 
	\put(-202,-0.5){$q_1$}
	\put(-109,-0.5){$q_2$}
	\put(-202,84){$q_4$}
	\put(-109,84){$q_3$}
	\end{picture}
	\begin{center}
		\caption{Two typical compact subsets in $\mathcal{V}^o$.}
	\end{center}
\end{figure}
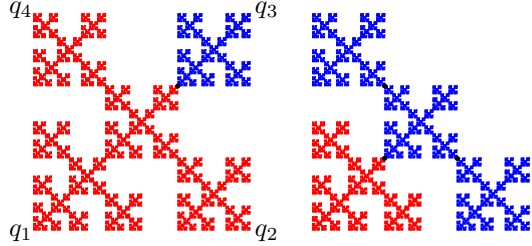

We can further divide $\mathcal{V}^o$ into smaller pieces by inductively using the ``cutting rule'' shown in Figure 4.1, see Figure 4.2 for the level-3 division. 

\begin{figure}[h]
 	\centering
 	\includegraphics[height=3.5cm]{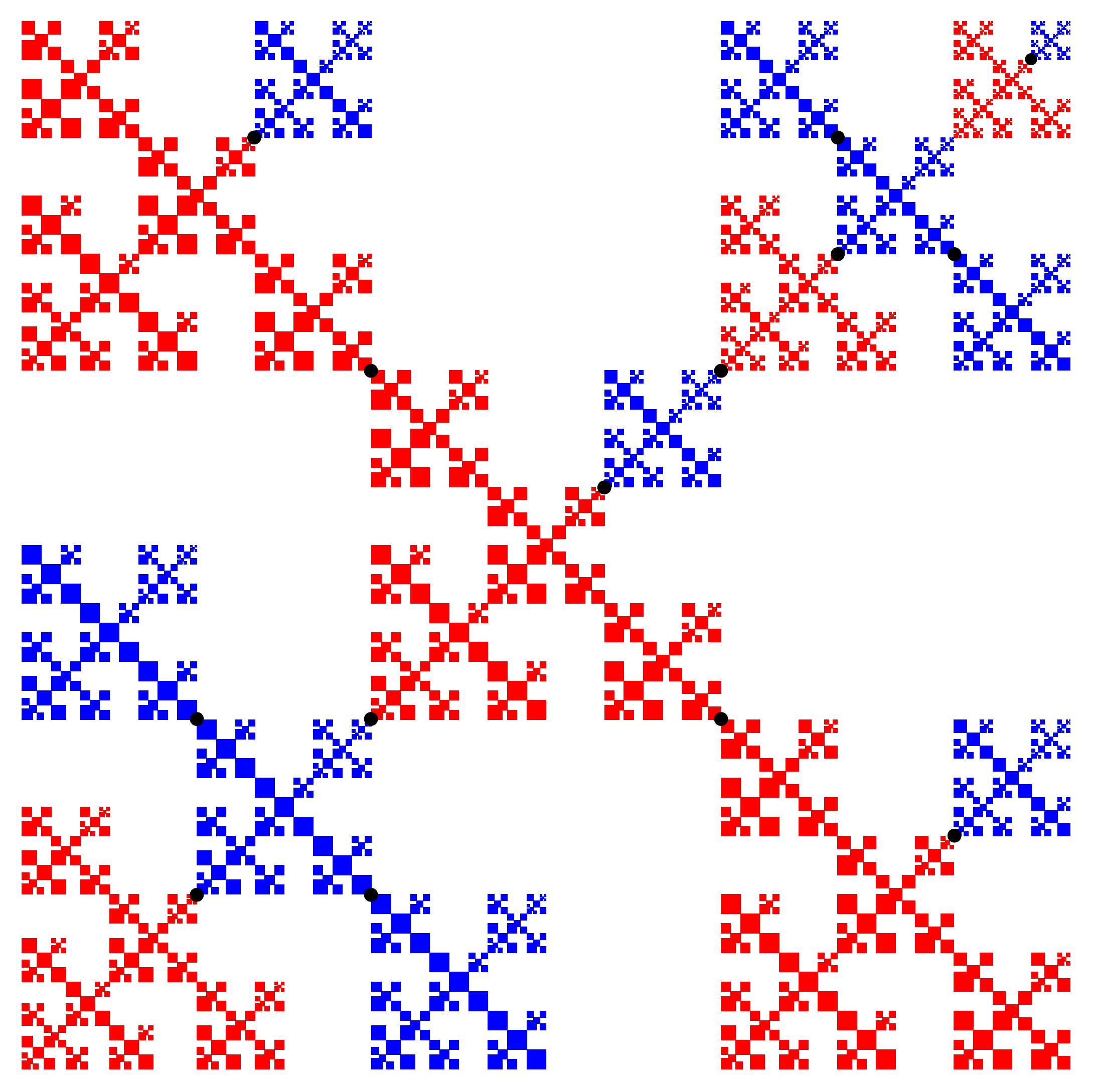}
 	\begin{center}
 		\caption{The level-$3$ division  of $\mathcal{V}^o$.}
 	\end{center}
\end{figure}

To be more precise, let  
\[K_{1}=\mathcal{V}^o\text{ and } K_{2}=\bigcup_{i=1,2,4,5} F_i\mathcal{V}^o,\]
then
\begin{equation}\label{eq31}
\begin{cases}
K_1=K_2\cup F_3K_1,\\
K_2=F_1K_2\cup F_2K_1\cup F_4K_1\cup F_5K_1.
\end{cases}
\end{equation}
Thus $\{K_1,K_2\}$ can be viewed as the invariant sets of a graph-directed construction $\mathcal{G}=(S, E, \Gamma)$ with $S=\{\mathcal{T}_1, \mathcal{T}_2\}$ being the state set, $E$ being the edge set consisting of $6$ edges, and $\Gamma=\{\psi_{i}\}_{i=1}^6$ being the collection of similitudes associated with $E$, that is, we could rewrite (\ref{eq31}) into
\begin{equation}\label{eq32}
\begin{cases}
K_1=\psi_1K_2\cup \psi_2K_1,\\
K_2=\psi_3K_2\cup \psi_4K_1\cup \psi_5K_1\cup \psi_6K_1.
\end{cases}
\end{equation}

For any edge $e\in E$, we write $\psi_e$ its associated similitude. Let $\Lambda$ be the collection of all finite walks in $G=(S,E)$ emanating from $\mathcal{T}_1$, including $\vartheta$ as the empty walk. Let $K_\vartheta=\mathcal{V}^o$, and for any walk $\bm{e}=e_1e_2\cdots e_n\in \Lambda$, denote $K_{\bm{e}}=\psi_{\bm{e}}K_{f(\bm{e})}$ a compact subset in $\mathcal{V}^o$. 
Then by Theorem 3.4, it is easy to check that the structure  $\{K_{\bm{e}},\Lambda_{\bm{e}}\}_{\bm{e}\in \Lambda}$ becomes an $f.r.f.t.$ nested structure of $\mathcal{V}^o$ with $M=2$, where $\Lambda_{\bm{e}}$ is the same notation introduced in the proof of Theorem 3.4. In addition, the boundary of $K_1$ is $\{q_1,q_2,q_3,q_4\}$ and the boundary of $K_2$ is $\{q_1,q_2,q_4,F_3q_1\}$, see Figure 4.3.

\begin{figure}[h]
	\centering
	\includegraphics[height=3cm]{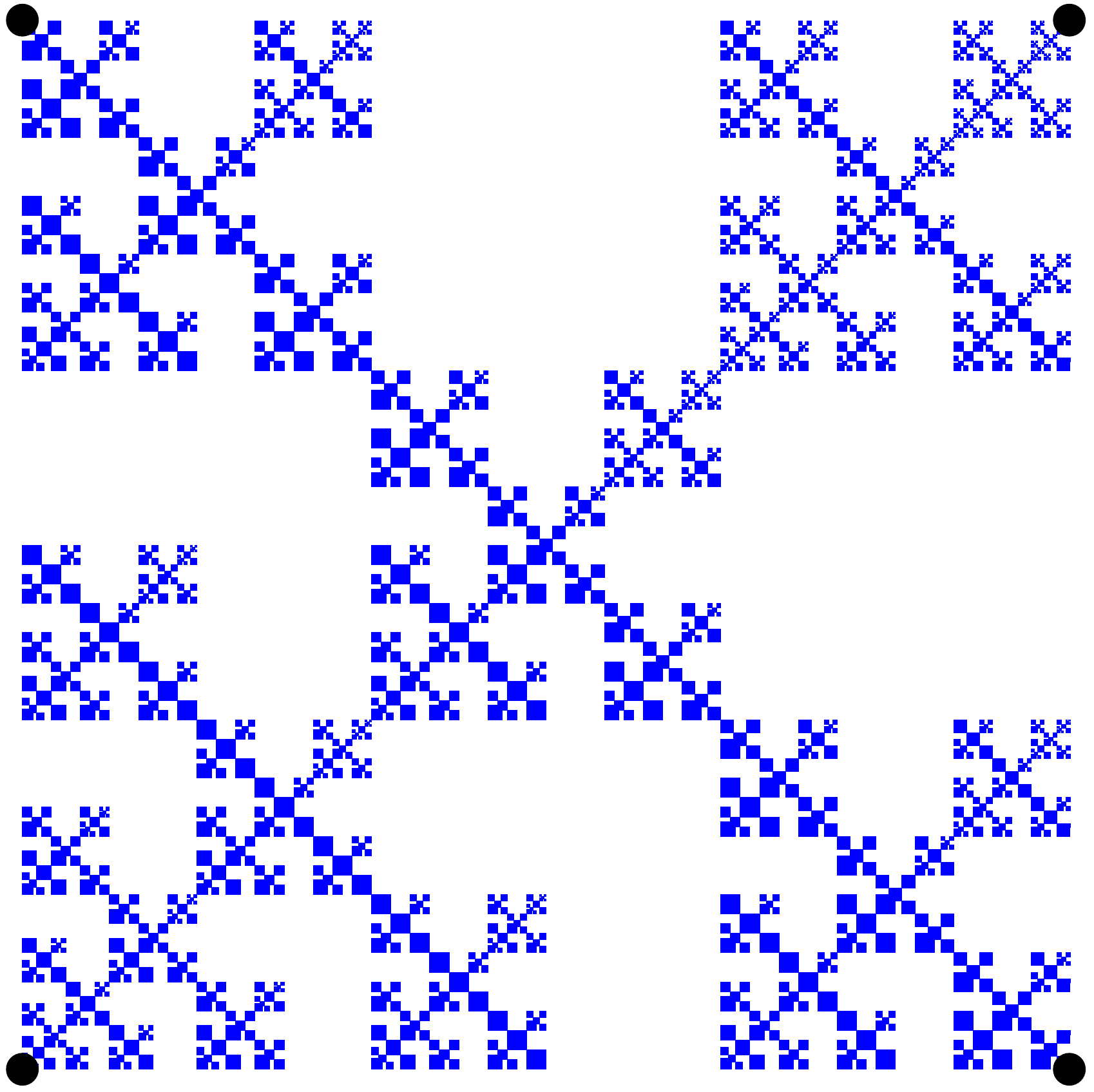}\qquad
	\includegraphics[height=3cm]{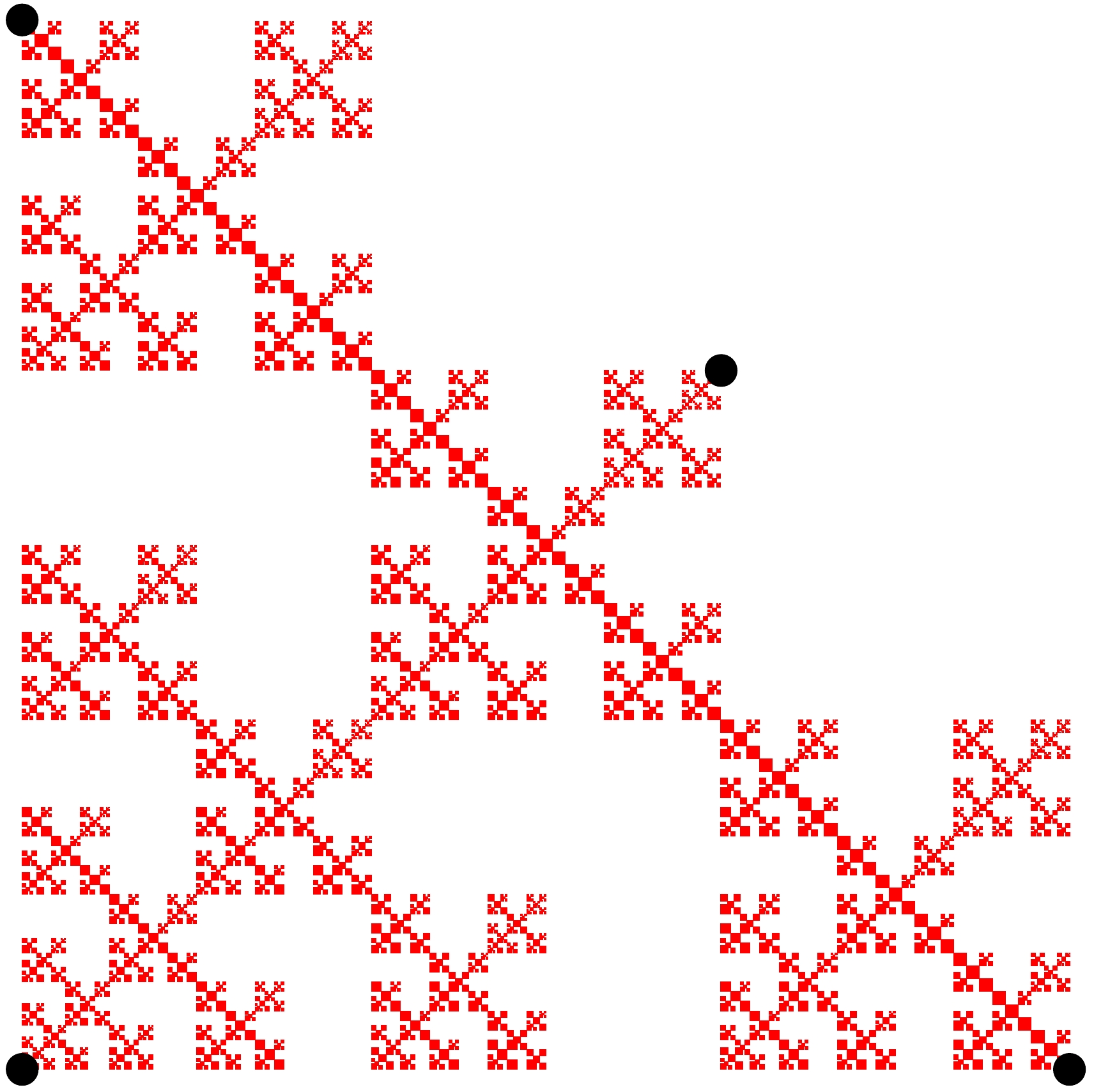}\qquad
	\includegraphics[height=3cm]{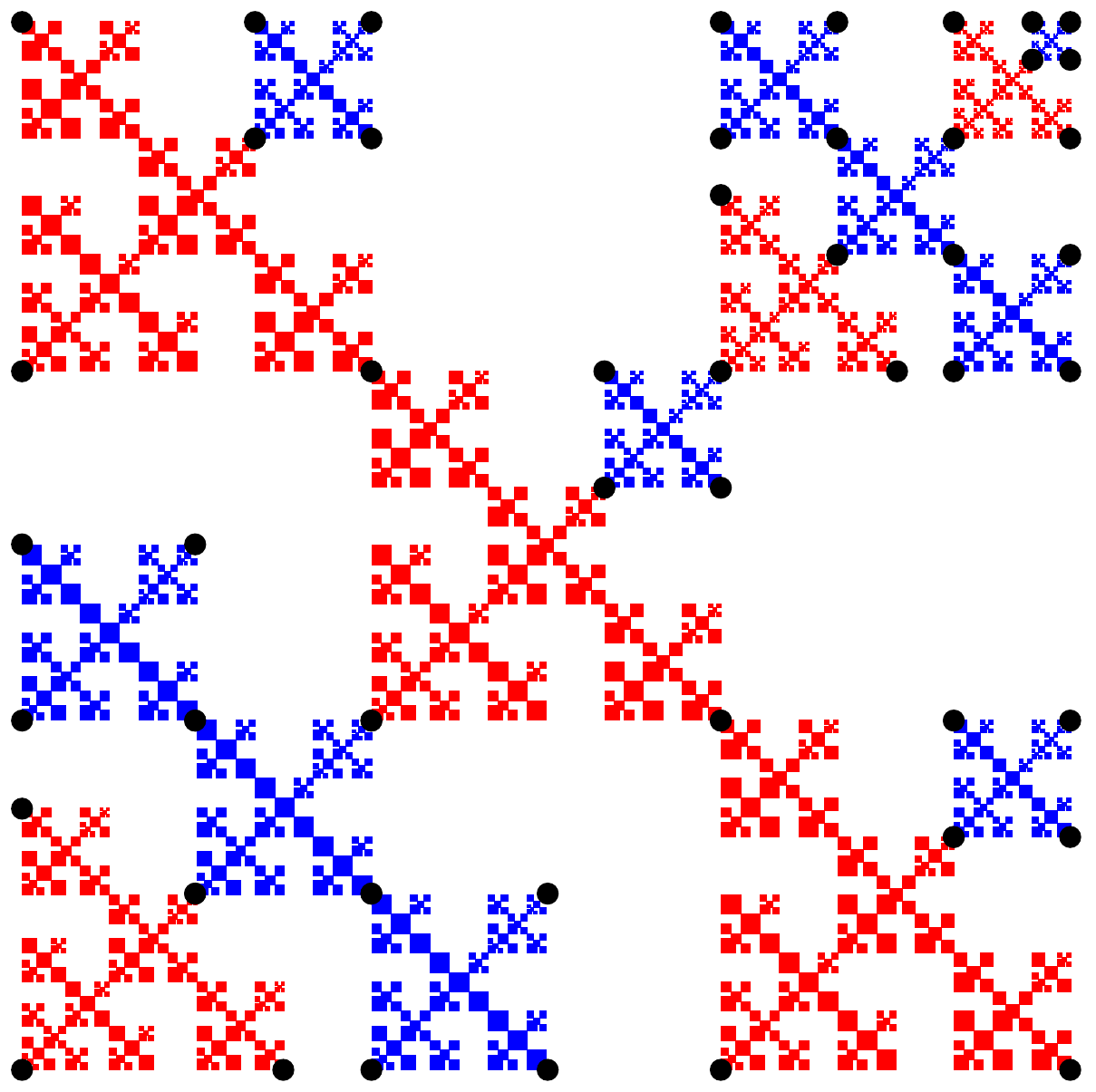}
	\begin{center}
		\caption{The boundary of $K_1,K_2$, and the level-$3$ vertices $V_3$.}
	\end{center}
\end{figure}

We would like to point out that the $f.r.f.t$  nested structure of $\mathcal{V}^o$ is not unique, which means that we could assign various $f.r.g.d.$ fractal  families that includes  $\mathcal{V}^o$ as a member. See Figure 4.4 for the illustration of another ``cutting rule'' of $\mathcal{V}^o$.

\begin{figure}[h]
	\centering
	\includegraphics[height=3cm]{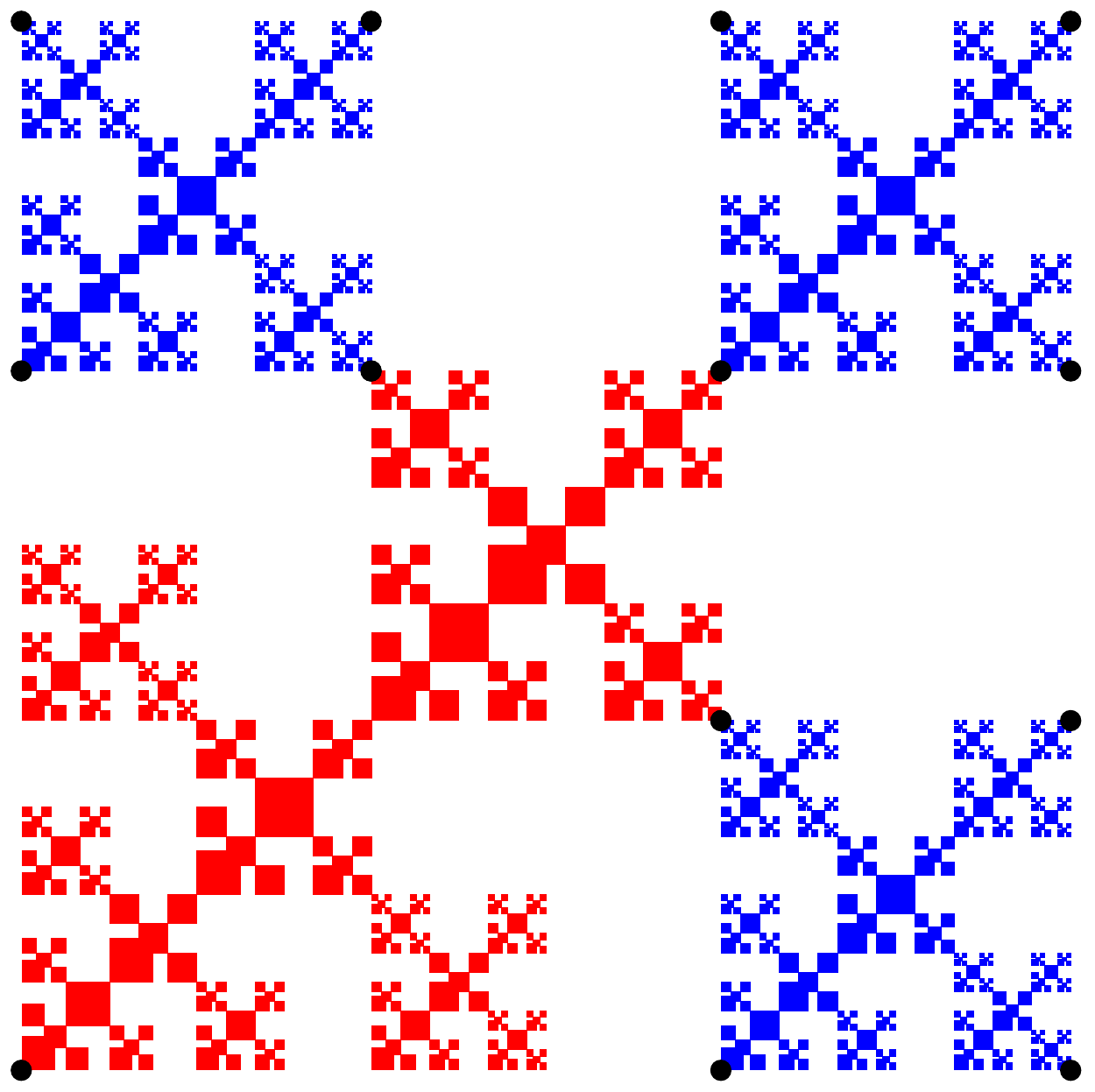}\qquad
	\includegraphics[height=3cm]{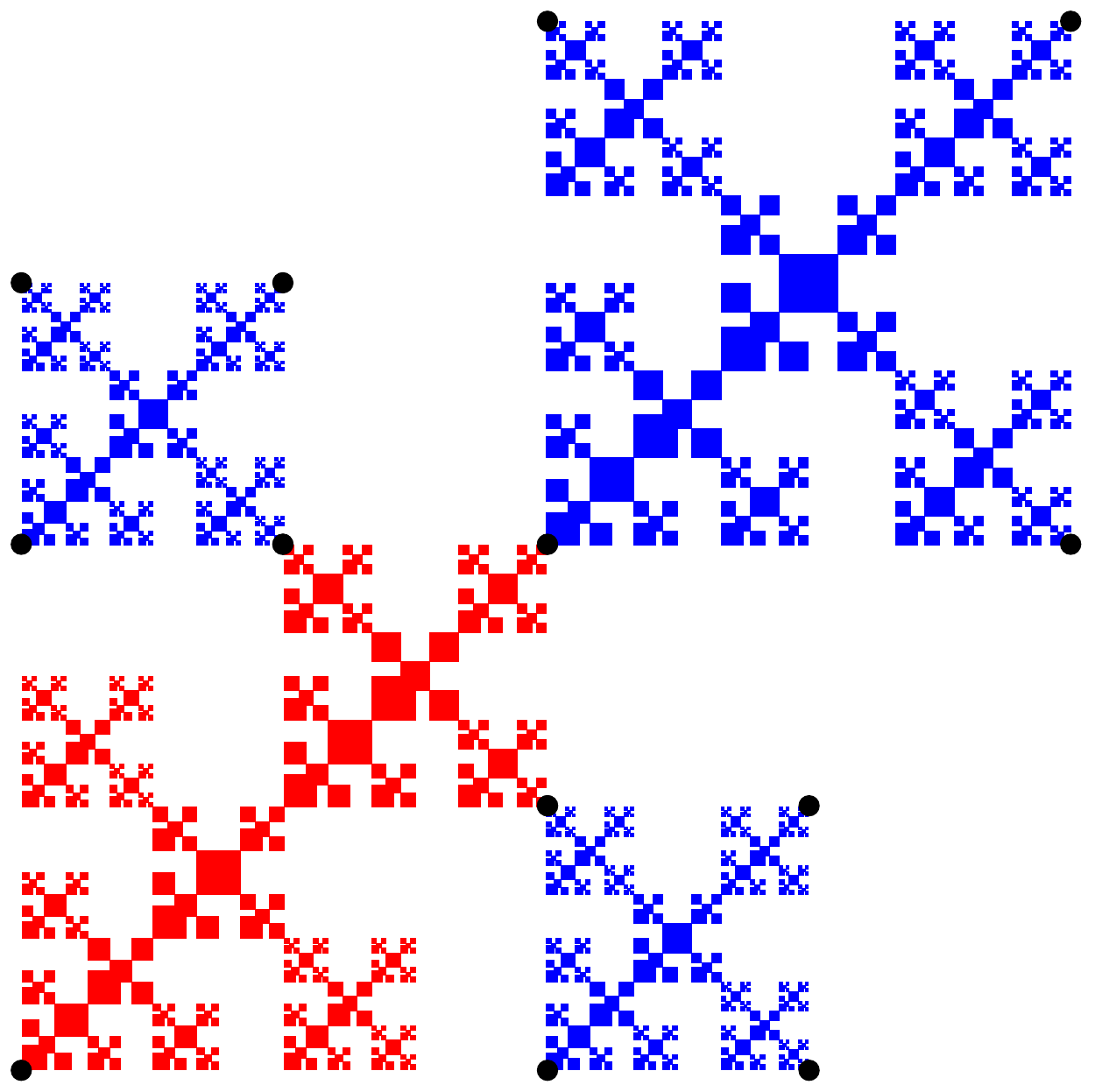}
	\begin{center}
		\caption{Another ``cutting rule'' of $\mathcal{V}^o$.}
	\end{center}
\end{figure}

\textbf{Example 2.(Overlapping gasket with open bottom)} 
Let $\{q_i\}_{i=1}^3$ be the vertices of an equilateral triangle, and $\{q_4,q_5\}$ be the centers of the line segments $\overline{q_1q_2}$ and $\overline{q_1q_3}$. The \textit{overlapping gasket with open bottom}, denoted by $\mathcal{SG}^o$, is the invariant set of the $i.f.s.$ $\{F_i\}_{i=1}^5$,
\[F_2: x\rightarrow\frac{1}{2} x+\frac{1}{2}q_2,\quad F_i: x\rightarrow \frac{1}{3} x+\frac{2}{3}q_i,i=1,3,4,5,\]
see Figure 1.3(left) for  $\mathcal{SG}^o$.

Let $K_1=\mathcal{SG}^o$ and  $K_2=(\mathcal{SG}^o\setminus F_2\mathcal{SG}^o)\cup \{F_2q_1\}$ be two connected compact subsets in $\mathcal{SG}^o$, then
\begin{equation}\begin{cases}
K_1=F_2K_1\cup K_2,\\
K_2=F_1K_1\cup F_3K_1\cup F_5K_1\cup F_4K_2.
\end{cases}
\end{equation}

It is easy to check that this provides an $f.r.f.t.$  nested structure of $\mathcal{SG}^o$. 
The boundary of $K_1$ is $\{q_1,q_2,q_3\}$, and the boundary of  $K_2$ is $\{q_1,q_3,q_4\}$.
See Figure 4.5 for an illustration.\\

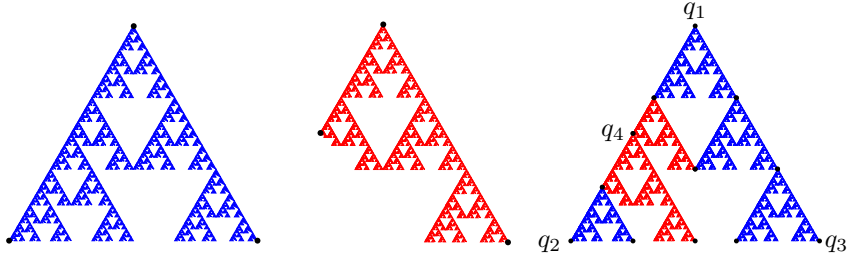
\begin{figure}[h]
	\centering
	\includegraphics[height=3cm]{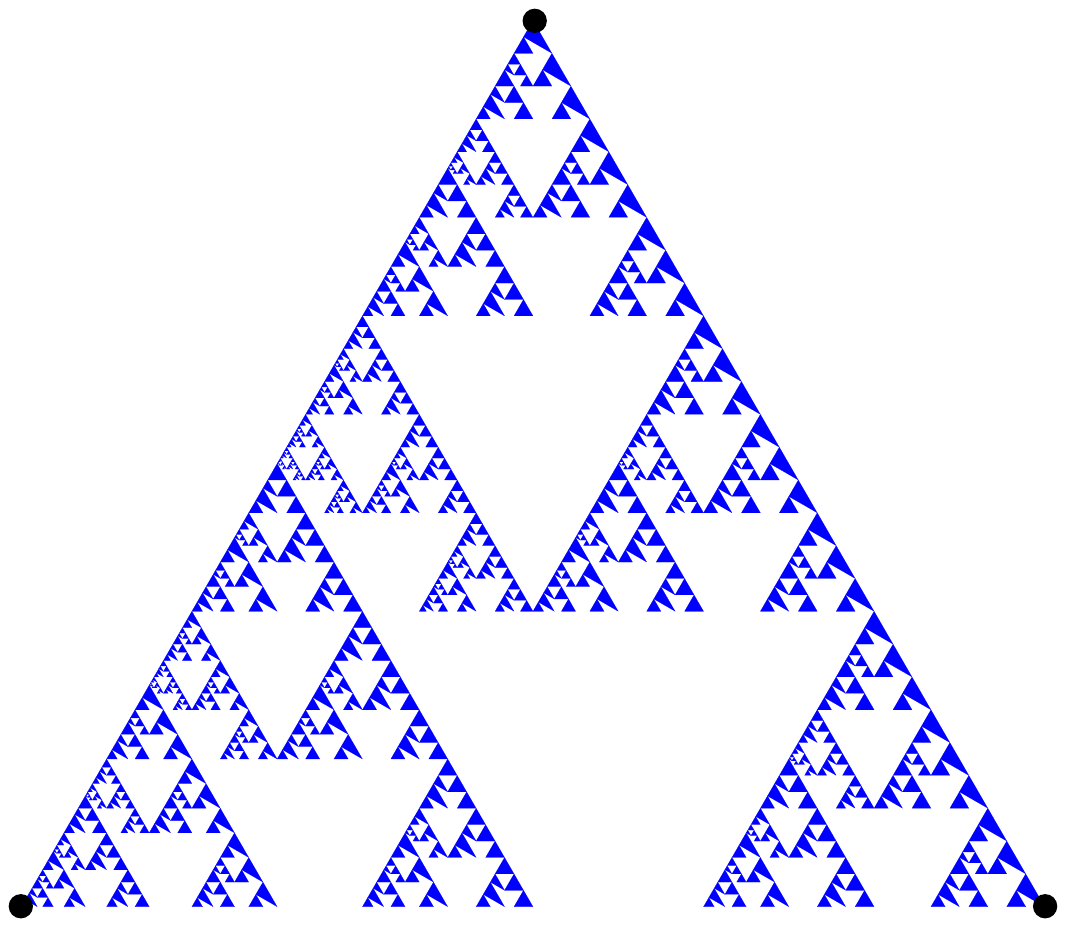}\qquad
	\includegraphics[height=3cm]{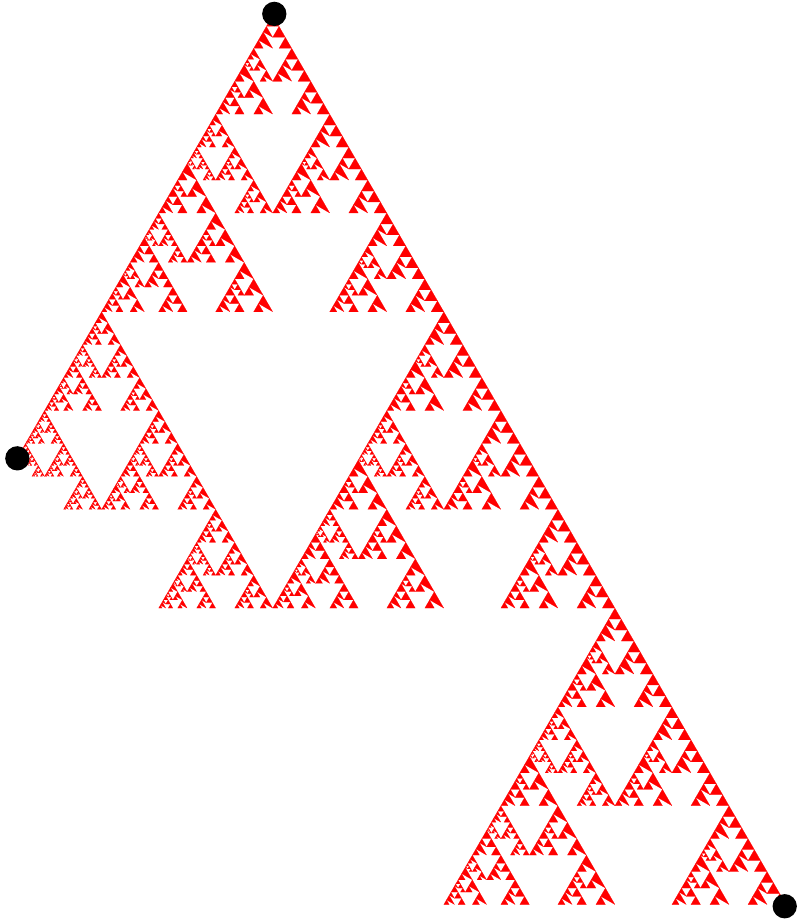}\qquad
	\includegraphics[height=3cm]{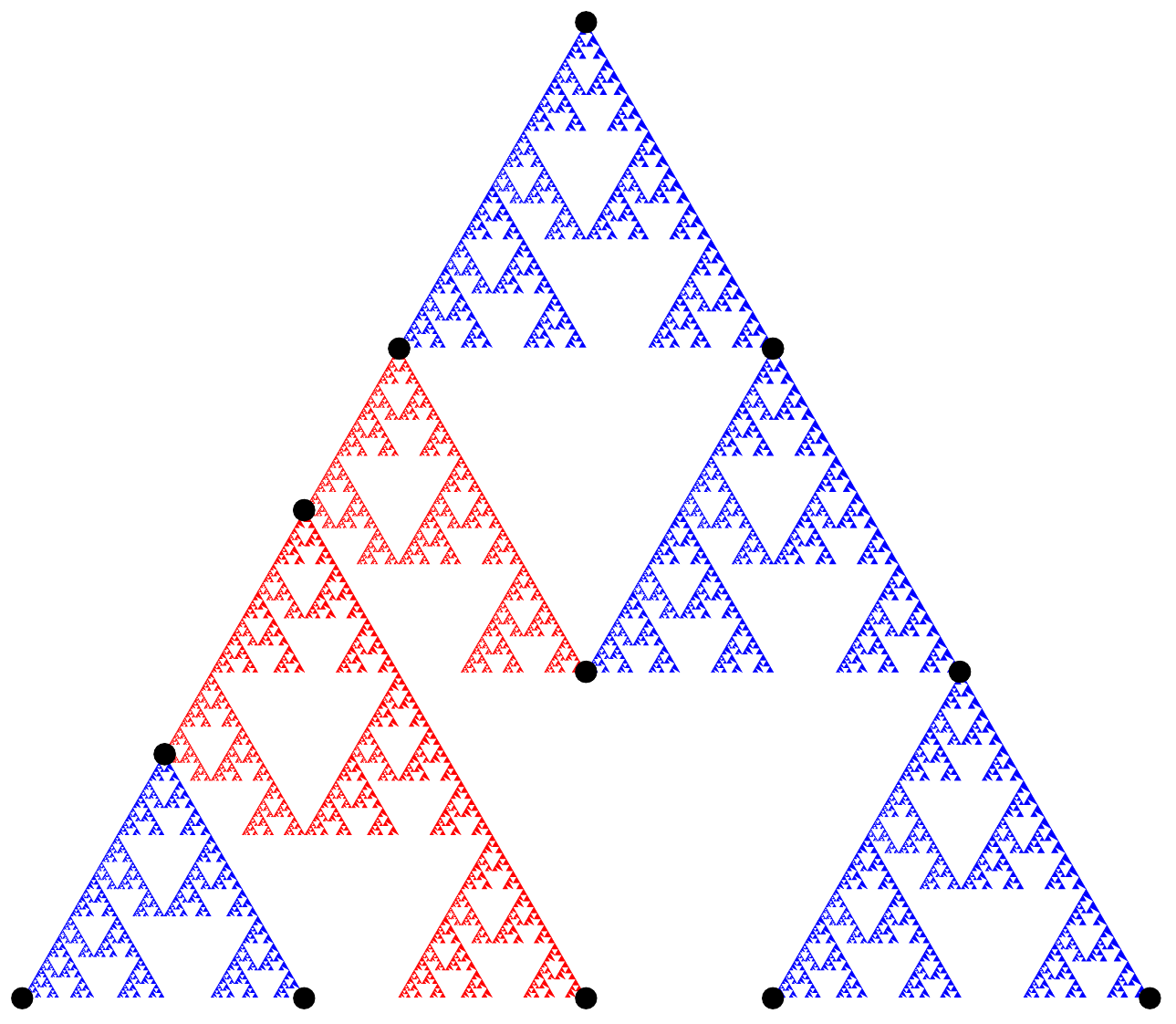}	
	\begin{picture}(0,0) \thicklines 
	\put(-57,88){$q_1$}
	\put(-112,0){$q_2$}
	\put(-3.5,0){$q_3$}
	\put(-88,42){$q_4$}
	\end{picture}
	\begin{center}
		\caption{ $K_1, K_2$ in $\mathcal{SG}^o$ and the level-$2$ division of $\mathcal{SG}^o$.}
	\end{center}
\end{figure}

\textbf{Example 3.(Vicsek windmill set)} 
Let $\{q_i\}_{i=1}^4$ be the vertices of a square in $\mathbb{R}^2$, say $\{(0,0),(1,0),(1,1),(0,1)\}$ for convenience. The \textit{Vicsek windmill set}, denoted by $\mathcal{V}^w$, is the invariant set of the $i.f.s.$ $\{F_i\}_{i=1}^8$,
\[\begin{aligned}
F_i(x)&=\frac{1}{4} x+\frac{3}{4}q_i,\text{ } i=1,2,3,4,\\
F_5(x)=\frac{1}{4} x&+(\frac{1}{4},0),
F_6(x)=\frac{1}{4} x+(\frac{1}{2},\frac{1}{4}),\\
F_7(x)=\frac{1}{4} x&+(\frac{1}{4},\frac{1}{2}),
F_8(x)=\frac{1}{4} x+(\frac{1}{2},\frac34),
\end{aligned}\]
see Figure 1.2(right) for $\mathcal{V}^w$. It is easy to check that $\mathcal{V}^w$ has an $f.r.f.t.$  nested structure with three types of islands, which are similar copies of
\[K_1=\mathcal{V}^w,\quad K_2=F_1\mathcal{V}^w\cup F_5\mathcal{V}^w,\quad K_3=F_1F_8\mathcal{V}^w\cup F_1F_3\mathcal{V}^w\cup F_5F_4\mathcal{V}^w,\]
where the latter two are the union of two or three copies of $\mathcal{V}^w$, see Figure 4.6. As the $f.r.g.d.$ construction of $\{K_1,K_2,K_3\}$ consists of long equations, we omit the exact expressions, but readers can get all the information from Figure 4.6. It is easy to see that the boundaries are 
\[\begin{aligned}
\partial K_1&=\{q_i\}_{i=1}^4,\quad \partial K_2=\{q_1,F_5q_2,F_5q_3,F_1q_4\},\\ &\partial K_3=\{F_1F_8q_1,F_5F_4q_2,F_5F_4q_3,F_1F_8q_4\},
\end{aligned}\]\\
which are the vertices of their associated rectangles.

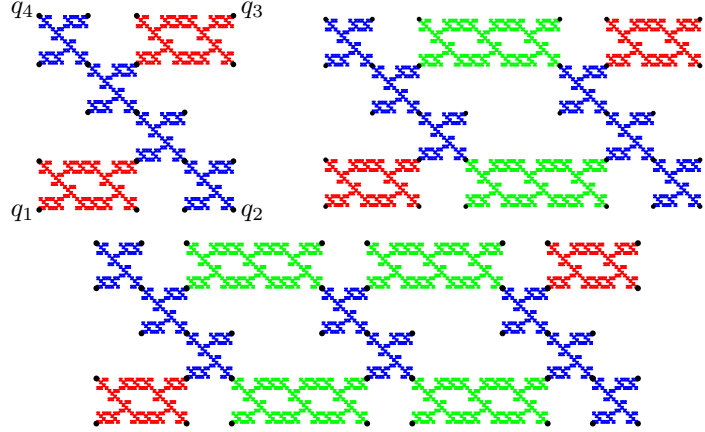
\begin{figure}[h]
	\centering
	\includegraphics[height=2.7cm]{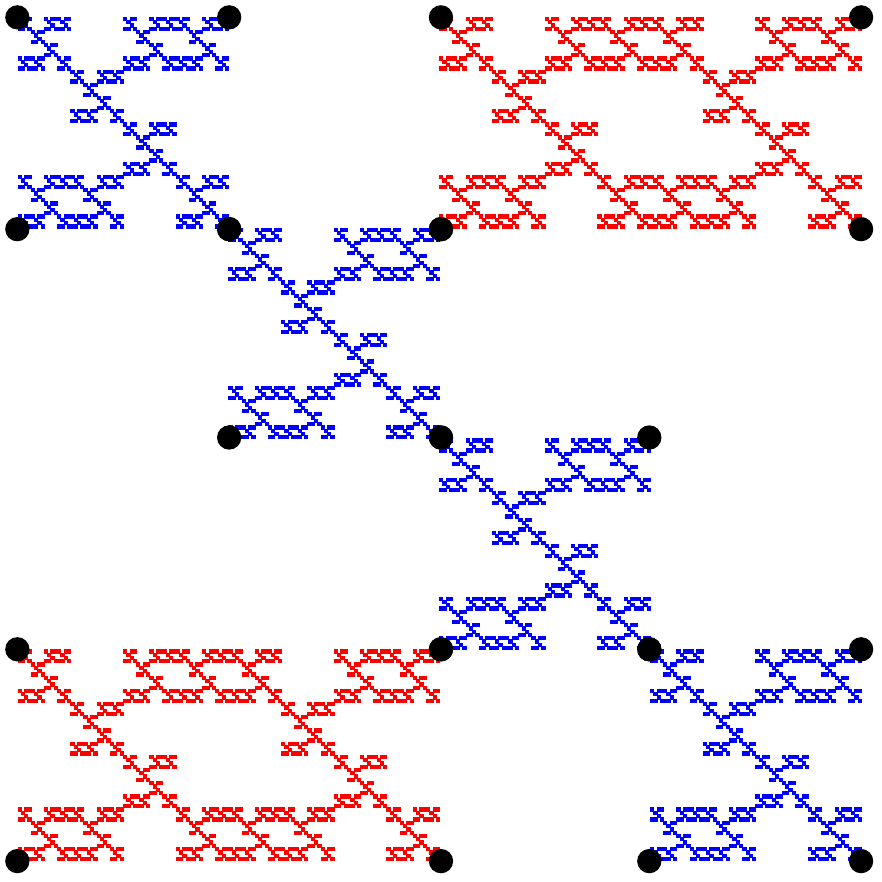}\qquad\quad
	\includegraphics[height=2.7cm]{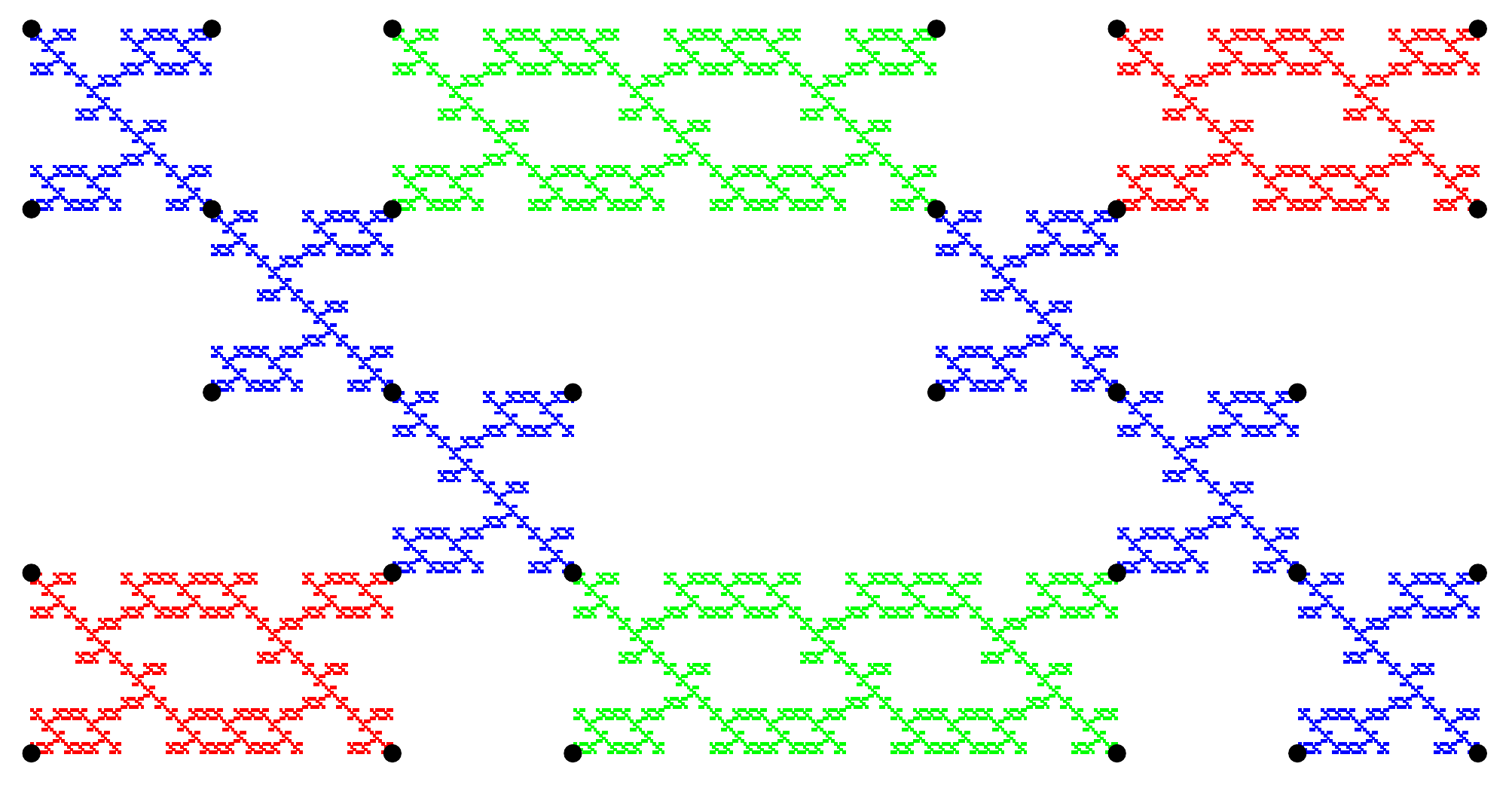}\\
	\vspace{0.2cm}
	\includegraphics[height=2.7cm]{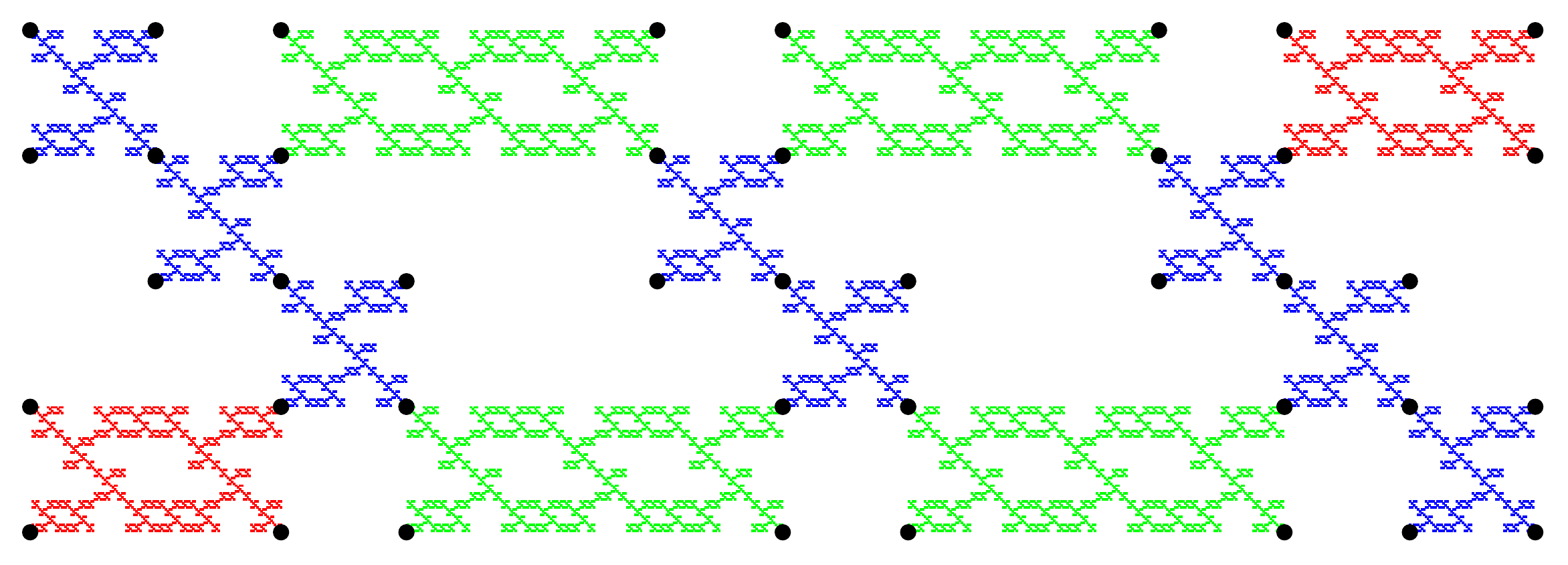}	
	\begin{picture}(0,0) \thicklines 
	\put(-245,83){$q_1$}
	\put(-158,83){$q_2$}
	\put(-245,160){$q_4$}
	\put(-158,160){$q_3$}
	\end{picture}
	\begin{center}
		\caption{Three typical islands in $\mathcal{V}^w$.}
	\end{center}
 
\end{figure}

\textbf{Example 4.(Symmetrical overlapping gasket with closed bottom)} 
Let $\{q_i\}_{i=1}^3$ be the vertices of an equilateral triangle, and $\{q_4,q_5\}$ be the centers of the line segments $\overline{q_1q_2}$ and $\overline{q_1q_3}$. The \textit{symmetrical overlapping gasket with closed bottom}, denoted by $\mathcal{SG}^c$, is the invariant set of the $i.f.s.$ $\{F_i\}_{i=1}^5$,
\[ F_i(x)=\frac{1}{3} x+\frac{2}{3}q_i,i=1,4,5, \quad F_i(x)=\frac{1}{2} x+\frac{1}{2}q_i, i=2,3,\]
see Figure 1.3(right) for $\mathcal{SG}^c$. 

Let $K_1=\mathcal{SG}^c, K_2=F_2\mathcal{SG}^c \cup F_4\mathcal{SG}^c$, then
\[\begin{cases}
K_1=F_1K_1\cup K_2\cup RK_2,\\
K_2=F_2K_2\cup F_2RK_2\cup F_4K_2\cup F_4RK_2\cup F_4F_1K_1,
\end{cases}\]
where $R$ is the reflection keeping $q_1$ fixed, and interchanging $q_2$ and $q_3$. This ``cutting rule'' provides an $f.r.f.t.$  nested structure of $\mathcal{SG}^c$. The boundaries of $K_1,K_2$ are
\[\partial K_1=\{q_1,q_2,q_3\},\quad\partial K_2=\{q_2,F_2q_3,F_4q_3,F_4q_1\}.\]
See Figure 4.7 for an illustration.\\

\begin{figure}[h]
	\centering
	\includegraphics[height=3cm]{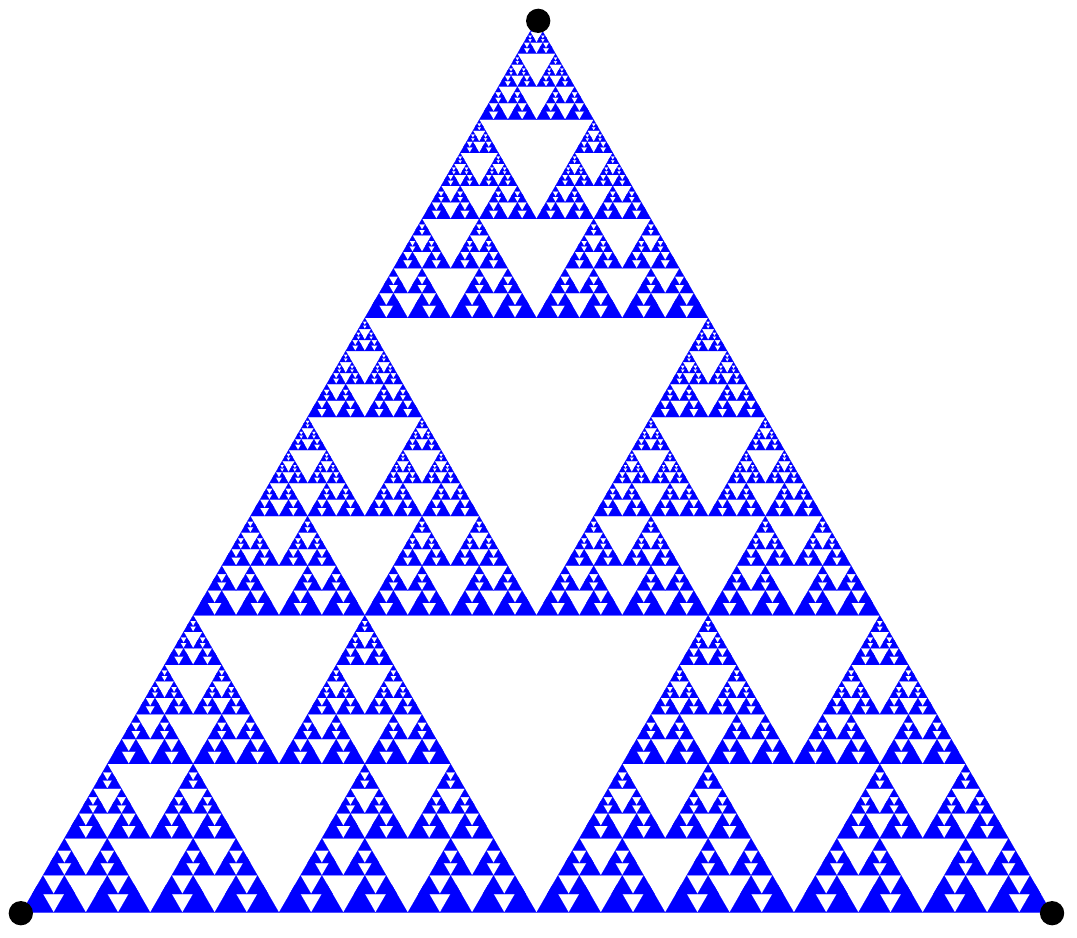}\qquad
	\includegraphics[height=3cm]{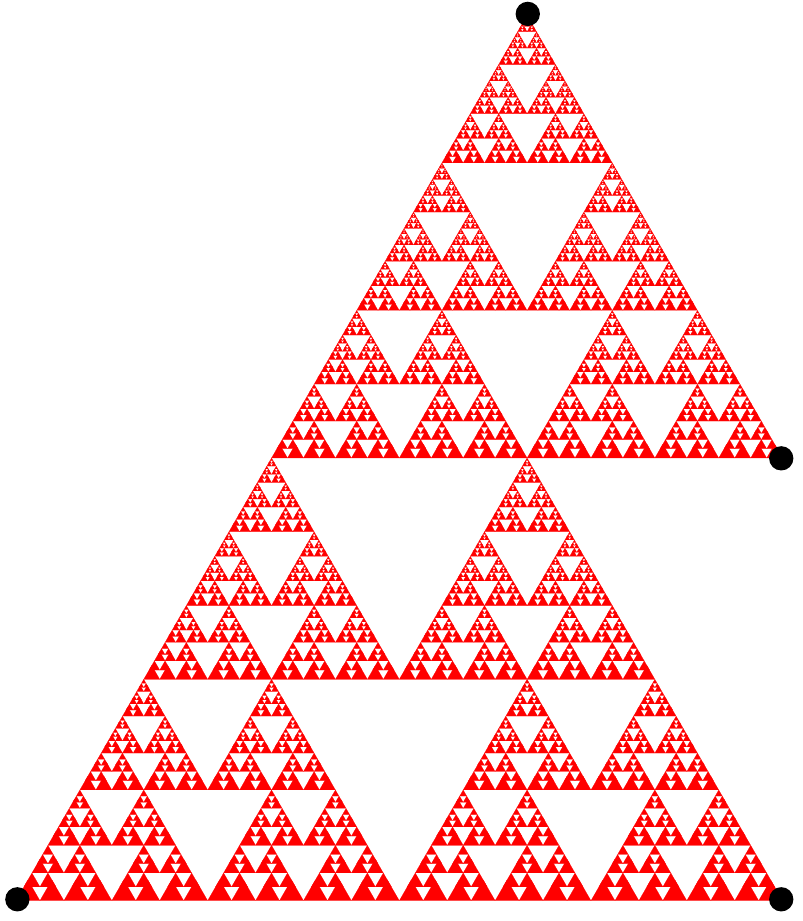}\qquad
	\includegraphics[height=3cm]{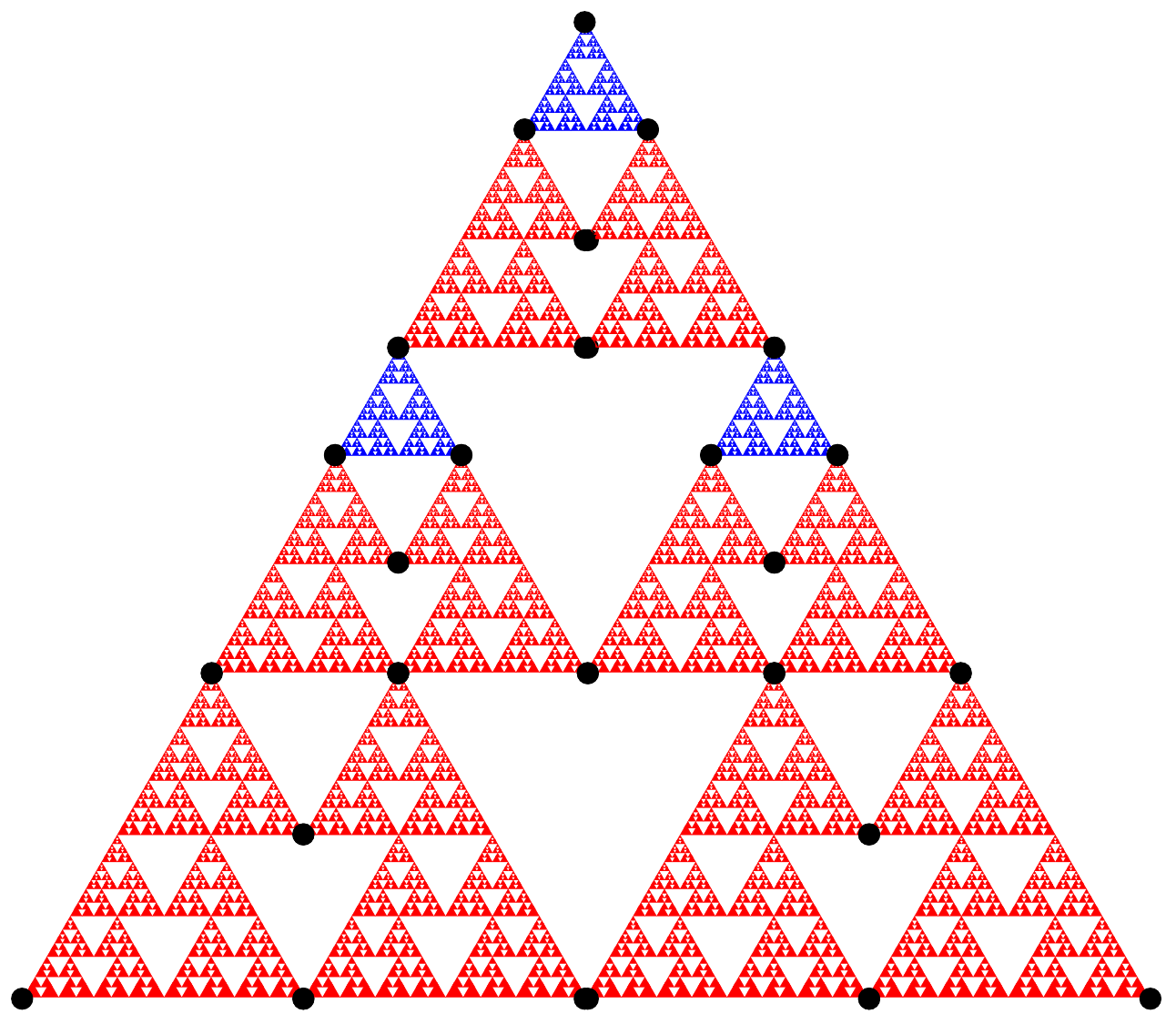}
	\begin{picture}(0,0) \thicklines 
	\put(-57,88){$q_1$}
    \put(-112,0){$q_2$}
    \put(-3.5,0){$q_3$}
    \end{picture}
		
	\begin{center}
		\caption{$K_1, K_2$ in $\mathcal{SG}^c$  and the level-$2$ division of $\mathcal{SG}^c$.}
	\end{center}
\end{figure}
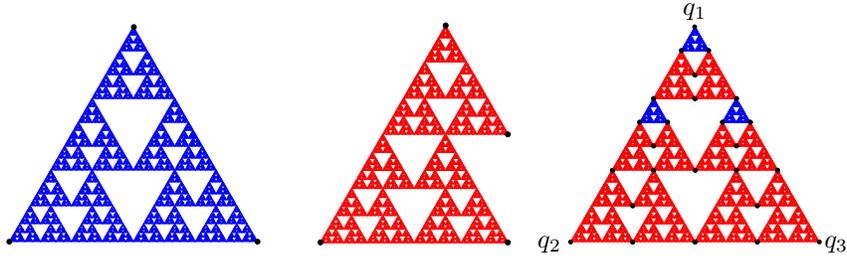

\section{Dirichlet forms on $f.r.f.t.$ self-similar sets}
In this section, we construct ``self-similar'' Dirichlet forms on self-similar sets with $f.r.f.t.$ nested structures, and show the solutions of typical examples. Here we use the term \textit{self-similar} to refer to the scaling invariance of the Dirichlet forms associated with  the $f.r.f.t.$ nested  structures. The construction essentially comes from the method that Kigami introduced to deal with $p.c.f.$ self-similar sets [K2]. 

Let's recall some basic concepts that will be used in this section.

 For a finite set $V$, we use  $l(V)$ to denote the collection of all maps from $V$ into $\mathbb{R}$.
A symmetric linear operator(matrix) $H: l(V)\to l(V)$ is called a \textit{(discrete) Laplacian} on $V$ if $H$ is non-positive definite, $Hu=0$ if and only if $u$ is a constant on $V$, and  $H_{xy}\geq 0$ for all $x\neq y\in V$. Write $x\sim y$ if $H_{xy}>0$, there is a symmetric bilinear form $\mathcal{E}_H(\cdot,\cdot)$ on $l(V)$, called the \textit{(discrete) Dirichlet form} associated with $H$, written as
\[\mathcal{E}_H(u,v)=\sum_{x\sim y}c_{x,y}\big(u(x)-u(y)\big)\big(v(x)-v(y)\big),\quad \forall u,v\in l(V),\]
with $c_{x,y}=H_{xy}$ called the \textit{conductance} between $x,y$.  Conversely, for a symmetric bilinear form $\mathcal{E}(\cdot,\cdot)$ on $l(V)$ as above with $c_{x,y}>0$ and $\mathcal{E}(u,u)=0$ if and only if $u$ is a constant on $V$, there is a unique Laplacian $H$ on $V$ such that $\mathcal{E}=\mathcal{E}_H$.
The pair $(V,H)$ is called a \textit{(finite) resistance network}. We write $\mathcal{E}_H(u)=\mathcal{E}_H(u,u)$ for short. 

{If $(V,H_1)$ and $(U,H_2)$ are two pairs of  resistance networks satisfying that $V\subset U$ and}
\begin{equation}\label{eq51}
\mathcal{E}_{H_1}(v)=\min\{\mathcal{E}_{H_2}(u):u\in l(U),u|_{V}=v\},\quad \forall v\in l(V),
\end{equation}
we say they are \textit{compatible} and write $(V,H_1)\leq (U,H_2)$.
Note that if $(V,H_1)\leq (U,H_2)$, then for any $v\in l(V)$, there exists a unique function $u\in l(U)$ attaining the minimum in (\ref{eq51}).

Given a compatible sequence $\{(V_n,H_n)\}_{n\geq 0}$, i.e. $(V_n,H_n)\leq (V_m,H_m)$ if $n\leq m$, there is a limit form $(\mathcal{E},\mathcal{F})$, called the \textit{resistance form} on $V_*=\bigcup_{n=0}^\infty V_n$, defined as
\[\mathcal{E}(u,v)=\lim_{n\to\infty} \mathcal{E}_{H_n}(u|_{V_n},v|_{V_n}),\quad \forall u,v\in \mathcal{F},\]
with
\[\mathcal{F}=\{u\in l(V_*): \mathcal{E}(u):=\lim_{n\rightarrow\infty}\mathcal{E}_{H_n}(u)<\infty\}.\]

The  form $(\mathcal{E},\mathcal{F})$ can be naturally  extended to  $l\big(cl_R(V_*)\big)$, still denoted by $(\mathcal{E},\mathcal{F})$,  where $cl_R(V_*)$  is the closure of $V_*$ with respect to the \textit{effective resistance metric} $R(\cdot,\cdot)$, \[R(x,y)=\inf\{\mathcal{E}(u)|u(x)=0,u(y)=1\}^{-1}, \quad  x\neq y\in V_*.\] 

For the resistance form $(\mathcal{E},\mathcal{F})$ on $cl_R(V_*)$, for any subset $V$, the \textit{trace} of $\mathcal{E}$ onto $V$, denoted by  $\mathcal{E}|_V$, is defined as the unique Dirichlet form on $l(V)$ satisfying
\[\mathcal{E}|_V(v)=\min\{\mathcal{E}(u):u|_{V}=v, u\in\mathcal{F}\},\quad \forall v\in l(V).\] 
In particular, the form $\mathcal{E}_{H_n}$ is the trace of $\mathcal{E}$ onto $V_n$, $\forall n\geq 0$. 

The above concepts can be found in details in [K5-K7].
In the following, for a self-similar set $K$ possessing an $f.r.f.t.$ nested structure, we will construct  a ``self-similar'' compatible  sequence of networks on $\{V_n\}_{n\geq 0}$, the sequence of finite vertices  introduced in Section 2. Providing this can be fulfilled and in addition $cl_R(V_*)=K$ in the sense of homeomorphism, then using standard argument, associated with any  Borel probability measure $\mu$(we often require $\mu$ has a scaling invariant property) on $K$, the resistance form $(\mathcal{E}, \mathcal{F})$ turns out to be a local regular (self-similar) Dirichlet form on $L^2_\mu(K)$. \\

Recall that as in Section 3, for a self-similar set $K$ with an $f.r.f.t.$ nested structure $\{K_\alpha,\Lambda_\alpha\}_{\alpha\in \Lambda}$, there is an associated graph-directed construction $\mathcal{G}=(S,E,\{\psi_e\}_{e\in E})$ and an $f.r.g.d.$ fractal family $\mathcal{K}=\{K_{\alpha_i}\}_{i=1}^M$, including $K=K_{\alpha_1}$ as a member, satisfying the identity (3.1) and (3.2). For $\gamma=P(\beta)$, we use the notation $e(\gamma,\beta)$ to specify the associated edge in $E$ from $\mathcal{T}_{t(\gamma)}$ to $\mathcal{T}_{t(\beta)}$. It is easy to check  that $\psi_{e(\gamma,\beta)}=\phi_{\gamma,\alpha_{t(\gamma)}}\circ\phi_{\alpha_{t(\beta)},\beta}$. For $\gamma=P^n(\beta)$ with $n\geq 1$, we write $\bm{e}(\gamma,\beta)=e_1e_2\cdots e_n$ with $e_i=e\big(P^{n+1-i}(\beta),P^{n-i}(\beta)\big)$, $1\leq i\leq n$, the walk from $\mathcal{T}_{t(\gamma)}$ to $\mathcal{T}_{t(\beta)}$ for simplicity.

For  $1\leq i\leq M$, Let  $D_{\alpha_i}$ be an initial Laplacian on $V_{\alpha_i}$ and denote $\mathcal{E}_{\alpha_i,0}(\cdot,\cdot)$ its associated Dirichlet form. Note that it may be possible that $V_{\alpha_1}=V_\vartheta=\emptyset$. If it happens, we just skip $D_{\alpha_1}$. Let $\bm{r}=(r_e)_{e\in E}$ be a vector of positive numbers. For $1\leq i\leq M$, as shown in Proposition 2.3(d),  $V_{\alpha_i}\subset V_{\alpha_i,1}$, we define a Dirichlet form on $l(V_{\alpha_i,1})$,
\begin{equation}
\mathcal{E}_{\alpha_i,1}(u,v)=\sum_{\beta\in\Lambda_{\alpha_i}} r_{e(\alpha_i,\beta)}^{-1} \mathcal{E}_{\alpha_{t(\beta)},0} (u\circ \phi_{\alpha_{t(\beta)},\beta},v\circ \phi_{\alpha_{t(\beta)},\beta}), \quad\forall u,v\in l(V_{\alpha_i,1}),
\end{equation}
and let $H_{\alpha_i,1}$ be its associated Laplacian.

In general, for $\alpha\in \Lambda$, $n\geq 1$, define $V_{\alpha,n}=\bigcup_{\beta\in P^{-n}(\alpha)}V_{\beta}$ (obviously $V_{\alpha_1,n}=V_n$), then by Proposition 2.3(d), $V_{\alpha,n-1}\subset V_{\alpha,n}$. For $n\geq 2$, $1\leq i\leq M$, we define a Dirichlet form $\mathcal{E}_{\alpha_i,n}(\cdot,\cdot)$ on $l(V_{\alpha_i,n})$ inductively,
\[\begin{aligned}
\mathcal{E}_{\alpha_i,n}(u,v)&=\sum_{\beta\in\Lambda_{\alpha_i}}r_{e(\alpha_i,\beta)}^{-1} \mathcal{E}_{\alpha_{t(\beta)},n-1} (u\circ \phi_{\alpha_{t(\beta)},\beta},v\circ \phi_{\alpha_{t(\beta)},\beta})\\
&=\sum_{\beta\in P^{-n}({\alpha_i})}r_{\bm{e}(\alpha_i,\beta)}^{-1} \mathcal{E}_{\alpha_{t(\beta)},0} (u\circ \phi_{\alpha_{t(\beta)},\beta},v\circ \phi_{\alpha_{t(\beta)},\beta}), \quad\forall u,v\in l(V_{\alpha_i,n}),
\end{aligned}\] 
and let $H_{\alpha_i,n}$ be the associated Laplacian, where we write $r_{\bm{e}}=r_{e_1}r_{e_2}\cdots r_{e_{|\bm{e}|}}$ for simplicity.

Analogous to that for the $p.c.f.$ self-similar sets, in order to make the sequence $\{(V_n, H_{n})\}_{n\geq 1}=\{(V_{\alpha_1,n}, H_{\alpha_1,n})\}_{n\geq 1}$ compatible, we only need to assume that the pair $(\{D_{\alpha_i}\}_{i=1}^M, \bm{r})$ satisfies the \textit{renormalization equations}
\begin{equation}
(V_{\alpha_i},D_{\alpha_i})\leq (V_{\alpha_i,1},H_{\alpha_i,1}),\quad\forall 1\leq i\leq M.
\end{equation}

\textbf{Definition 5.1.} \textit{Call the pair $(\{D_{\alpha_i}\}_{i=1}^M,\bm{r})$ a} harmonic structure \textit{with respect to the $f.r.f.t.$  nested structure $\{K_\alpha,\Lambda_\alpha\}_{\alpha\in \Lambda}$ if the renormalization equations (5.3) are satisfied. In addition, we say  $(\{D_{\alpha_i}\}_{i=1}^M,\bm{r})$ is} regular\textit{ if $r_{\bm{e}}<1$ for every cycle $\bm{e}$ in the directed graph $G=(S,E)$.}

For a harmonic structure $(\{D_{\alpha_i}\}_{i=1}^M, \bm{r})$, the elements in $\bm{r}$ are called \textit{renormalization factors}. Providing a harmonic structure exists, we actually get a compatible sequence of networks $\{(V_{\alpha_i,n}, H_{\alpha_i,n})\}_{n\geq 1}$ for each $1\leq i\leq M$. The limit forms are denoted by $\mathcal{E}_{\alpha_i}, 1\leq i\leq M$. In particular, we write $\mathcal{E}=\mathcal{E}_{\alpha_1}$, $$
\mathcal{F}=\{u\in l(V_*): \mathcal{E}(u)<\infty\},$$ and call $(\mathcal{E},\mathcal{F})$ a \textit{self-similar resistance form} induced by $(\{D_{\alpha_i}\}_{i=1}^M, \bm{r})$.

\textbf{Proposition 5.2.} \textit{Suppose $(\{D_{\alpha_i}\}_{i=1}^M, \bm{r})$ is a regular harmonic structure of an $f.r.f.t.$ self-similar set $K$. There exist two constants $0<\rho_1<\rho_2<1$ and $C_1, C_2>0$, such that 
\[C_1\rho_1^{|\bm{e}|}\leq r_{\bm{e}}\leq C_2\rho_2^{|\bm{e}|}\]
holds for any cycle $\bm{e}$ in the directed graph $G=(S,E)$.}

\textit{Proof.} The lower bound estimate is obvious as we just need to take $\rho_1=\min\{r_e\}_{e\in E}$ and $C_1=1$. Conversely, assume $\bm{e}=e_1 e_2\cdots e_n$ is a cycle whose length $
n\geq M$, then it will contain at least a cycle $\bm{e}_1=e_{k}e_{k+1}\cdots e_{k+m}$, with $m<M$, $k\geq 1$. By deleting this cycle from $\bm{e}$ we will get a new cycle $\bm{e'}=e_1\cdots e_{k-1}e_{k+m+1}\cdots e_n$ with $r_{\bm{e}'}=r_{\bm{e}}/r_{\bm{e}_1}$. Repeat the same operation on $\bm{e}'$ and continue, until the length of new cycle is less than $M$. We can extract at least  $[\frac{n}{M}]$ cycles from $\bm{e}$ whose lengths are all no more than $M$. Note that $\lambda:=\max\{r_{\bm{e}}:|\bm{e}|\leq M, \bm{e}\text{ is a cycle}\}<1$ as $(\{D_{\alpha_i}\}_{i=1}^M, \bm{r})$ is regular. Thus by letting $\rho_2=\lambda^{\frac1M}$ and choose $C_2$ appropriately, we get the upper bound estimate.
\hfill$\square$

Thus if the harmonic structure is regular,  by Proposition 5.2, using routine arguments(see for example [K5, K7]), we can find that the induced effective resistance metric $R(\cdot,\cdot)$ of $(\mathcal{E},\mathcal{F})$  on $V_*$ is equivalent to the Euclidean metric which yields $cl_R(V_*)=K$. Hence the form $(\mathcal{E},\mathcal{F})$ turns out to be a local regular Dirichlet form on $L^2_\mu(K)$ for any Borel probability measure $\mu$ on $K$. Without causing any confusion, we still denote the resulting Dirichlet form by $(\mathcal{E},\mathcal{F})$. 
Even if the harmonic structure is not regular, it may still be possible to generate a local regular Dirichlet form on $L^2_\mu(K)$ providing more assumptions are made on the measure $\mu$. The argument is similar to the $p.c.f.$ case [K5], see also a discussion for the $f.r.g.d.$ fractals in [HN], we omit it.  From now on we will only interest in those $f.r.f.t.$ self-similar sets possessing  regular harmonic structures.

In the following parts, we will construct regular harmonic structures on the examples shown in previous sections. 

Before proceeding, we will mention an interesting fact about cut-points in the fractals that will simplify the calculation. Let $K$ be an $f.r.f.t.$ self-similar set, $\mathcal{K}=\{K_{\alpha_i}\}_{i=1}^M$ be its associated $f.r.g.d.$ fractal family including $K=K_{\alpha_1}$ as a member. For $1\leq i\leq M$, $p\in K_{\alpha_i}$, we call $p$ a \textit{cut-point} of $K_{\alpha_i}$ if $K_{\alpha_i}\setminus\{p\}$ is disconnected. Note that since $K_{\alpha_i}$ is connected, $K_{\alpha_i}\setminus \{p\}$ is a locally arcwise connected set.

\textbf{Lemma 5.3.} \textit{Let $\{K_\alpha,\Lambda_\alpha\}_{\alpha\in\Lambda}$ be an $f.r.f.t.$  nested structure with a regular harmonic structure $(\{D_{\alpha_i}\}_{i=1}^M,\bm{r})$. For $1\leq i\leq M$, let $p$ be a cut-point of $K_{\alpha_i}$, and $\{p_k\}_{k=1}^m$ be a finite set in $K_{\alpha_i}\setminus\{p\}$. Write the restriction of $\mathcal{E}_{\alpha_i}$ onto $V=\{p\}\cup\{p_k\}_{k=1}^m$ by}
\[\mathcal{E}_{\alpha_i}|_{V}(u,v)=\sum_{x\neq y\in V}c_{x,y}\big(u(x)-u(y)\big)\big(v(x)-v(y)\big),\quad \forall u,v\in l(V),\]
then $c_{p_k,p_l}>0$ if and only if $p_k,p_l$ belong to a same connected component of $K_{\alpha_i}\setminus\{p\}$.

\textit{Proof. }  The proof is obvious and routine by a standard discussion of harmonic structures. We omit it. \hfill$\square$




\subsection{Vicsek set with overlaps}
In this first example, we will give all the regular harmonic structures associated with the $f.r.f.t.$ nested structure of $\mathcal{V}^o$ given by (4.2) in Section 4. The same notations introduced  there will be used.

\textbf{Lemma 5.4.} \textit{The point $q_5$ is a cut-point of $\mathcal{V}^o$. In addition, $\mathcal{V}^o\setminus \{q_5\}$ has $4$ connected components, and each contains one element in $\{q_1,q_2,q_3,q_4\}$.}

\textit{Proof.} For any $n\geq 1$, $\mathcal{V}^o\setminus F_5^n\mathcal{V}^o$ is disconnected, with four connected components, each containing one element in $\{q_1,q_2,q_3,q_4\}$. Denote by $C_{i,n}$ the component in $\mathcal{V}^o\setminus F_5^n\mathcal{V}^o$ containing $q_i$. 

Then for each path $\gamma$ from $q_i$ to $q_j$(a continuous curve $\gamma: [0,1]\rightarrow \mathcal{V}^o$ such that $\gamma(0)=q_i$, $\gamma(1)=q_j$) with $i\neq j$, it intersects $F_5^n\mathcal{V}^o$ for any $n\geq 1$, so the path $\gamma$ contains the point $q_5$.  Thus $q_i$, $q_j$ belong to different connected components of $\mathcal{V}^o\setminus \{q_5\}$ since it is locally arcwise connected. 
On the other hand, for $1\leq i\leq 4$, let $C_i=\bigcup_{n=1}^\infty C_{i,n}$, then $C_i$ is obviously  connected and contains $q_i$. So $C_i$'s are all the connected components of $\mathcal{V}^o\setminus \{q_5\}$.\hfill$\square$\\

For the consistency of the notations, we write $K_{\alpha_1}=K_1$, $K_{\alpha_2}=K_2$ and denote $V_{\alpha_1}$, $V_{\alpha_2}$ as their boundaries. Taking into account the above lemma, we include the cut-point $q_5$ when we deal with $V_{\alpha_1}$, $V_{\alpha_2}$ for convenence. Let \[\tilde{V}_{\alpha_1}=V_{\alpha_1}\cup\{q_5\}=\{q_1,q_2,q_3,q_4,q_5\},\text{ and }\tilde{V}_{\alpha_2}=V_{\alpha_2}\cup \{q_5\}=\{q_1,q_2,q_4,q_5,F_3q_1\}.\]
It is easy to check that
\[\begin{cases}
\tilde{V}_{\alpha_1}\subset \psi_1\tilde{V}_{\alpha_2}\cup \psi_2\tilde{V}_{\alpha_1},\\
\tilde{V}_{\alpha_2}\subset \psi_3\tilde{V}_{\alpha_2}\cup \psi_4\tilde{V}_{\alpha_1}\cup \psi_5\tilde{V}_{\alpha_1}\cup \psi_6\tilde{V}_{\alpha_1}.
\end{cases}\]
Define the analogy of a harmonic structure on $\tilde{V}_{\alpha_i}$'s, denoted by $(\{\tilde{D}_{\alpha_i}\}_{i=1}^M,\bm{r})$, by requiring the analogy of  $(5.3)$. By the elementary operation of resistance networks, we have

\textbf{Proposition 5.5.} \textit{There is a one to one correspondence between regular harmonic structures $(\{D_{\alpha_i}\}_{i=1}^M, \bm{r})$ and pairs $(\{\tilde{D}_{\alpha_i}\}_{i=1}^M,\bm{r})$ such that }
\[(V_{\alpha_i},D_{\alpha_i})\leq(\tilde{V}_{\alpha_i},\tilde{D}_{\alpha_i}), \quad \forall 1\leq i\leq M.\]

So we instead  study the renormalization equations $(5.3)$ on $\tilde{V}_{\alpha_i}$'s. By Lemma 5.3 and Lemma 5.4, for $p\neq q$ in $\tilde{V}_{\alpha_i}$, $(\tilde{D}_{\alpha_i})_{pq}> 0$ if and only if one of $p,q$ is $q_5$. See Figure 5.1 for an illustration, where we connect an edge between vertices $p,q$ when $(\tilde{D}_{\alpha_i})_{pq}> 0$. 

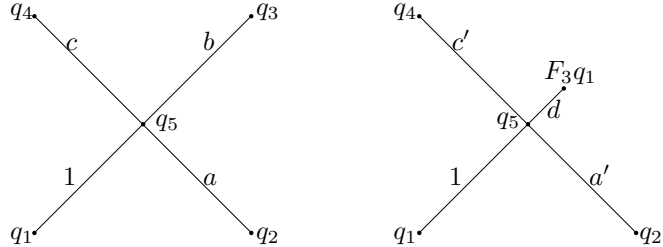
\begin{figure}[h]
	\centering
	\includegraphics[height=3cm]{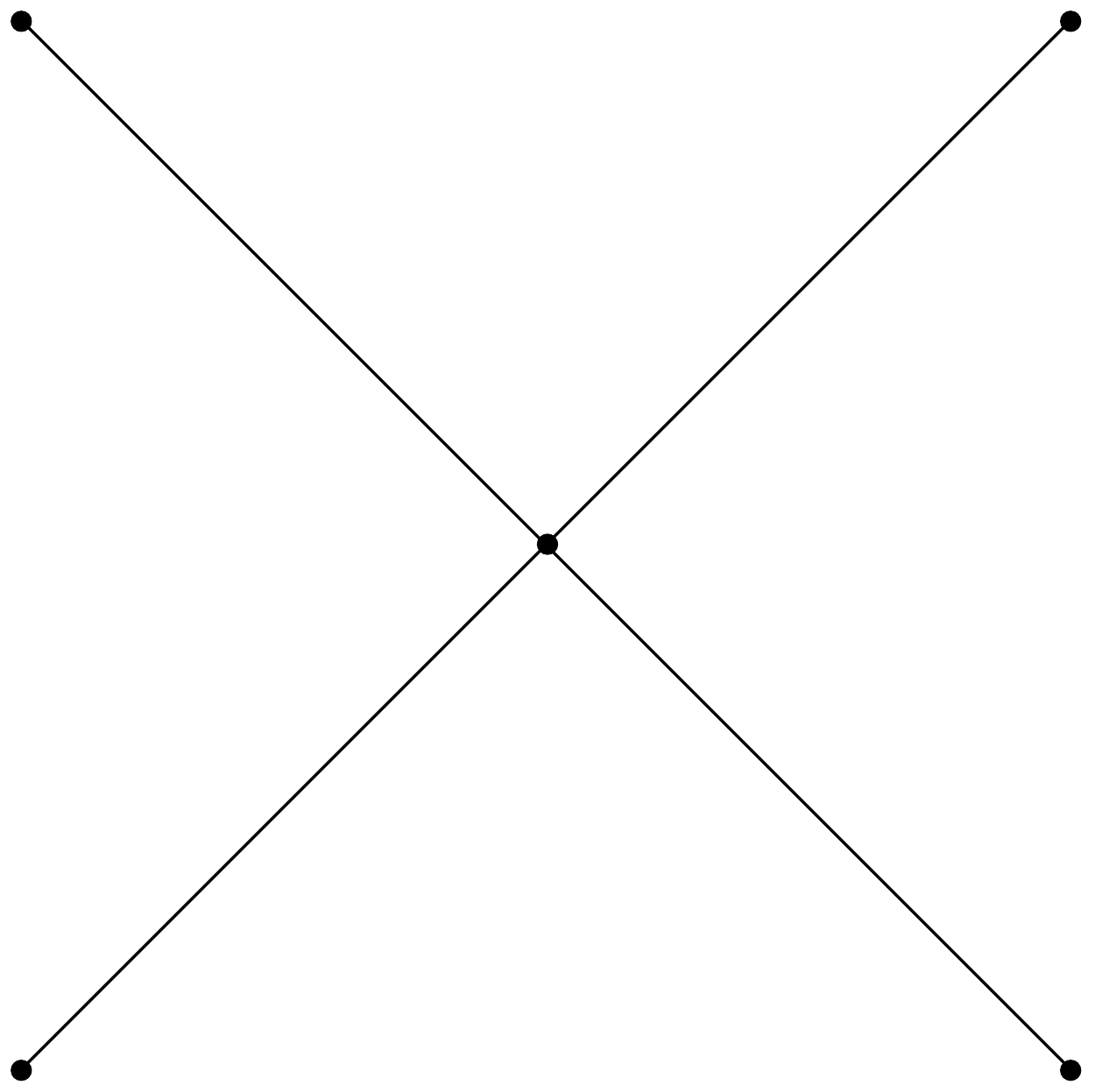}\qquad\qquad\qquad
	\includegraphics[height=3cm]{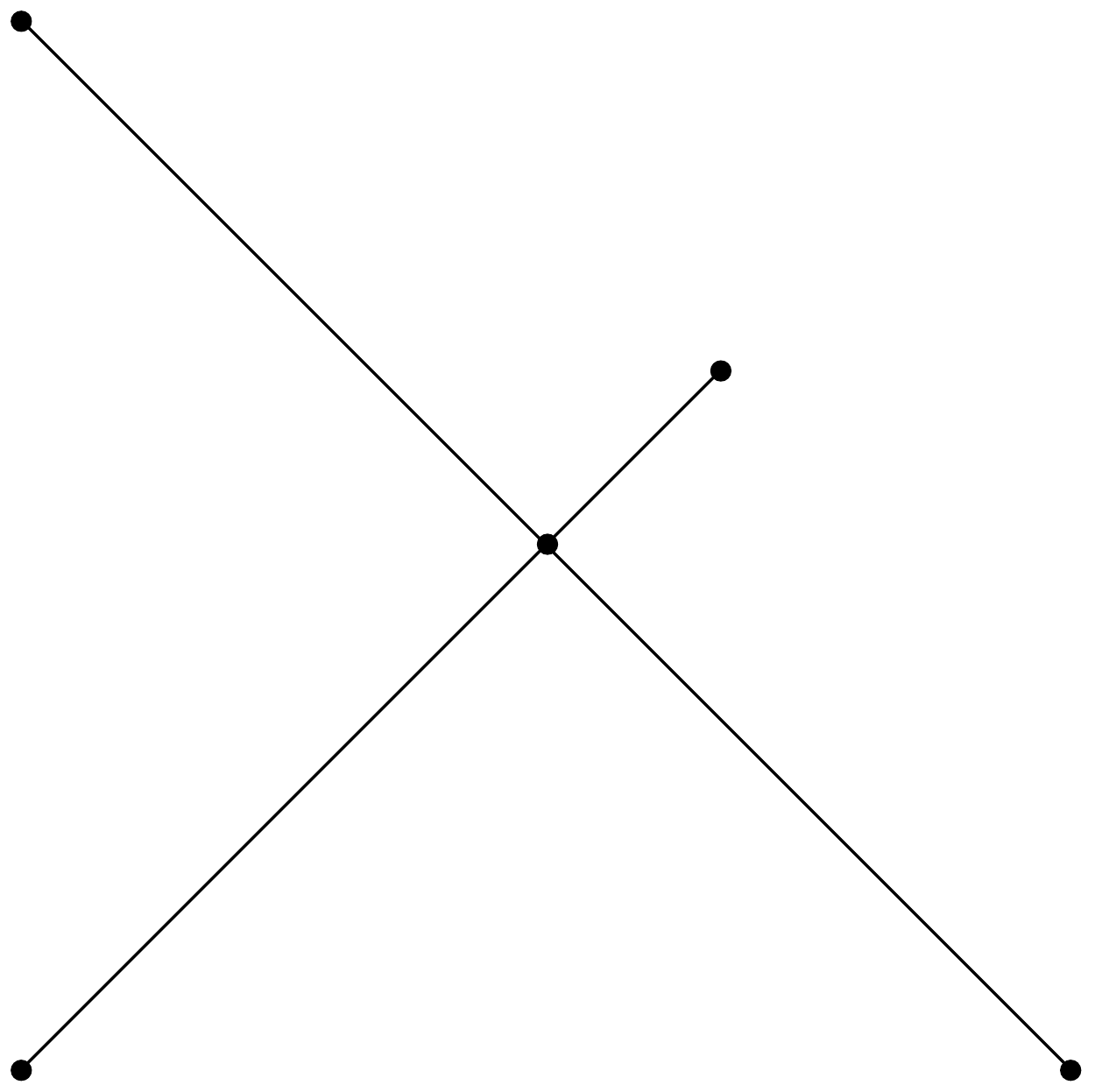}
	\begin{picture}(0,0) \thicklines 
	\put(-242,-0.5){$q_1$}
	\put(-149,-0.5){$q_2$}
	\put(-242,84){$q_4$}
	\put(-149,84){$q_3$}
	\put(-187,42){$q_5$}
	
	\put(-222,20){$1$}
	\put(-169,20){$a$}
	\put(-169,71){$b$}
	\put(-221,71){$c$}
	
	\put(-97,-0.5){$q_1$}
    \put(-4,-0.5){$q_2$}
    \put(-97,84){$q_4$}
    \put(-40,60){$F_3q_1$}
    \put(-58,43){$q_5$}
    
	\put(-76,20){$1$}
    \put(-23,20){$a'$}
    \put(-75,71){$c'$}
    \put(-39,45){$d$}    
	\end{picture}
	\begin{center}
	\caption{The resistance network $(\tilde{V}_{\alpha_i}, \tilde{D}_{\alpha_i})$, $i=1,2$.}
	\end{center}
\end{figure}

To simplify the notations, we write $r^{(i)}_{p,q}=1/(\tilde{D}_{\alpha_i})_{pq}$ the \textit{resistance} between $p,q$ in the network $(\tilde{V}_{\alpha_i}, \tilde{D}_{\alpha_i})$. Without loss of generality, we denote
\[r^{(1)}_{q_1,q_5}=1,\quad r^{(1)}_{q_2,q_5}=a, \quad r^{(1)}_{q_3,q_5}=b,\quad r^{(1)}_{q_4,q_5}=c,\]
and 
\[r^{(2)}_{q_1,q_5}=1,\quad r^{(2)}_{q_2,q_5}=a', \quad r^{(2)}_{F_3q_1,q_5}=d,\quad r^{(2)}_{q_4,q_5}=c'.\]
For convenience of readers, we mark them in Figure 5.1. As we have already labelled the similitudes $\{\psi_i\}_{i=1}^6$ by numbers $1,2,\cdots, 6$, we write $r_1, r_2,\cdots,r_6$ for their associated renormalization factors. See Figure 5.2 for an illustration.

\begin{figure}[h]
	\centering
	\includegraphics[height=3cm]{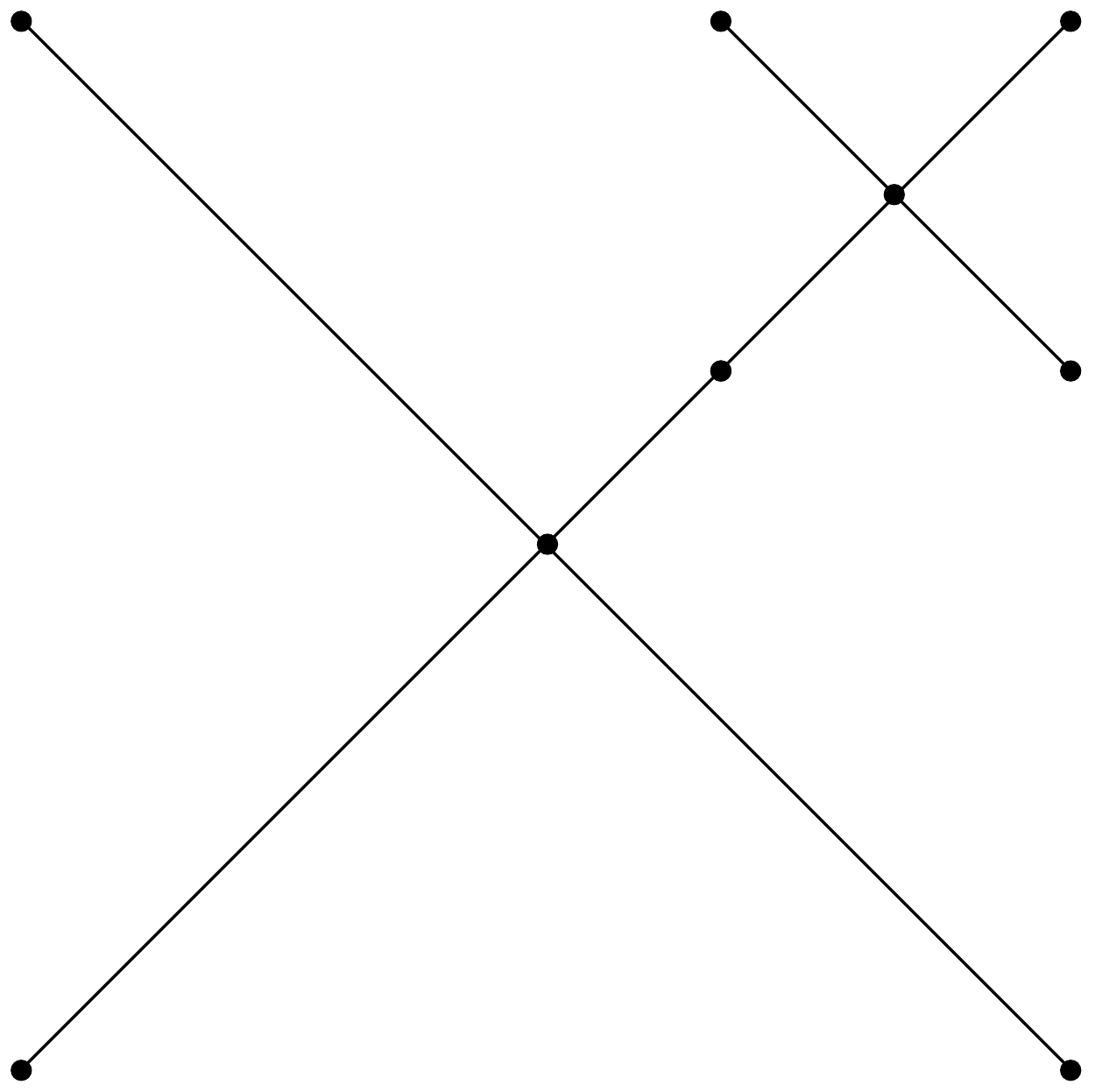}\qquad\qquad\qquad
	\includegraphics[height=3cm]{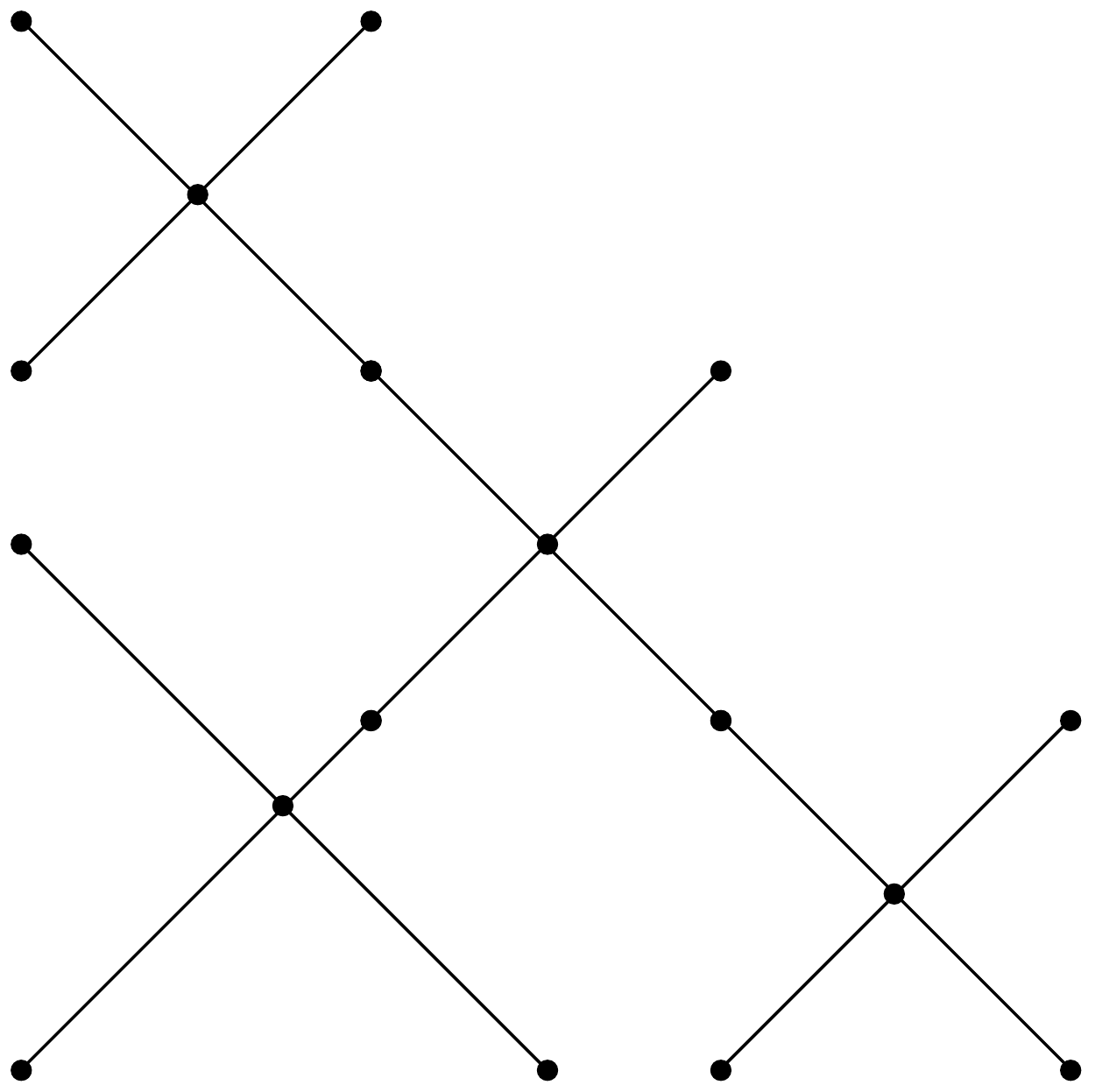}
	\begin{picture}(0,0) \thicklines 
	\put(-222,22){$r_1$}
	\put(-155,74){$r_2$}
	\put(-72,10){$r_3$}
	\put(-10,20){$r_4$}
	\put(-35,45){$r_6$}
	\put(-65,70){$r_5$}
	\end{picture}
	\begin{center}
		\caption{The renormalization factors $r_1,r_2,\cdots,r_6$.}
	\end{center}
\end{figure}
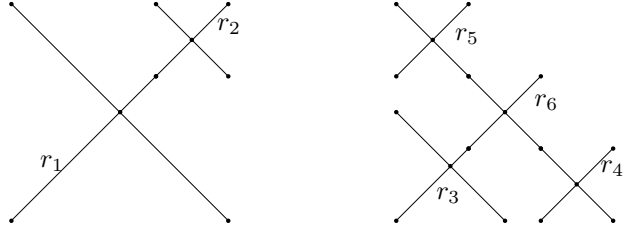

By easy operation on resistors in series, the renormalization equations are equivalent to the following equations,
\begin{equation}
\begin{cases}
a=a',c=c',r_1=1,r_6b=d,\\
r_2(1+b)+d=b,\\
r_3(1+d)+r_6=1,\\
r_4(a+c)+r_6a=a,\\
r_5(a+c)+r_6c=c.
\end{cases}
\end{equation}
It is easy to solve the equations to get that
\begin{equation}
r_1=1,\text{ }r_2=\frac{b-d}{1+b},\text{ }r_3=\frac{b-d}{b+bd},\text{ }r_4=\frac{ab-ad}{ab+cb},\text{ }r_5=\frac{cb-cd}{ab+cb},\text{ }r_6=\frac{d}{b}.
\end{equation}
Obviously, to get a regular harmonic structure, we only need to assume $a,b,c,d>0$ and $b>d$. This gives all the regular harmonic structures associated with the given $f.r.f.t.$ nested structure. The solutions depend on $4$ parameters.

Since the fractal $\mathcal{V}^o$ possesses obvious geometric symmetry and local symmetries, we may demand the resulting Dirichlet forms possess the same symmetries. We just need to require  $a=c$, and this will give all the \textit{symmetric} regular harmonic structures. So the symmetric solutions depend on $3$ parameters. Furthermore, it is also reasonable to demand that the islands of same size have the same energy. To be precise, for any two islands $K_\alpha$, $K_\beta$ with the same type and same size, i.e., $\phi_{\alpha,\beta}$ is an isometric mapping, we want 
\[\mathcal{E}(u)=\mathcal{E}(u\circ \phi_{\alpha,\beta})\]
for any function $u$ supported in $K_\beta$. We need to require that $r_2=r_4=r_5=r_6$ in addition. This will gives  $b=1,d=\frac{1}{3}$ along with $a=c$.  So the solutions depend on only $1$ parameters.  We call such solutions \textit{homogeneous} harmonic structures, see further discussion in Section 6.

We summarize what we have got in the following.

\textbf{Theorem 5.6.} \textit{For the $f.r.f.t.$ nested structure associated with (4.2) of $\mathcal{V}^o$, the full solution of the regular harmonic structures is as shown in (5.5). It depends on $4$ parameters. The solution of symmetric ones depends on $3$ parameters and the solution of homogeneous ones depends on $1$ parameters.}

Before ending this example, we would like to briefly discuss the harmonic structures of the other $f.r.f.t.$ nested structure of $\mathcal{V}^o$ mentioned in Section 4, see Figure 4.4. Analogously, we define $\tilde{V}_{\alpha_1}:=V_{\alpha_1}\cup\{q_5\}=\{q_1,q_2,q_3,q_4,q_5\}$ and $\tilde{V}_{\alpha_2}:=V_{\alpha_2}\cup\{q_5\}=\{q_1,F_2q_4,F_3q_1,F_4q_2,q_5\}$. See Figure 5.3 for the resistance network of $(\tilde{V}_{\alpha_i}, \tilde{D}_{\alpha_i})$ with parameters $a,b,c,d,e,f$ marked there.

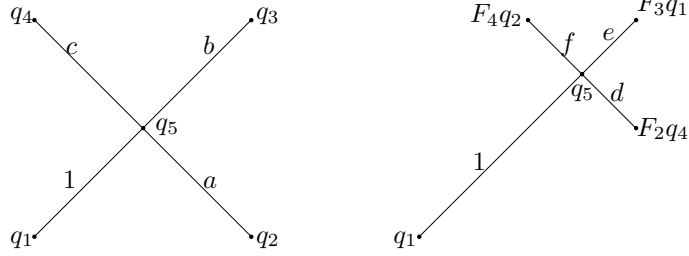
\begin{figure}[h]
	\centering
	\includegraphics[height=3cm]{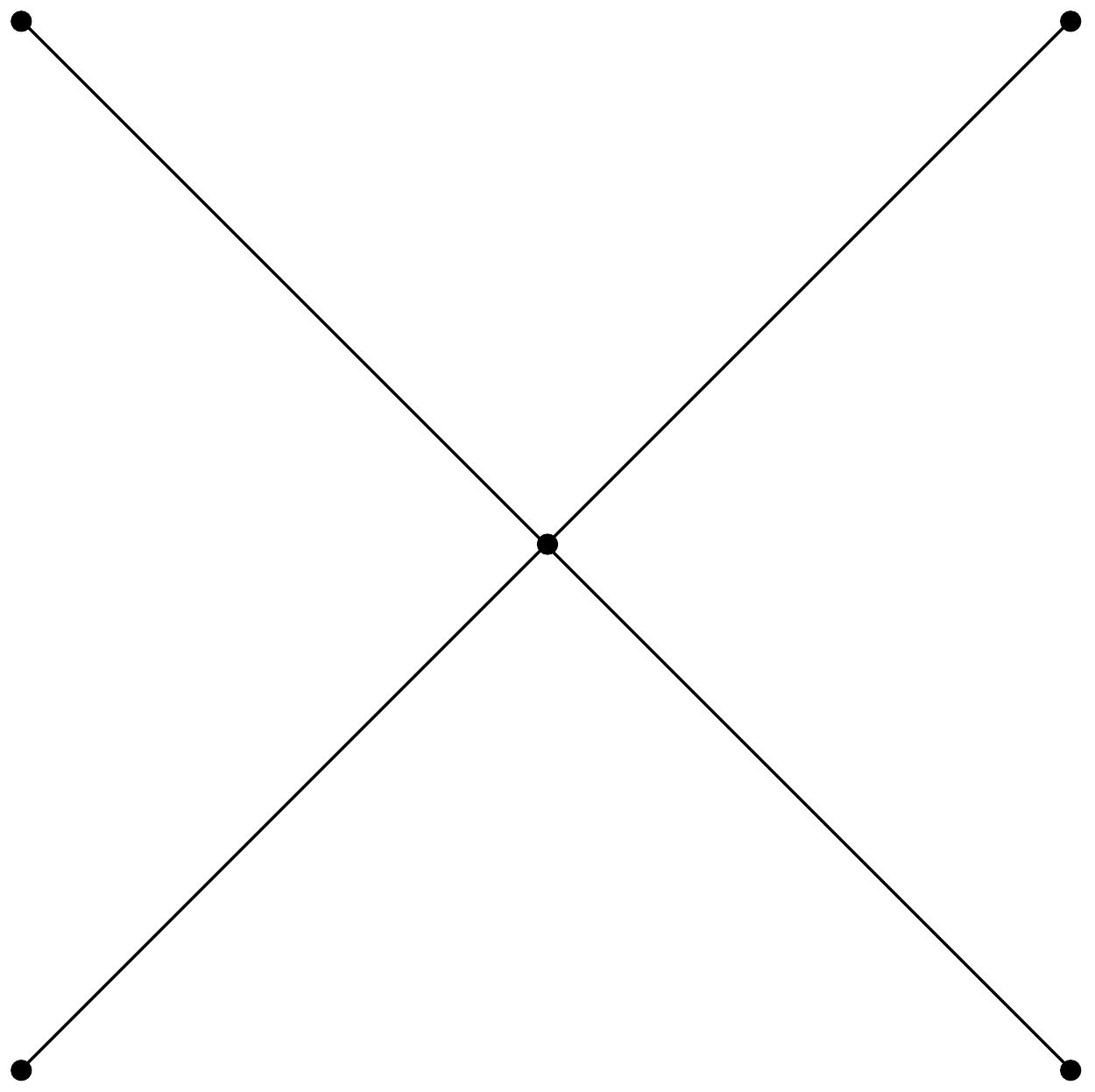}\qquad\qquad\qquad
	\includegraphics[height=3cm]{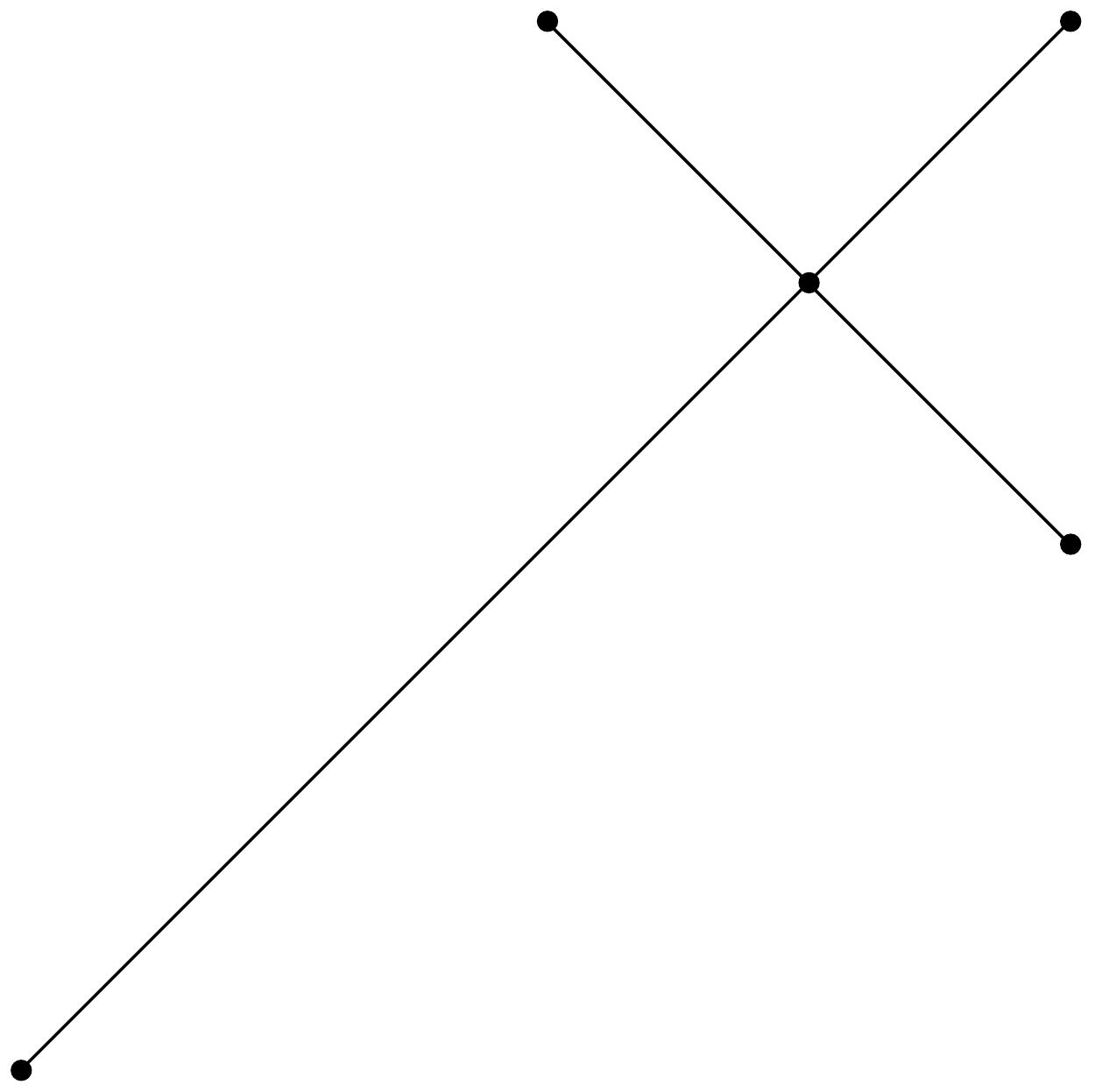}
	\begin{picture}(0,0) \thicklines 
	\put(-242,-0.5){$q_1$}
	\put(-149,-0.5){$q_2$}
	\put(-242,84){$q_4$}
	\put(-149,84){$q_3$}
	\put(-187,42){$q_5$}
	
	\put(-222,20){$1$}
	\put(-169,20){$a$}
	\put(-169,71){$b$}
	\put(-221,71){$c$}
	
	\put(-97,-0.5){$q_1$}
    \put(-5,40){$F_2q_4$}
    \put(-68,83){$F_4q_2$}
    \put(-5,87){$F_3q_1$}
    \put(-30,55){$q_5$}
    
    \put(-67,27){$1$}
    \put(-15,53){$d$}
    \put(-18,76){$e$}
    \put(-34,72){$f$}    
	\end{picture}
	\begin{center}
		\caption{The resistance network $(\tilde{V}_{\alpha_i}, \tilde{D}_{\alpha_i})$, $i=1,2$.}
	\end{center}
\end{figure}

We only give the associated homogeneous regular harmonic structures. For this purpose,  some renormalization factors are required to coincide. After some simplification, we can use two parameters $r,s$ to present all the renormalization factors, see Figure 5.4.

\begin{figure}[h]
	\centering
	\includegraphics[height=3cm]{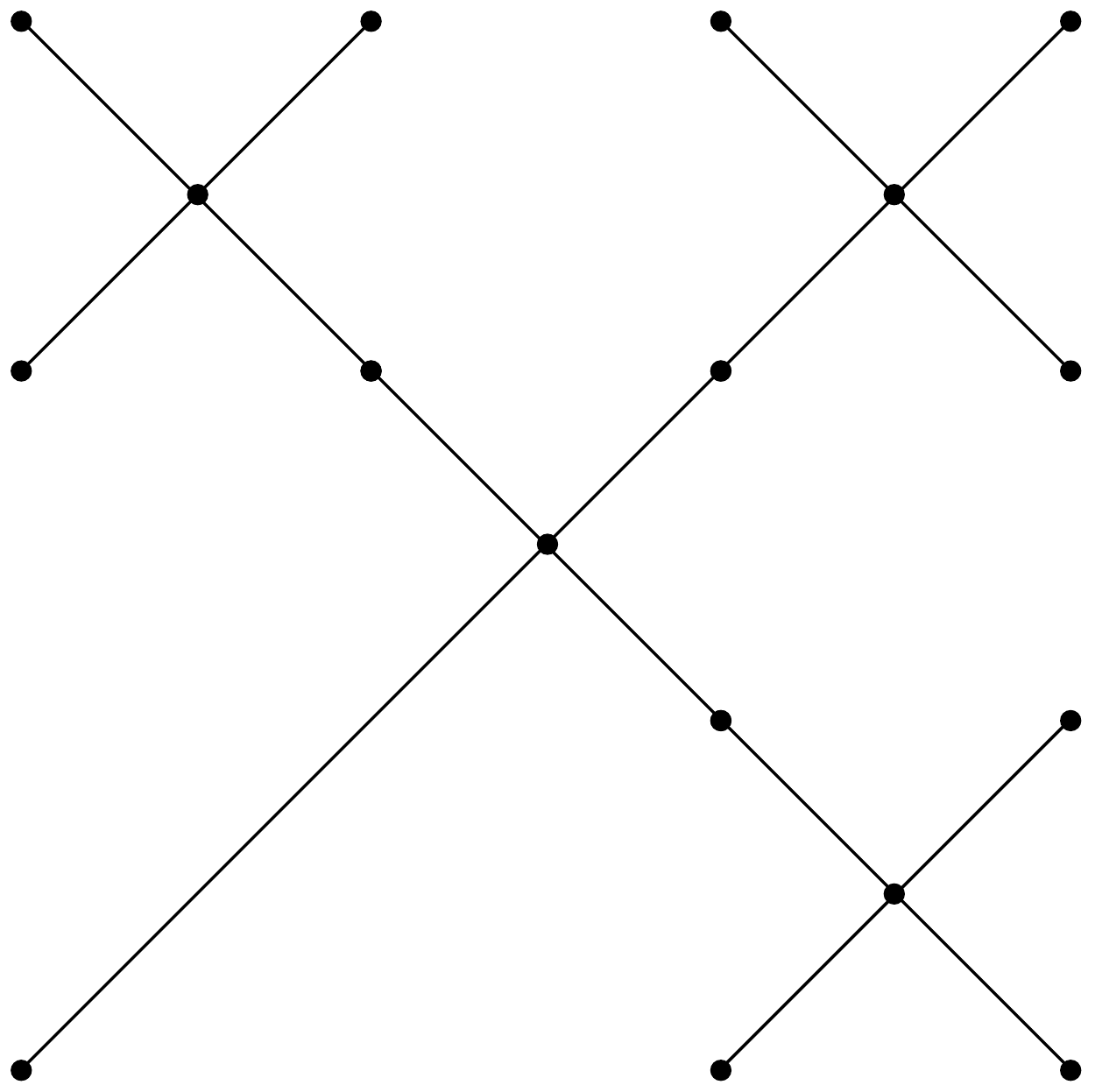}\qquad\qquad\qquad
	\includegraphics[height=3cm]{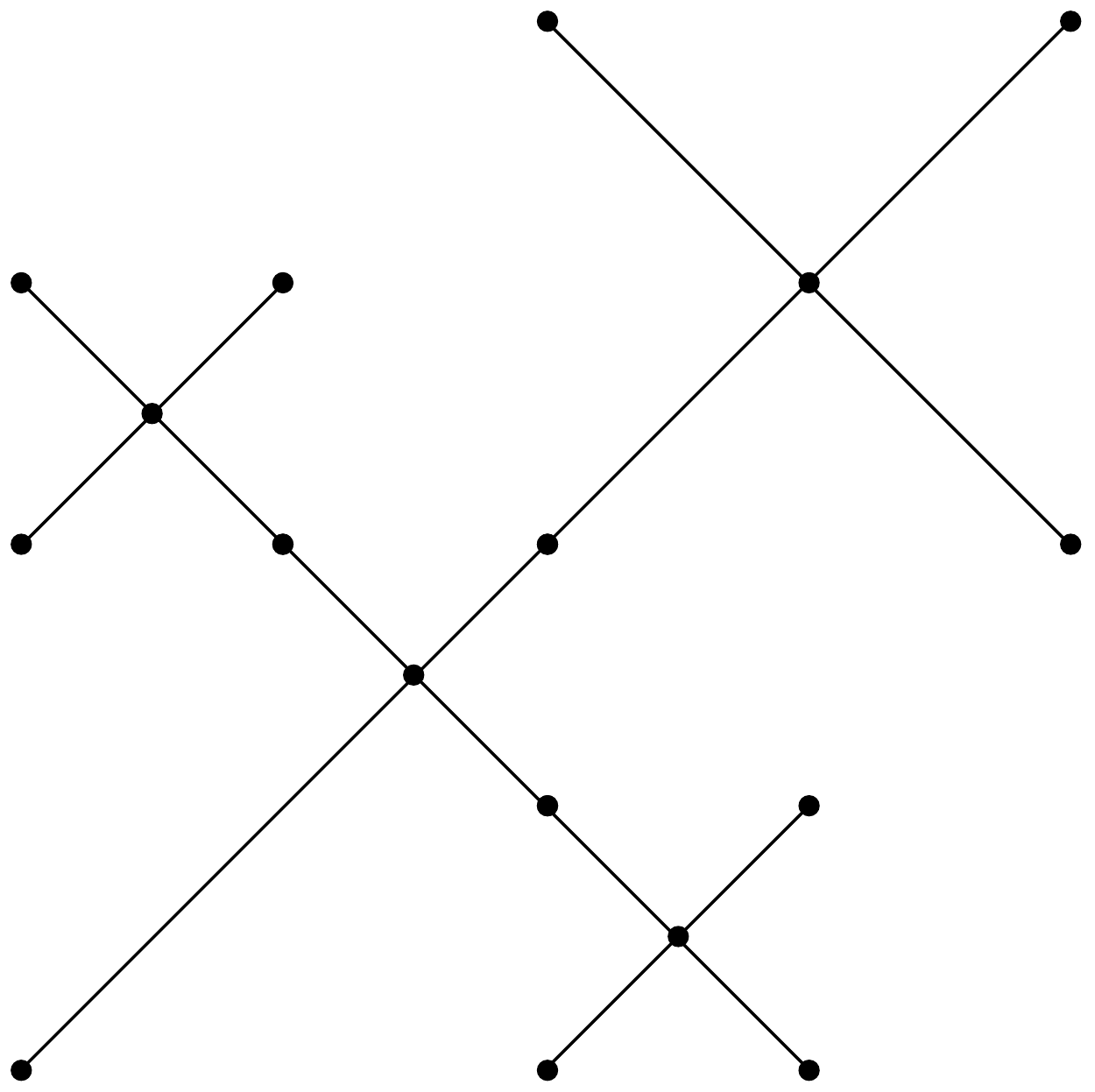}
	\begin{picture}(0,0) \thicklines 	
	\put(-222,20){$1$}
	\put(-169,20){$r$}
	\put(-165,74){$r$}
	\put(-220,74){$r$}
	
	\put(-27,67){$r$}
	\put(-72,22){$s$}
	\put(-81,57){$rs$}
	\put(-40,17){$rs$}

	\end{picture}
	\begin{center}
		\caption{The renormalization factors.}
	\end{center}
\end{figure}
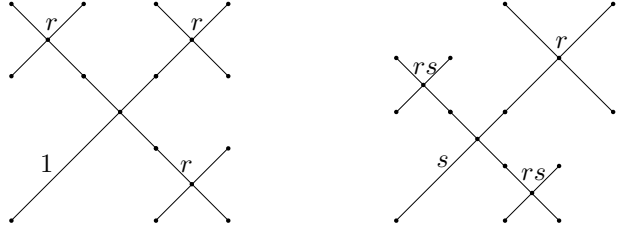
The solution is $r=\frac{1}{3},s=\frac{1}{2},a=c, b=1, e=\frac{1}{3}, d=f=\frac{1}{3}a$, depending on parameter $a$. This gives the same homogeneous regular harmonic structures as that of the  $f.r.f.t.$ nested  structure given by (4.2).

\subsection{Overlapping gasket with open bottom}In this example $q_4$ is a cut-point by a same argument as in the previous example. For the uniformity of notations, we take 
\[\psi_1=F_2,\psi_2=id,\psi_3=F_1,\psi_4=F_3,\psi_5=F_5,\psi_6=F_4,\]
then we have the following  graph-directed construction with edge set $\{e_i\}_{i=1}^6$,
\[\begin{cases}
K_1=\psi_1K_1\cup \psi_2K_2,\\
K_2=\psi_3K_1\cup \psi_4K_1\cup \psi_5K_1\cup\psi_6K_2.
\end{cases}\]
See Figure 5.5 for an illustration.
\begin{figure}[h]
	\centering
	\includegraphics[height=3.5cm]{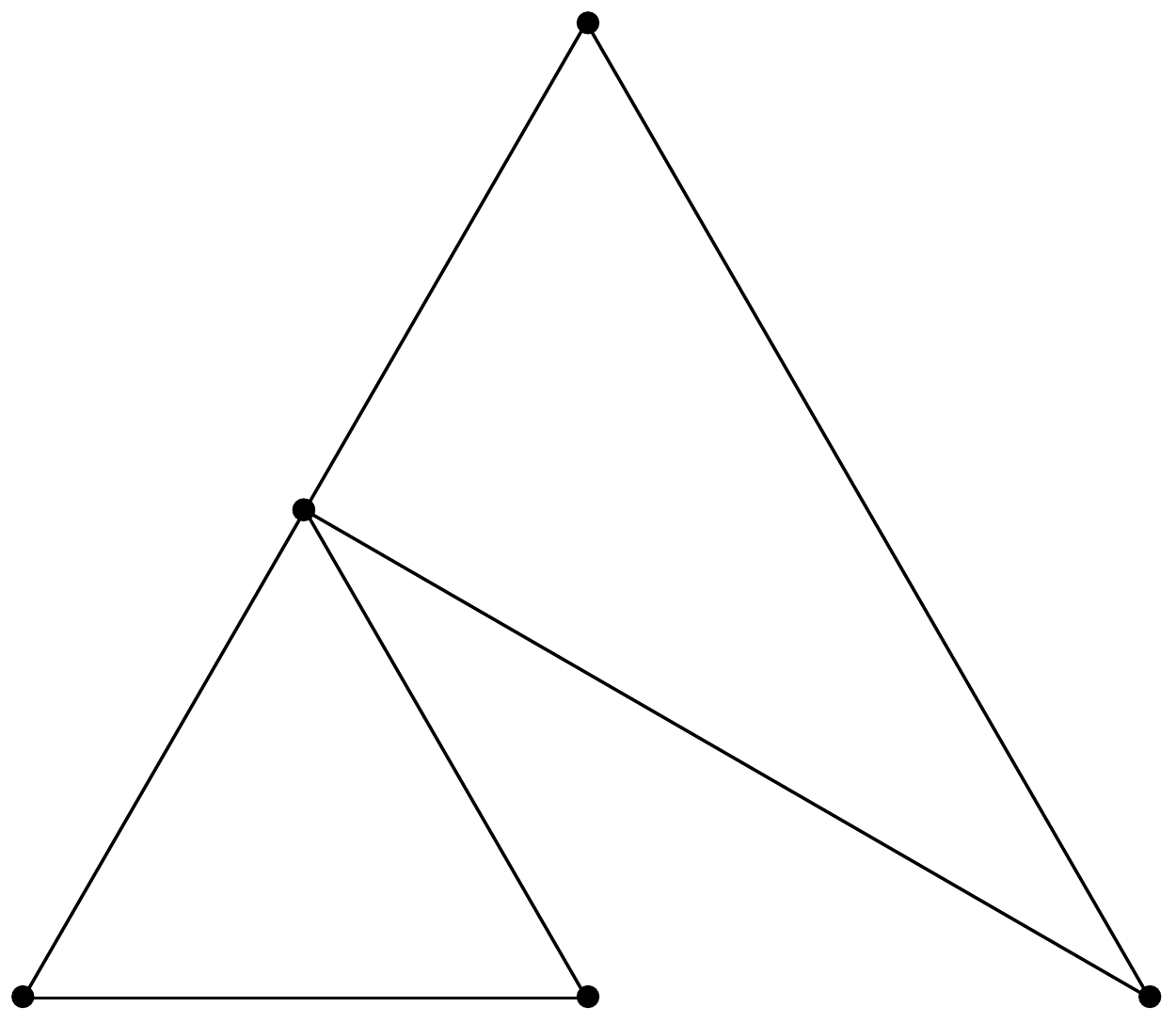}\qquad\qquad\qquad
	\includegraphics[height=3.5cm]{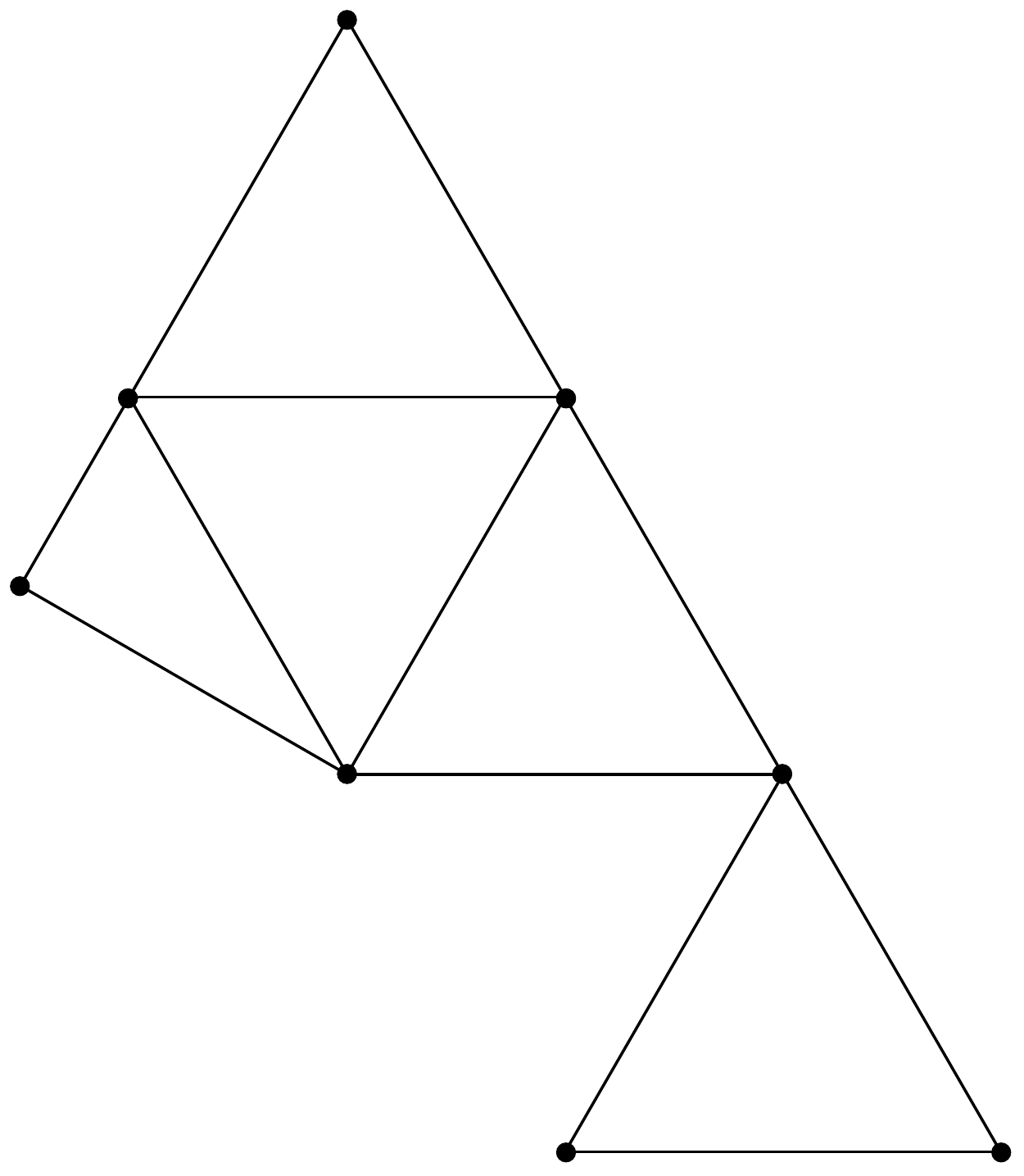}
	\begin{picture}(0,0) \thicklines 
	\put(-247,15){$\psi_1K_1$}
	\put(-218,55){$\psi_2K_2$}
	\put(-88,47){$\psi_6K_2$}
	\put(-34,6){$\psi_4K_1$}
	\put(-52,39){$\psi_5K_1$}
	\put(-71,70){$\psi_3K_1$}
	\end{picture}
	\begin{center}
		\caption{The graph-directed construction of $\mathcal{SG}^o$.}
	\end{center}
\end{figure}
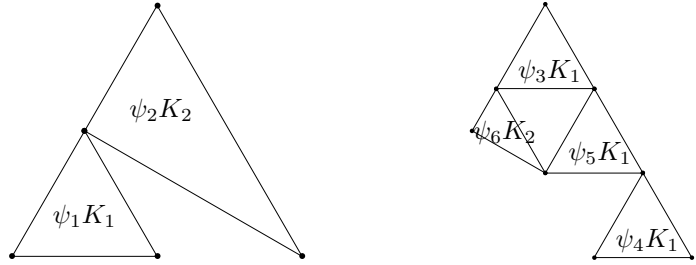
By Theorem 3.4 we get a $f.r.f.t.$ nested structure with $K_{\alpha_1}=K_1=\mathcal{SG}^o, K_{\alpha_2}=K_2$. For simplicity, take 
\[r^{(1)}_{q_1,q_2}=a,r^{(1)}_{q_1,q_3}=b,r^{(1)}_{q_2,q_3}=c, r^{(2)}_{q_1,q_4}=d, r^{(2)}_{q_1,q_3}=e,r^{(2)}_{q_4,q_3}=f,\]
where $r^{(i)}_{p,q}$'s are the resistances in the resistance networks $(V_{\alpha_i},D_{\alpha_i}), i=1,2$, see Figure 5.6.
\begin{figure}[h]
	\centering
	\includegraphics[height=3.5cm]{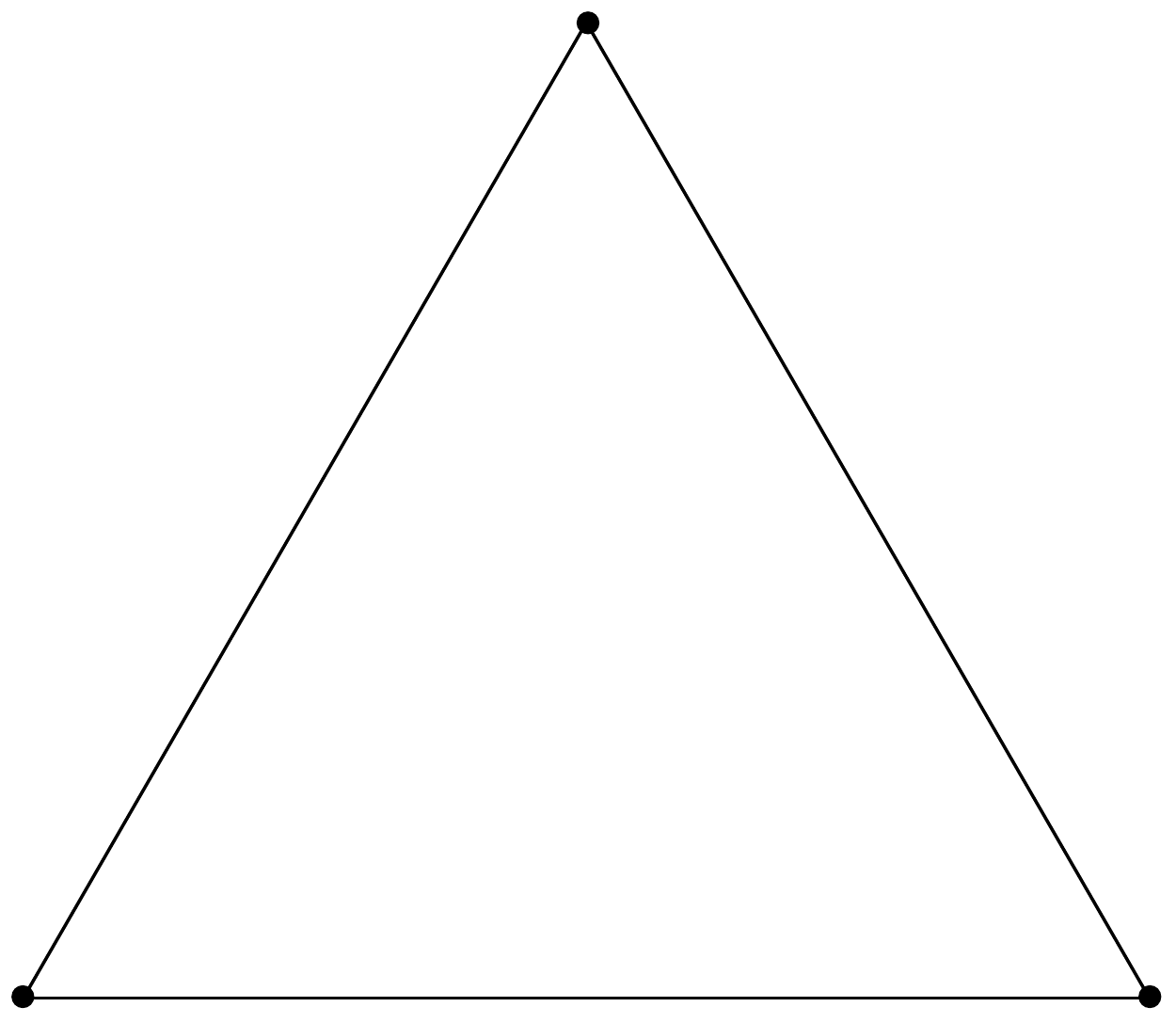}\qquad\qquad\qquad
	\includegraphics[height=3.5cm]{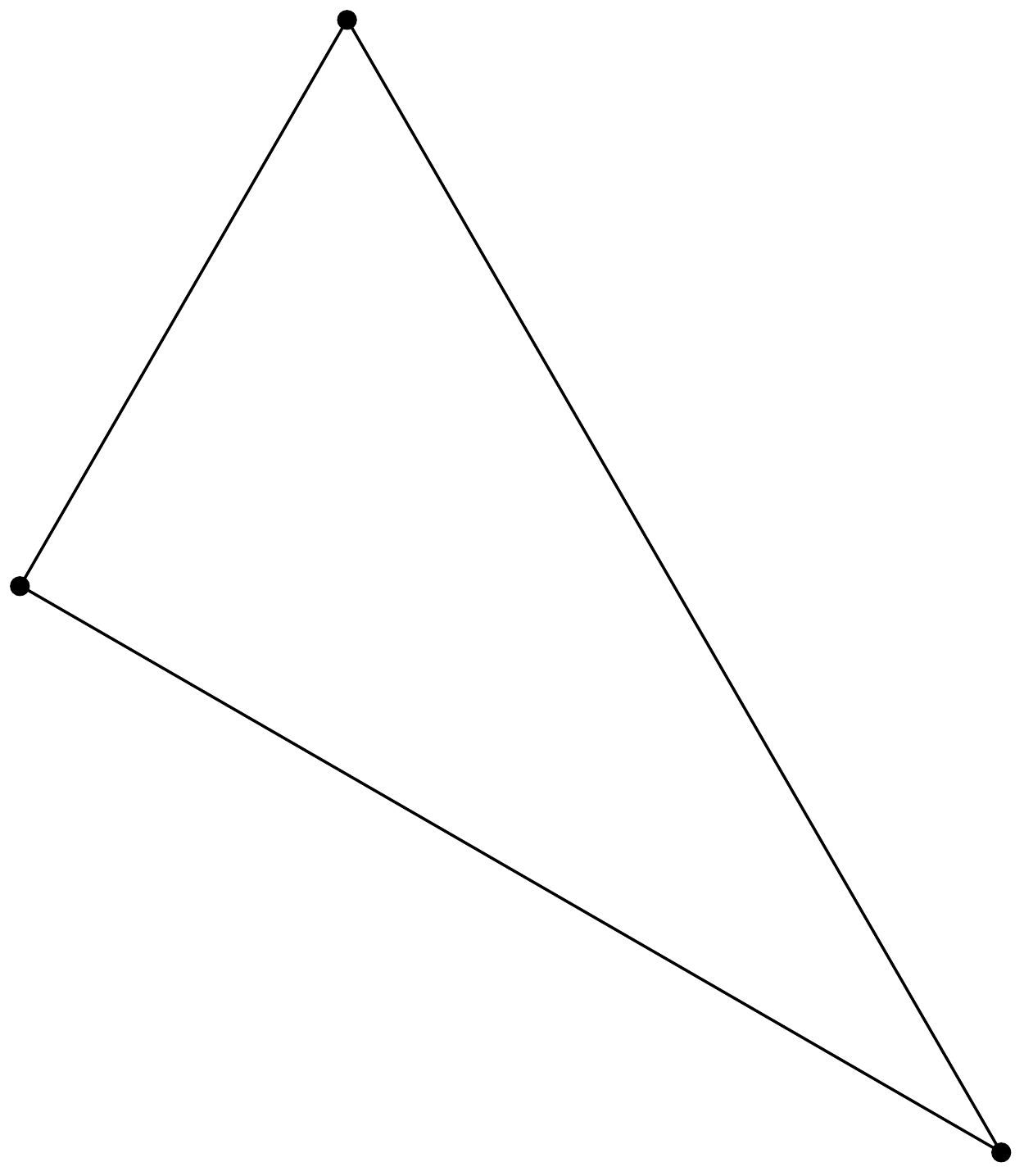}
	\begin{picture}(0,0) \thicklines 
	\put(-272,0){$q_2$}
	\put(-212,101){$q_1$}
	\put(-151,0){$q_3$}
	\put(-240,50){$a$}	
	\put(-180,50){$b$}
	\put(-210,5){$c$}
	
	\put(-98,49){$q_4$}
	\put(-3,0){$q_3$}
	\put(-65,101){$q_1$}
	\put(-80,73){$d$}
	\put(-30,46){$e$}
	\put(-60,24){$f$}
	\end{picture}
	\begin{center}
		\caption{The resistance networks $(V_{\alpha_i},D_{\alpha_i}),i=1,2$.}
	\end{center}
\end{figure}

For the renormalization factors, we used the simplified notations $r_i$ instead of $r_{e_i}$ as in the previous example, with $r_2=1$ as $\psi_2=id$. We focus on the homogeneous regular harmonic structures. It is easy to verify that it is necessary to require $r_3=r_4=r_5=r_6$, and we use the symbol $s$ to denote  them  hereafter.

It is convenient to use the $\Delta-Y$ transformation for resistance networks here, see Figure $5.7$ for an illustration of the $\Delta-Y$ transformation. See Figure 5.8 for the transformations between the first two level resistance networks.

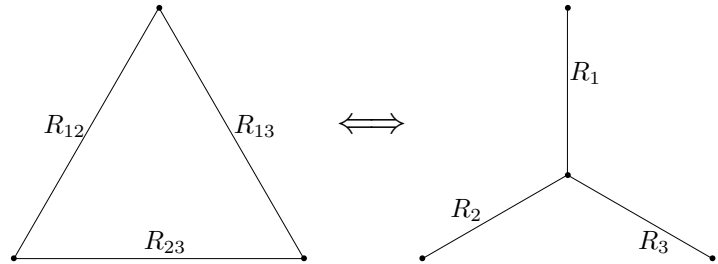
\begin{figure}[h]
	\centering
	\includegraphics[height=3.5cm]{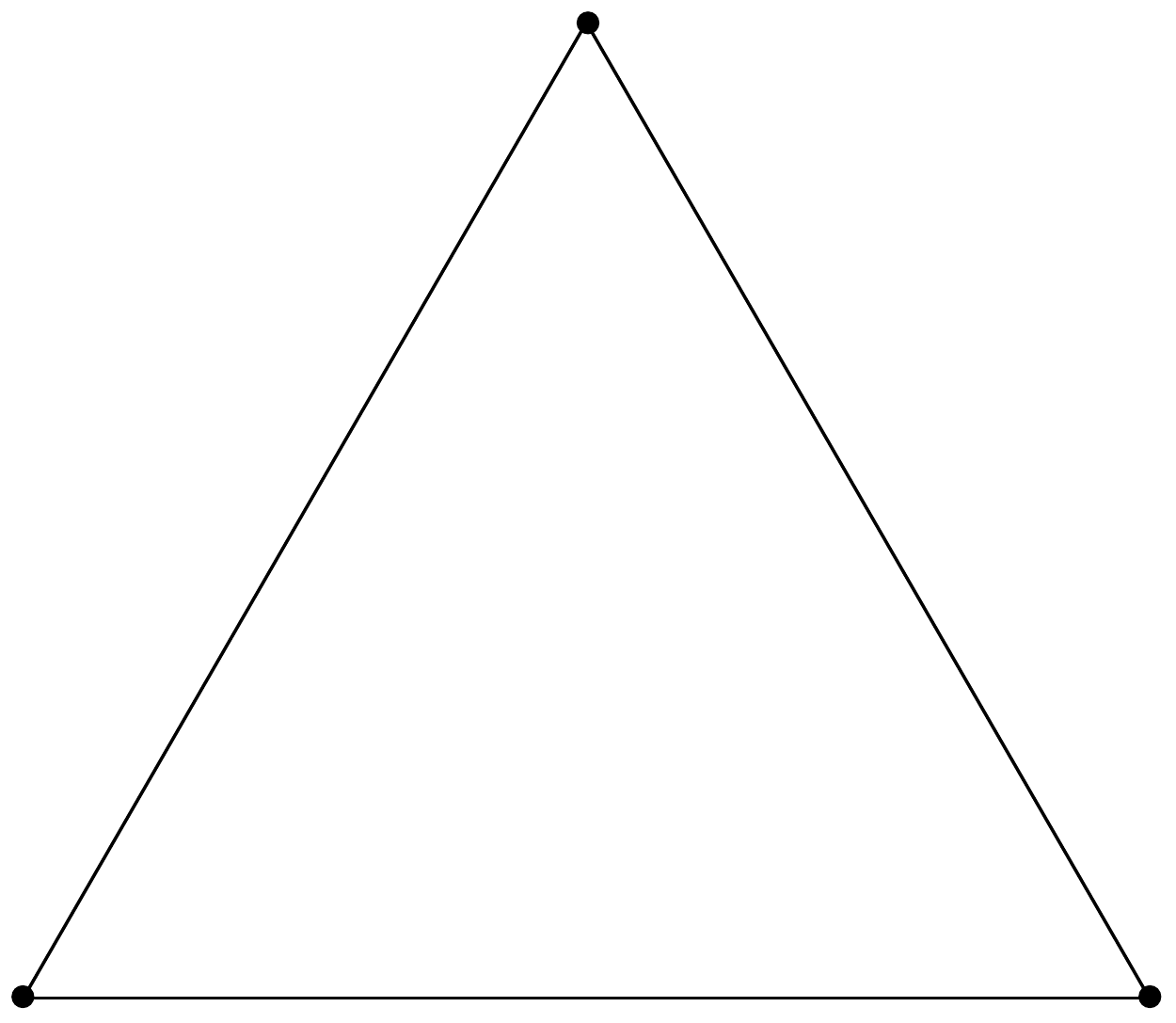}\qquad\qquad
	\includegraphics[height=3.5cm]{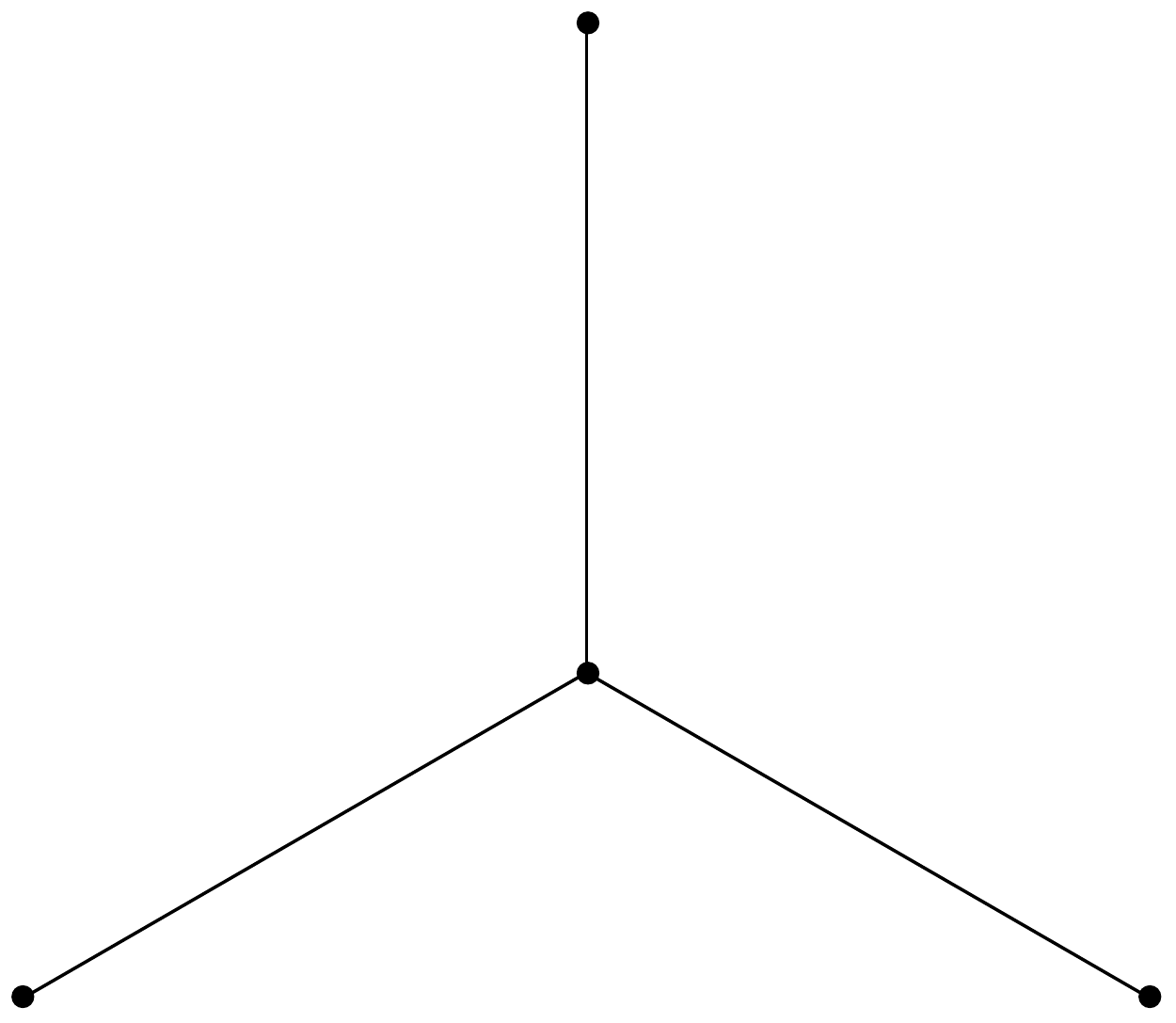}
	\begin{picture}(0,0) \thicklines 
	\put(-259,50){$R_{12}$}	
	\put(-187,50){$R_{13}$}
	\put(-221,6){$R_{23}$}
	\put(-34,6){$R_3$}
	\put(-105,18){$R_2$}
	\put(-60,70){$R_1$}
	\put(-148,50){\LARGE{$\Longleftrightarrow$}}
	\end{picture}
	\begin{center}
		\caption{An illustration of the  $\Delta-Y$ transformation, with $R_i=\frac{R_{ij}R_{ik}}{R_{12}+R_{23}+R_{13}}$, $\{i,j,k\}=\{1,2,3\}$. All $R_i$'s and $R_{ij}$'s are resistances. }
	\end{center}
\end{figure}

\begin{figure}[h]
	\centering
	\includegraphics[height=3cm]{graphso01.pdf}\qquad
	\includegraphics[height=3cm]{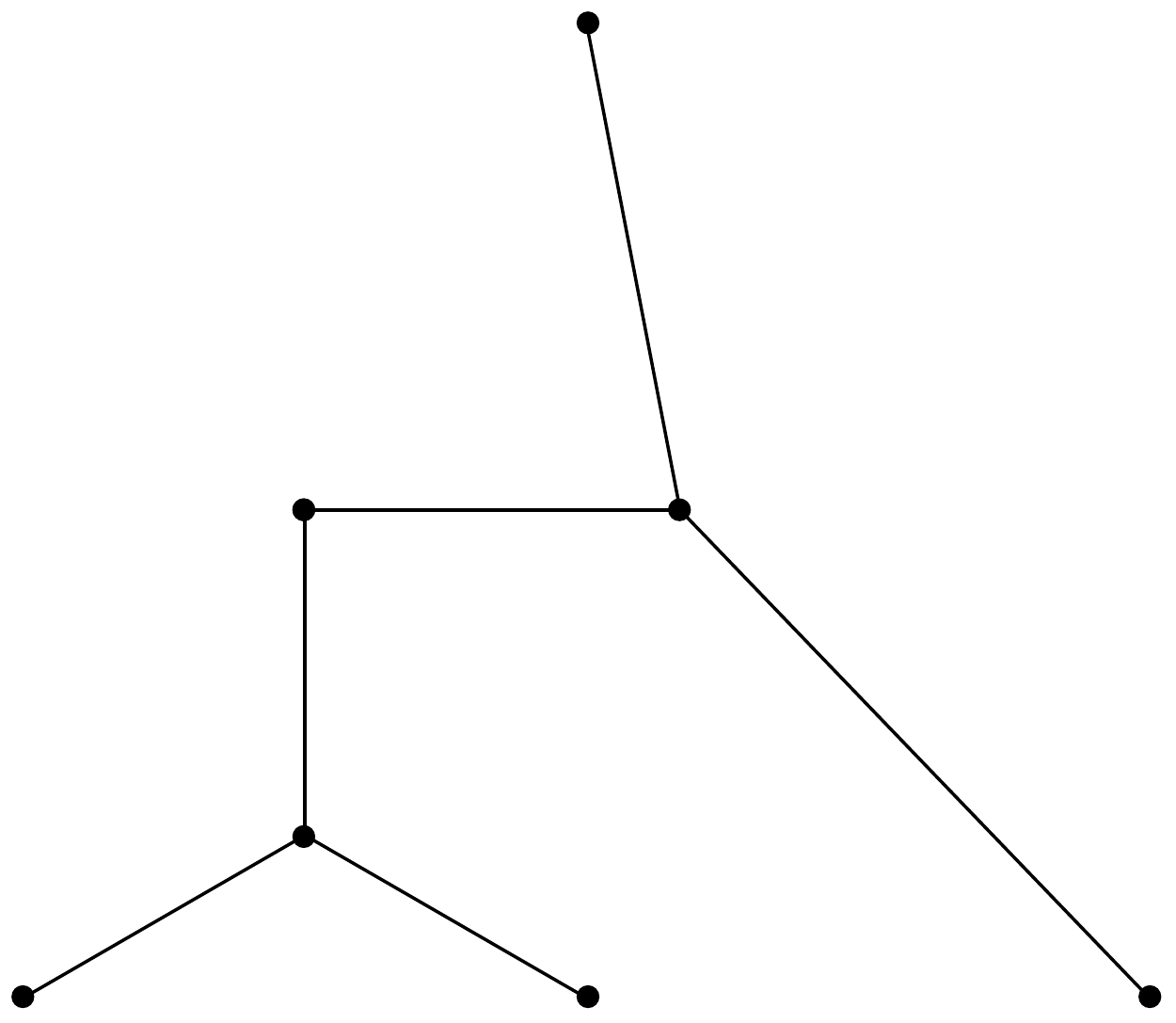}\qquad
	\includegraphics[height=3cm]{y.pdf}\\\vspace{0.4cm}
	\includegraphics[height=3cm]{graphso11.pdf}\qquad
	\includegraphics[height=3cm]{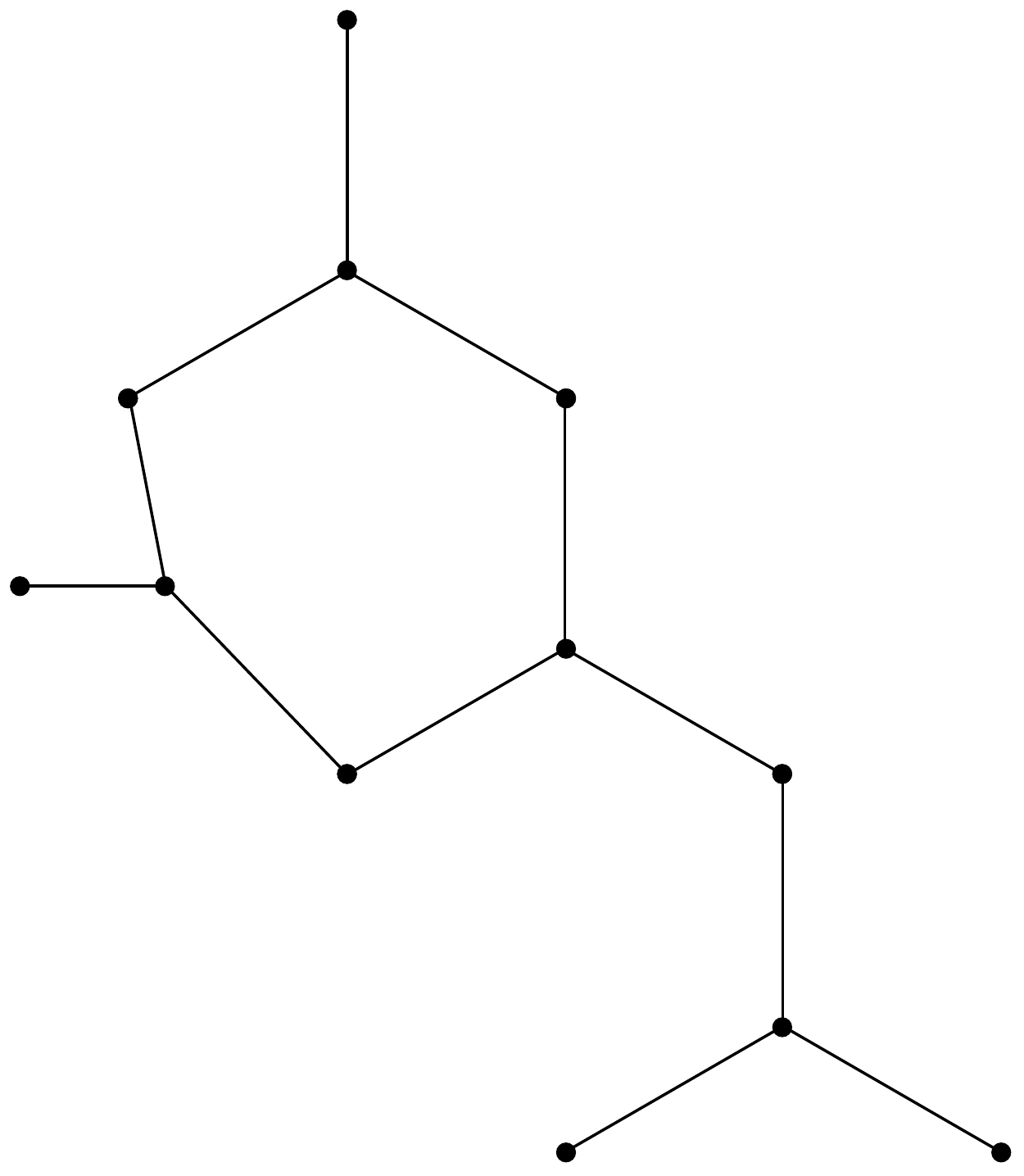}\qquad
	\includegraphics[height=3cm]{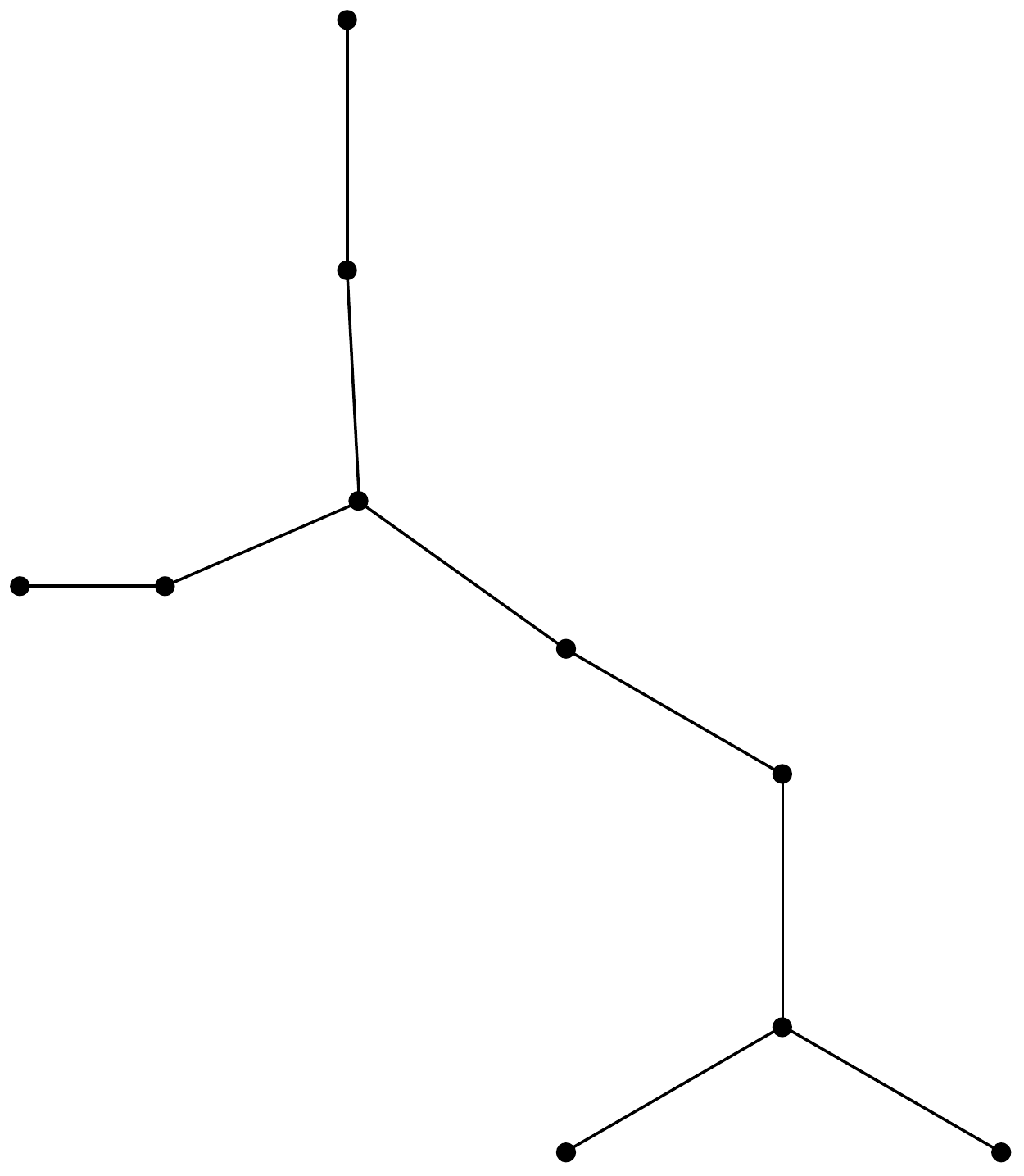}\qquad
	\includegraphics[height=3cm]{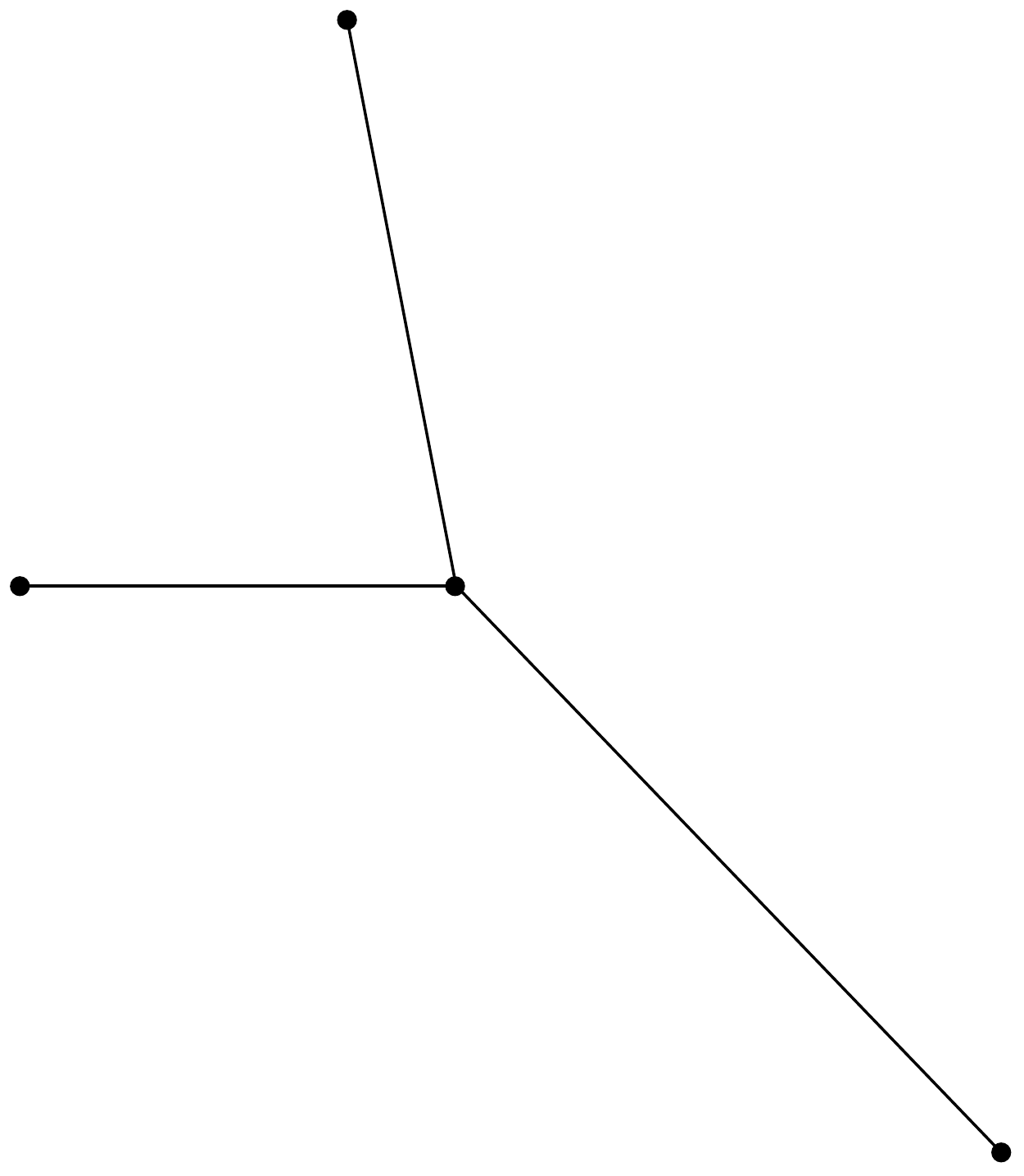}
	\begin{picture}(0,0) \thicklines 
	\put(-68,145){$\Longrightarrow$}
	\put(50,145){$\Longrightarrow$}
	\put(-15,40){$\Longrightarrow$}
	\put(-109,40){$\Longrightarrow$}
	\put(82,40){$\Longrightarrow$}
	\end{picture}
	\begin{center}
		\caption{Transformations between the first two level resistance networks.}
	\end{center}
\end{figure}

In view of Figure 5.8, it is more convenient to use  the $Y$-shaped networks for the calculation, then solve the $\Delta$-shaped networks by doing the inverse $\Delta-Y$  transformation. See 
Figure 5.9 for the $Y$-shaped networks with resistances $a', b', e', d', f'$ marked there.
\begin{figure}[h]
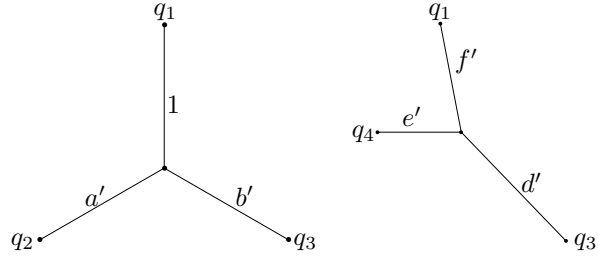

	\centering
	\includegraphics[height=3cm]{y.pdf}\qquad\quad
	\includegraphics[height=3cm]{graphso10t1.pdf}
	\begin{picture}(0,0) \thicklines
	\put(-161,87){$q_1$}
	\put(-215,0){$q_2$}
	\put(-108,0){$q_3$}
	\put(-156,50){$1$}
	\put(-130,15){$b'$}
	\put(-187,15){$a'$}
	
	\put(-57,87){$q_1$}
	\put(-86,41){$q_4$}
	\put(-2,0){$q_3$}
	\put(-47,67){$f'$}
	\put(-67,45){$e'$}
	\put(-22,20){$d'$}
	\end{picture}
	\begin{center}
		\caption{The $Y$-shaped resistance network of $(V_{\alpha_i}, D_{\alpha_i})$, i=1,2.}
	\end{center}
\end{figure}

Put the resistances and renormalization factors into the transformations shown in Figure 5.8. We get
\begin{equation}f'=1,d'=b',e'+r_1+r_1a'=a',\end{equation} 
by the transformation on $V_{\alpha_1,1}$, and
\[\begin{cases}
s+\frac{s(1+a')(1+b')}{2(1+a'+b')}=1,\\
se'+\frac{s(1+a')(a'+b')}{2(1+a'+b')}=e',\\
2sb'+s+\frac{s(1+b')(a'+b')}{2(1+a'+b')}=b',
\end{cases}\]
by the transformation on $V_{\alpha_2,1}$, using equations in (5.6). Solving these equations, we get
\begin{equation}\begin{cases}
b'=d'=\frac{2+3a'}{a'-2},e'=\frac{1}{2a'}+\frac{1}{4}+\frac{a'}{4}, f'=1,\\
s=\frac{2+a'}{4+3a'}, r_1=\frac{-2-a'+3(a')^2}{4a'+4(a')^2}.
\end{cases}
\end{equation}
The solution depends on the parameter $a'$ and gives us  the homogeneous regular harmonic structures when $a'>2$. 

So the result is 

\textbf{Theorem 5.7.} \textit{For the $f.r.f.t.$ nested structure associated with (4.3) of $\mathcal{SG}^o$, the full solution of the homogeneous regular harmonic structures is as shown in (5.7). It depends on $1$ parameters.}

\subsection{Vicsek windmill set.} The Vicsek windmill set $\mathcal{V}^w$ possesses an obvious rotational  symmetry. It is reasonable to require that the homogeneous regular harmonic structures also possess the same symmetry.

Similar as previous examples, we write $K_{\alpha_1}=K_1=\mathcal{V}^w$, $K_{\alpha_2}=K_2$, $K_{\alpha_3}=K_3$, and ${G}=(S,E)$ the directed graph associated with the $f.r.g.d.$ construction  illustrated in Figure 4.6. At the first glance, to determine a homogeneous regular harmonic structure $(\{D_{\alpha_i}\}_{i=1}^3,\bm{r})$ for the associated  $f.r.f.t.$ nested structure, we need $8$ parameters to represent the renormalization factors. For $i,j\in\{1,2,3\}$, we use $r_{ij}$ to denote the renormalization factor associated with the edge in $E$ starting from $K_{\alpha_i}$ and ending to $K_{\alpha_j}$, see Figure 5.10 for an illustration. Note that there is no $r_{13}$ since no such edge exists in $E$. 

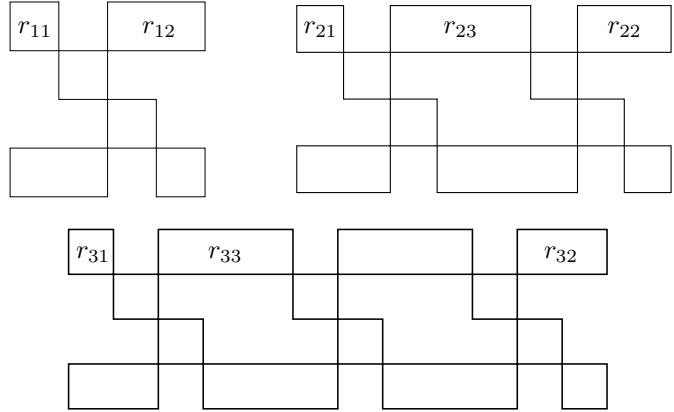
\begin{figure}[h]
	\centering
	\includegraphics[height=2.7cm]{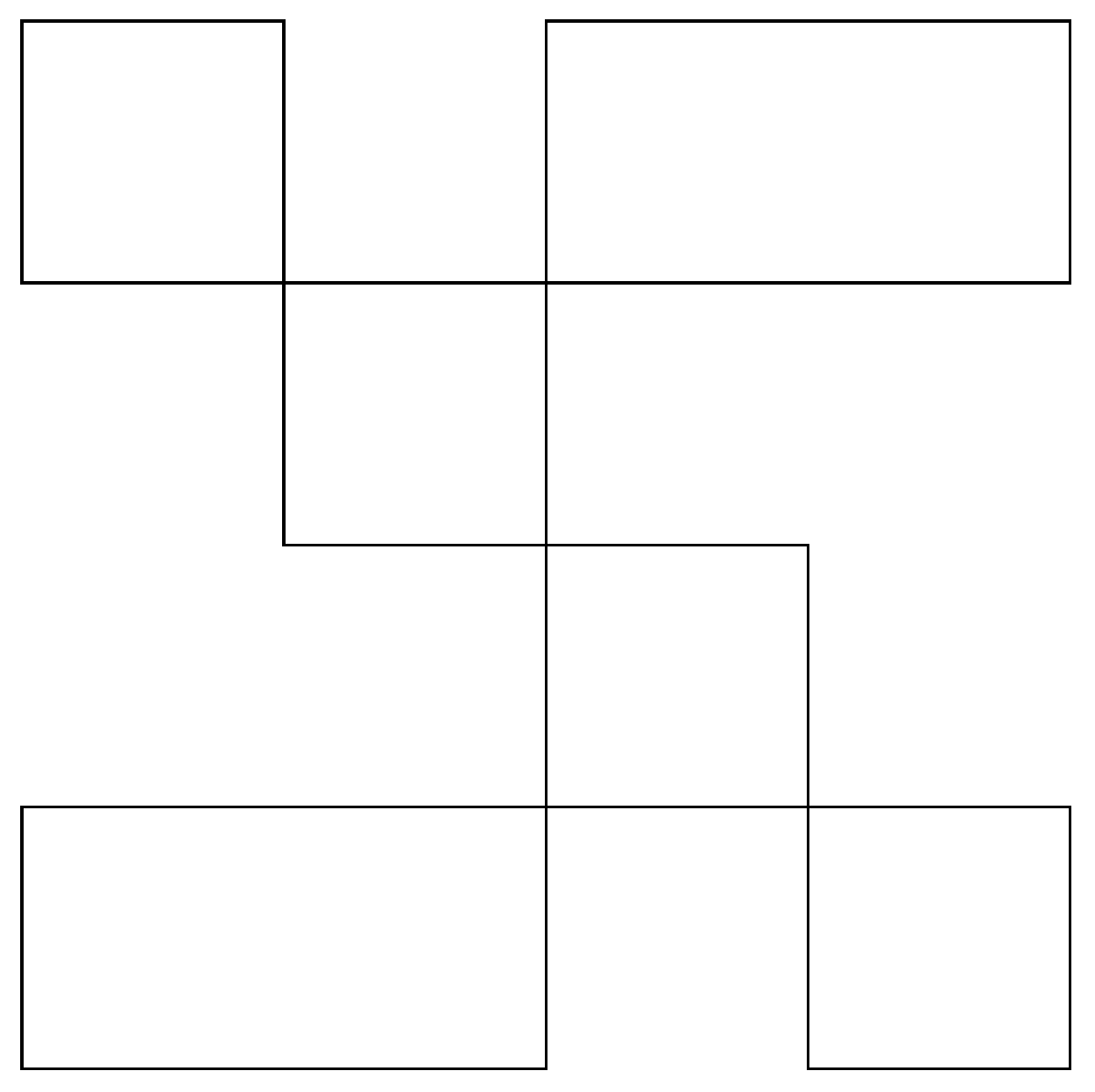}\qquad\quad
	\includegraphics[height=2.7cm]{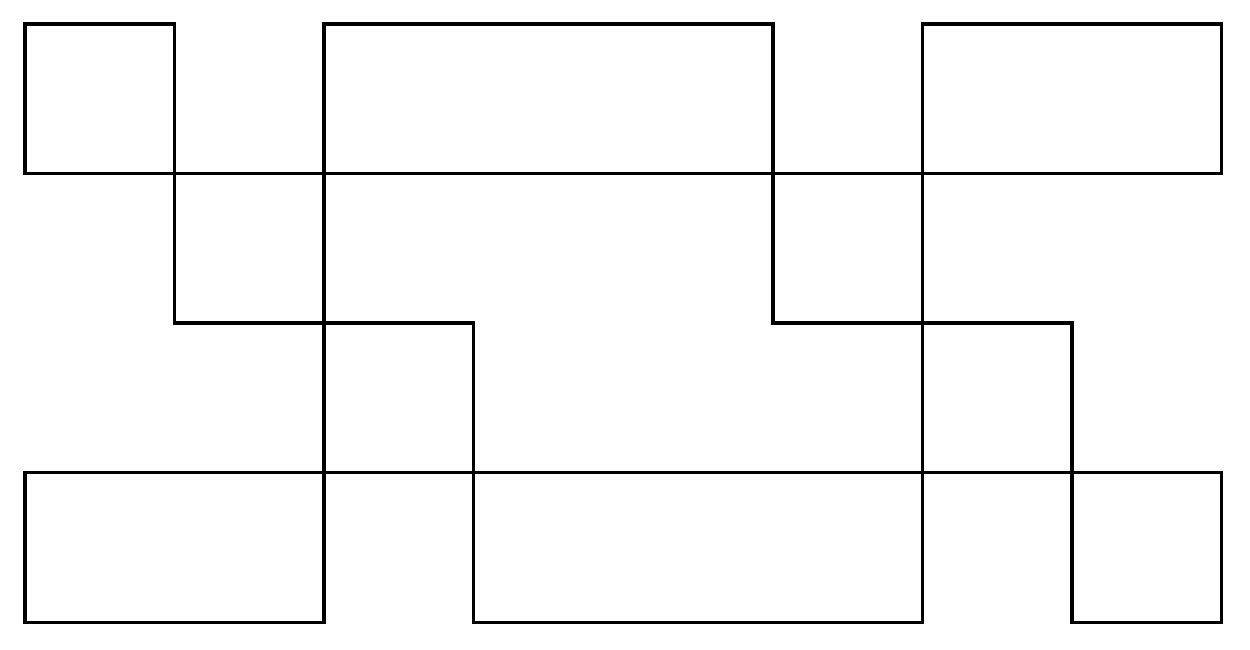}\\
	\vspace{0.2cm}
	\includegraphics[height=2.7cm]{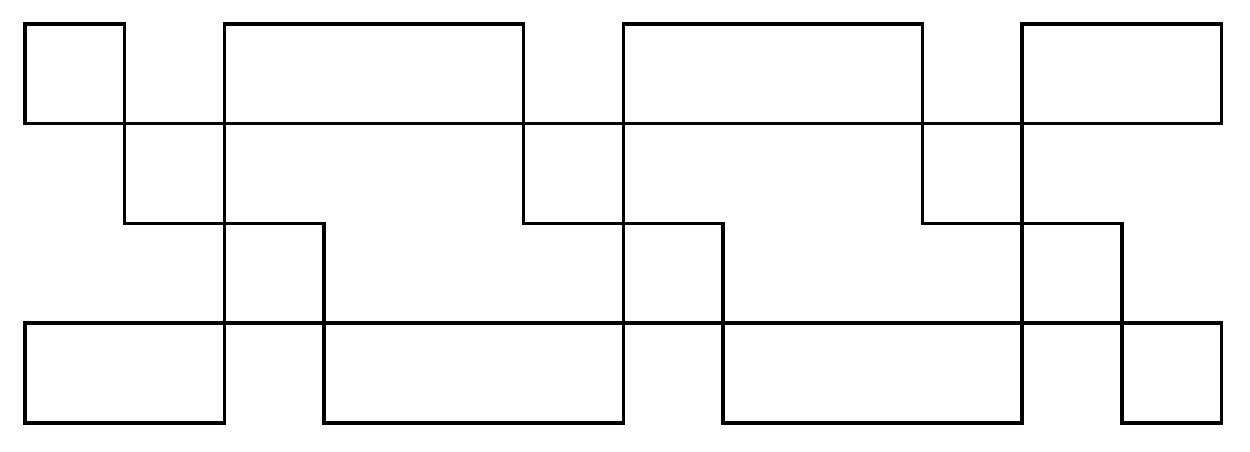}	
	\begin{picture}(0,0) \thicklines 
	\put(-231,147){$r_{11}$}
	\put(-184,147){$r_{12}$}
	\put(-123,147){$r_{21}$}
	\put(-70,147){$r_{23}$}
		\put(-8,147){$r_{22}$}
\put(-209,62){$r_{31}$}
\put(-32,62){$r_{32}$}
\put(-159,62){$r_{33}$}
	\end{picture}
	\begin{center}
		\caption{The renormalization factors.}
	\end{center}
 
\end{figure}

Noticing that the homogeneity requirement of $(\{D_{\alpha_i}\}_{i=1}^3,\bm{r})$ also implies that
\[\prod_{k=1}^{n}r_{i_{k-1}i_k}=r^n_{11},\]
for any finite head-to-tail sequences of factors $r_{i_0i_1}, r_{i_1i_2},\cdots, r_{i_{n-1}i_n}$, with $i_0=i_n=1$. An immediate observation is $r_{12}r_{21}=r_{11}^2,r_{12}r_{23}r_{31}=r^3_{11}$.
We  see that $r_{11}=r_{22}=r_{33}$ if we look at the pair of admissible sequences $121,1221$, and the pair of admissible sequences $1231,12331$. In  a similar way, we can get that  $r_{23}r_{32}=r^2_{22}=r^2_{11}$. So there are only three free parameters $r_{11},r_{21},r_{31}$. 
Furthermore, if $(\{D_{\alpha_i}\}_{i=1}^3,\bm{r})$ is a homogeneous regular harmonic structure with $\bm{r}=\{r_{ij}\}$ as above, then by letting $D'_{\alpha_1}=D_{\alpha_1}$, $D'_{\alpha_2}=\frac{r_{11}}{r_{21}}D_{\alpha_2}$, $D'_{\alpha_3}=\frac{r_{11}}{r_{31}}D_{\alpha_3}$, and $\bm{r'}$ the vector of constant $r_{11}$, it is easy to check that $(\{D'_{\alpha_i}\}_{i=1}^3,\bm{r}')$ is also a homogeneous regular harmonic structure, which yields the same resistance form induced by $(\{D_{\alpha_i}\}_{i=1}^3,\bm{r})$. 

In view of the above discussion, we only need to consider the  homogeneous regular harmonic structures $(\{D_{\alpha_i}\}_{i=1}^3,\bm{r})$ with $\bm{r}$ being a constant vector. To simplify the notations, we denote the common renormalization factor by $r$.

On the other hand, unlike the previous examples, it is easy to observe that for each island $K_\beta$ in the $f.r.f.t.$ nested structure, its boundary $V_\beta$ is not fully involved when $K_\beta$ intersects other islands, i.e., $V_\beta\setminus \bigcup_{\alpha\in\Lambda, K_\alpha\nsubseteq K_\beta}K_\alpha\neq\emptyset$, $\forall\beta\in\Lambda$. In fact, for a $K_{\alpha_1}$ type island, there are two manners up to the rotational symmetry, with $2$ or $3$ boundary points involved, and for a $K_{\alpha_2}$ or $K_{\alpha_3}$ type island, there is only one manner for it to intersect others up to the rotational symmetry, with $2$ boundary points involved. 

Thus, regarding the rotational symmetry, we only need to consider certain restricted resistance networks of $(V_{\alpha_i}, D_{\alpha_i})$'s. Firstly, we choose a $Y$-shaped restricted network on $\{q_2,q_3,q_4\}$ of $(V_{\alpha_1}, D_{\alpha_1})$, and denote the resistances to be $a,b,c$. There is also a symmetric version of restricted network on $\{q_1,q_2,q_4\}$. By simple series connection, the effective resistance between $q_2$ and $q_4$ is always  $a+c$. Secondly, for $(V_{\alpha_2},D_{\alpha_2})$, we restrict the network onto $\{q_1,F_5q_3\}$, two of the diagonal vertices in $V_{\alpha_2}$, and set the resistance between them to be $d$. Lastly, for $(V_{\alpha_3},D_{\alpha_3})$, we restrict it onto $\{F_1F_8q_1,F_5F_4q_2\}$ or $\{F_1F_8q_4,F_5F_4q_3\}$, two of the four vertices in $V_{\alpha_3}$, lying on a long side, and denote the resistance by $e$, which are same by the rotational symmetry. See Figure 5.11 for the three types of restricted networks.

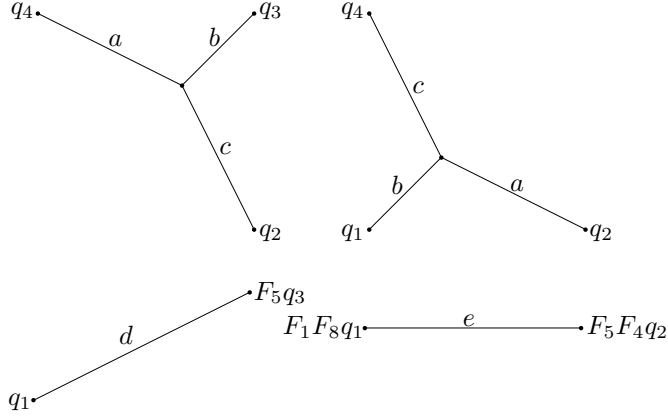
\begin{figure}[h]
	\centering
	\includegraphics[width=3cm]{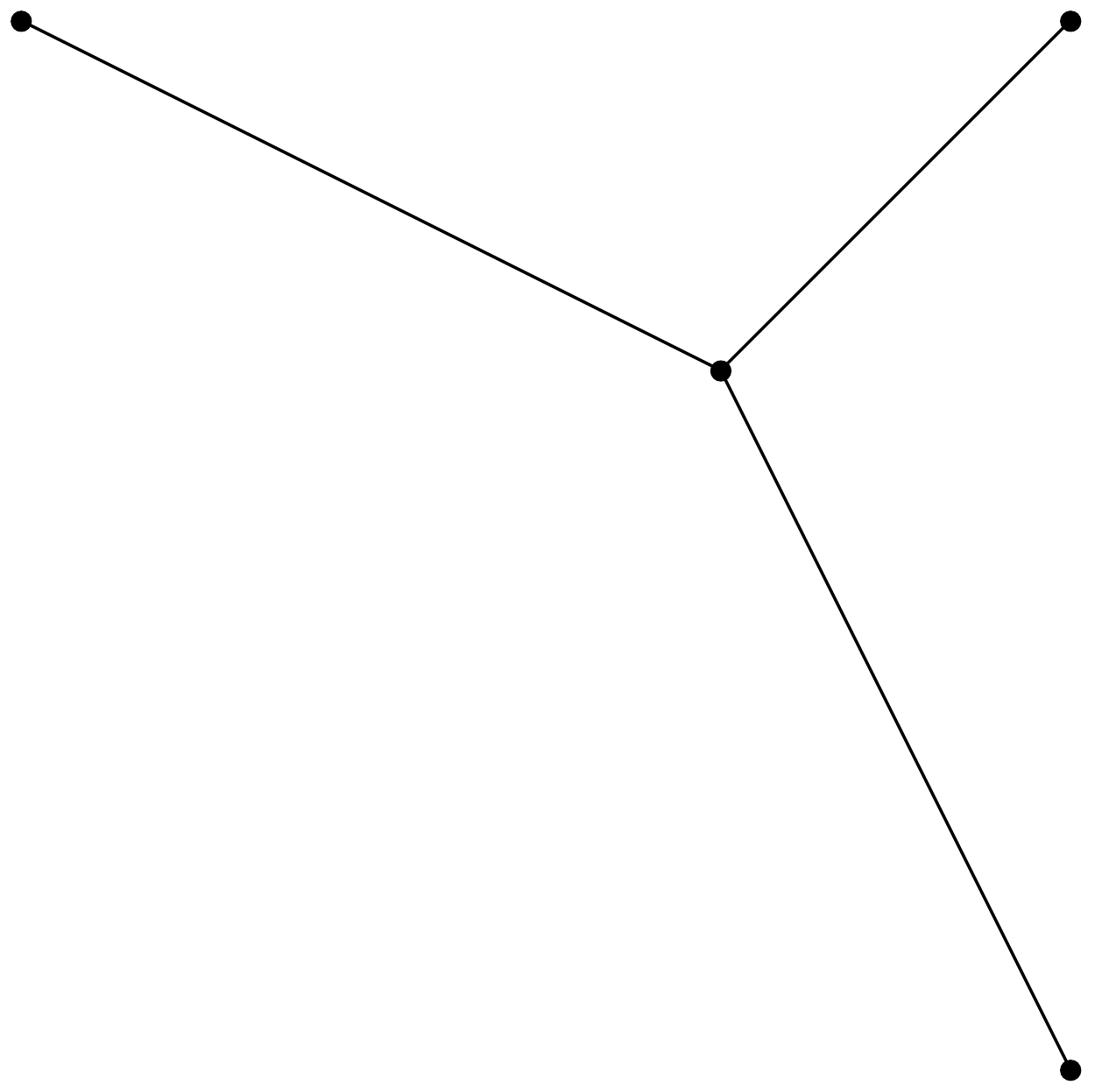}\qquad\qquad
	\includegraphics[width=3cm]{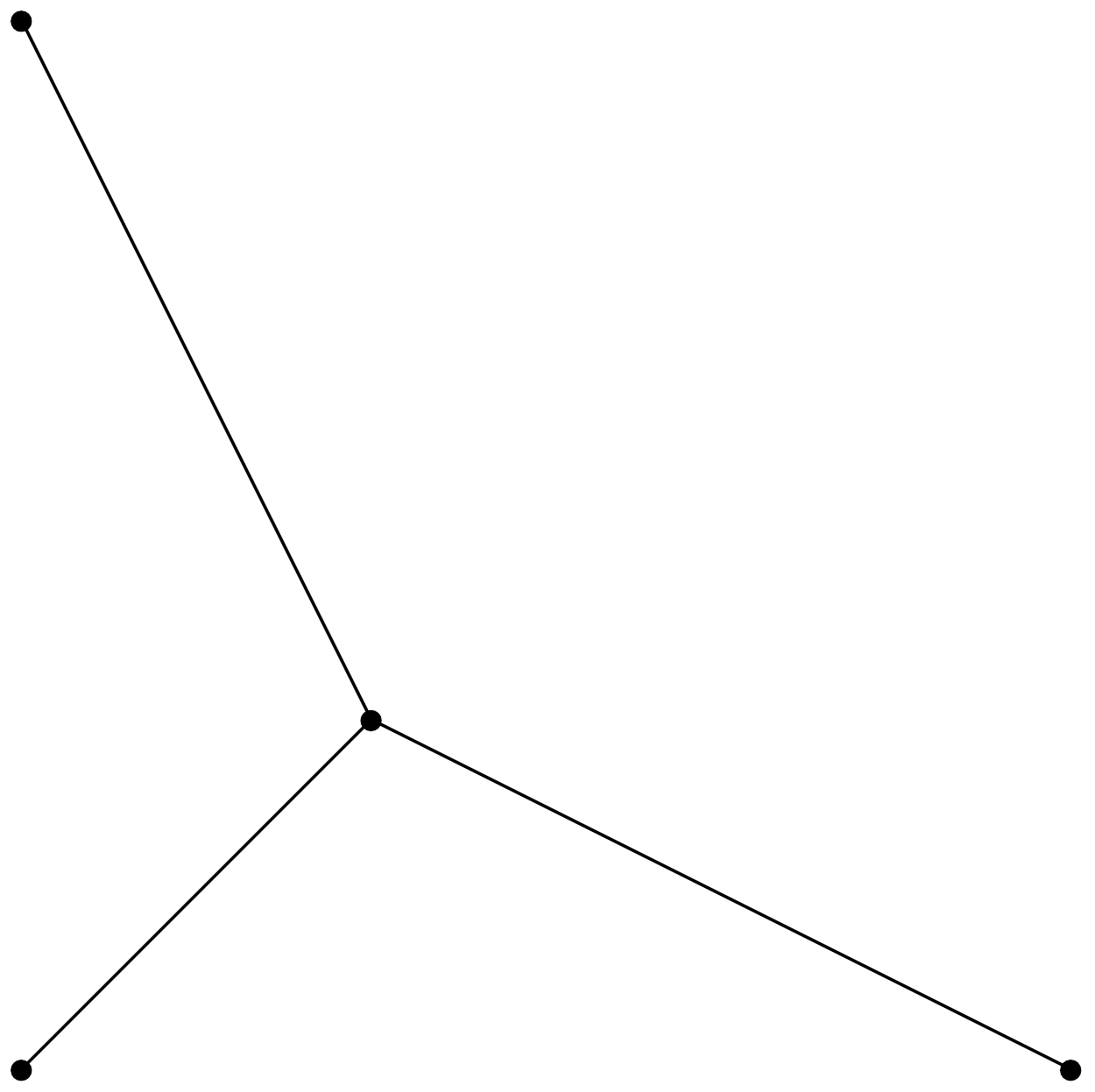}\\\vspace{0.2cm}
	\includegraphics[width=3cm]{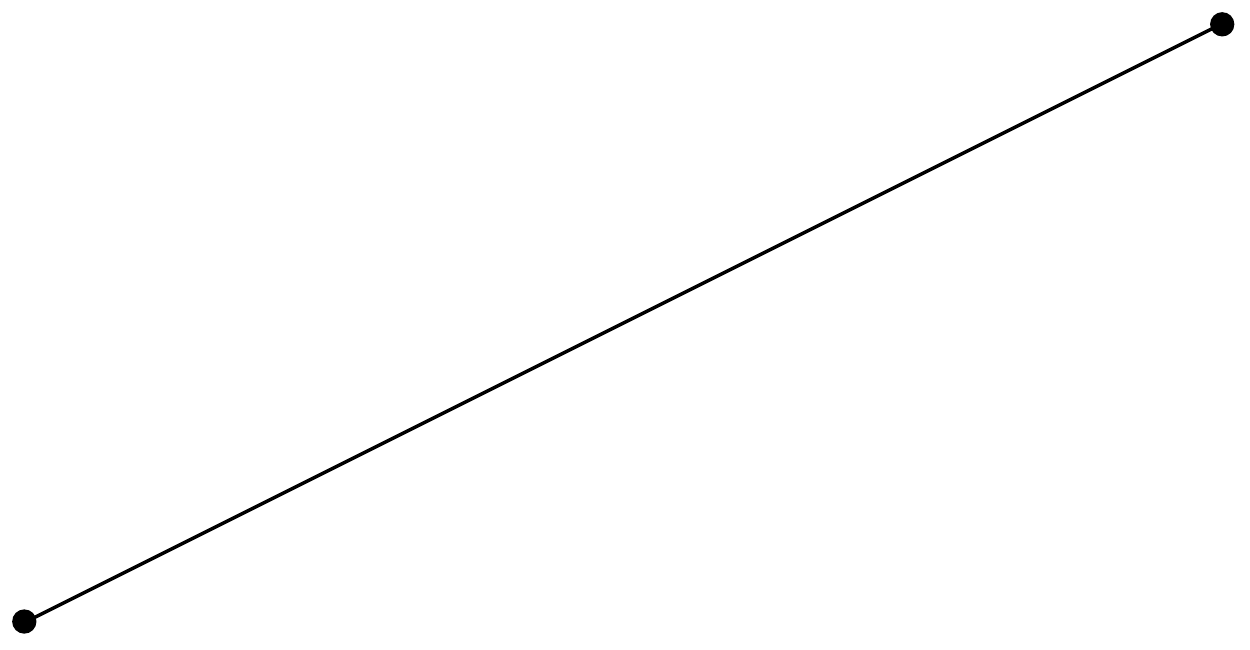}\qquad\qquad
	\includegraphics[width=3cm]{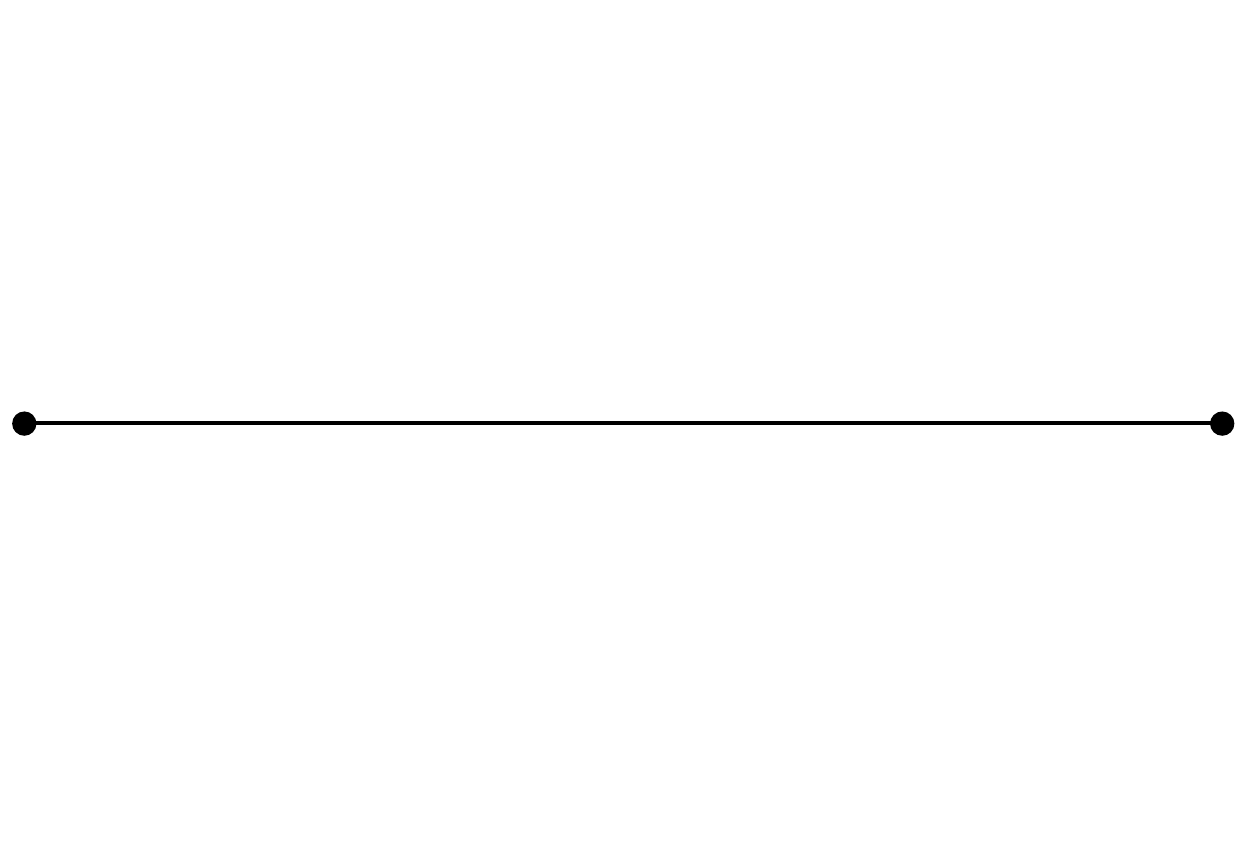}
	\begin{picture}(0,0) \thicklines 
	\put(-221,148){$q_4$}
	\put(-127,148){$q_3$}
	\put(-127,64){$q_2$}
	\put(-184,136){$a$}
	\put(-146,136){$b$}
	\put(-142,95){$c$}
		
	\put(-96,148){$q_4$}
	\put(-2,64){$q_2$}
	\put(-96,64){$q_1$}
	\put(-32,81){$a$}
	\put(-77,80){$b$}
	\put(-69,118){$c$}
	
	\put(-222,0){$q_1$}
	\put(-180,23){$d$}
	\put(-129,40){$F_5q_3$}
	
	\put(-118,27){$F_1F_8q_1$}	
	\put(-50,30){$e$}
	\put(-3,27){$F_5F_4q_2$}		
	\end{picture}
	\begin{center}
		\caption{The three types of restricted resistance networks in $(V_{\alpha_i}, D_{\alpha_i})$'s.}
	\end{center}
\end{figure}

Now we come to the calculations.

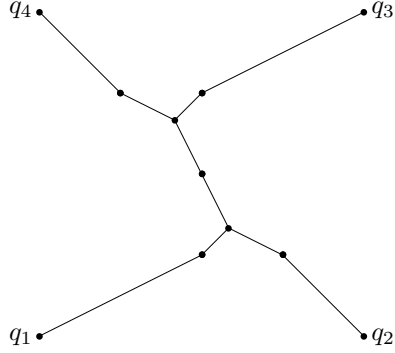
\begin{figure}[h]
	\centering
	\includegraphics[width=4.5cm]{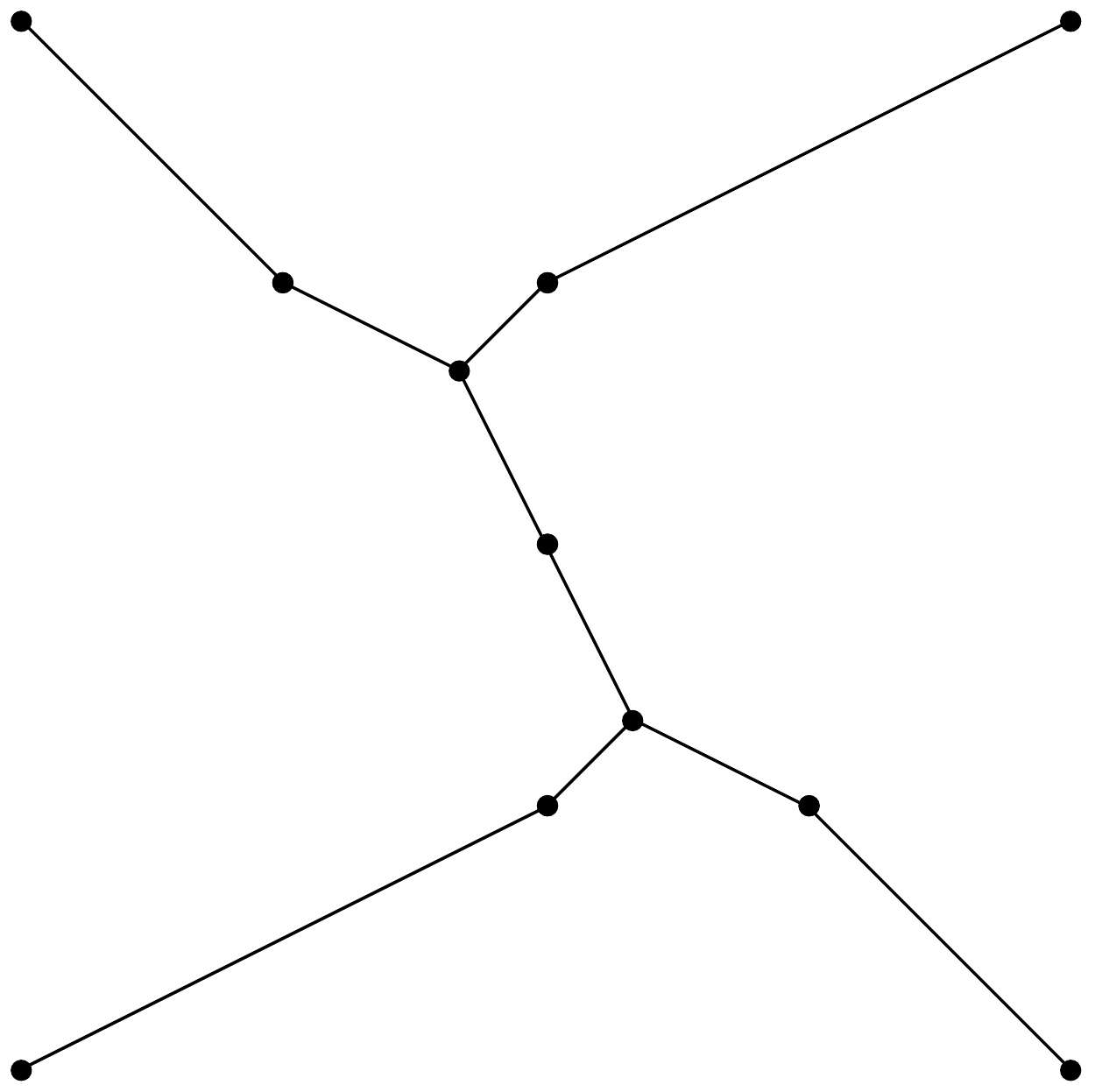}\qquad\qquad
	\begin{picture}(0,0) \thicklines 
	\put(-177,125){$q_4$}
	\put(-40,125){$q_3$}
	\put(-40,1){$q_2$}	
	\put(-177,1){$q_1$}		
	\end{picture}
	\begin{center}
		\caption{The Level-$1$ resistance network associated with $(V_{\alpha_1}, D_{\alpha_1})$.}
	\end{center}
\end{figure}

First, look at the level-$1$ resistance network associated with $(V_{\alpha_1}, D_{\alpha_1})$, generated by the above level-$0$ restricted networks, see Figure 5.12. By comparing the 
effective resistances between $q_i,q_j$ with that of $(V_{\alpha_1}, D_{\alpha_1})$ for distinct $(i,j)$'s, using series connection, we have

\[\begin{cases}
r(2a+b+c+d)=a+b,\\
r(4a+4c)=a+c,\\
r(2a+b+3c+d)=b+c.
\end{cases}\]
Solving the equations, we get 
\begin{equation}
r=\frac{1}{4},\quad c=2a,\quad d=3b.
\end{equation}

\begin{figure}[h]
	\centering
	\includegraphics[width=8cm]{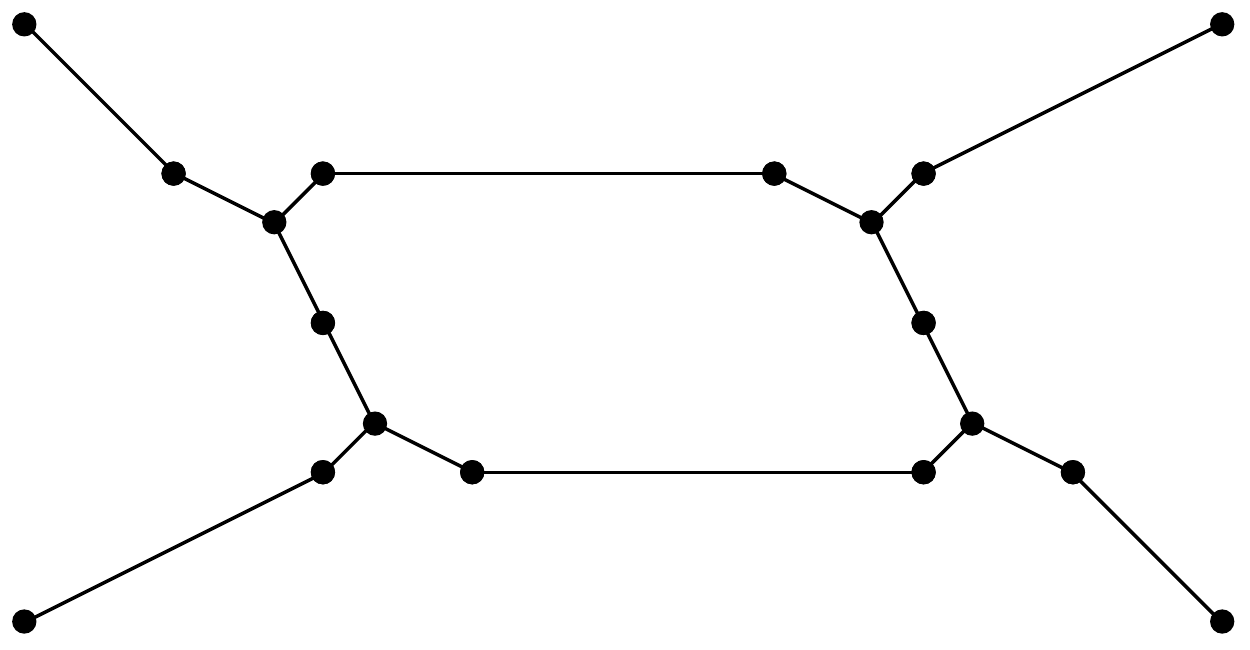}\qquad\qquad
		\begin{picture}(0,0) \thicklines 
	\put(-286,112){$F_1q_4$}
	\put(-42,112){$F_5q_3$}
	\put(-42,2){$F_5q_2$}	
	\put(-275,3){$q_1$}		
	\end{picture}
	
	\begin{center}
		\caption{The Level-$1$ resistance network associated with $(V_{\alpha_2}, D_{\alpha_2})$.}
	\end{center}
\end{figure}
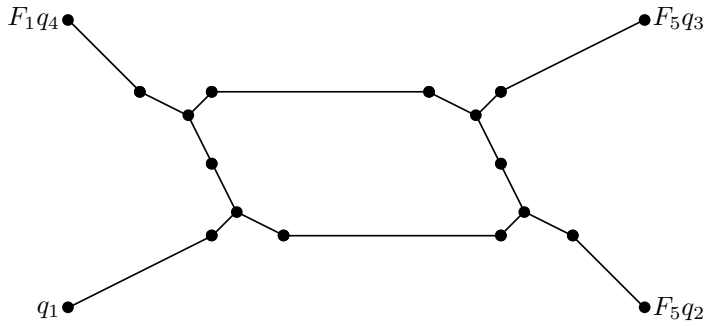

Next, we look at the level-$1$ resistance network associated with $(V_{\alpha_2}, D_{\alpha_2})$, shown in Figure 5.13. We just need to compare the effective resistance between $q_1$ and $F_5q_3$ with that of $(V_{\alpha_2}, D_{\alpha_2})$. Using series and parallel connection of resistors, we get 
\[d=r\big(2b+2d+\frac{1}{2}(a+b+2c+e)\big).\]
Substituting (5.8) into the above equation, we have
\begin{equation}\label{eq520}
e=7b-5a.
\end{equation}

\begin{figure}[h]
	\centering
	\includegraphics[width=8cm]{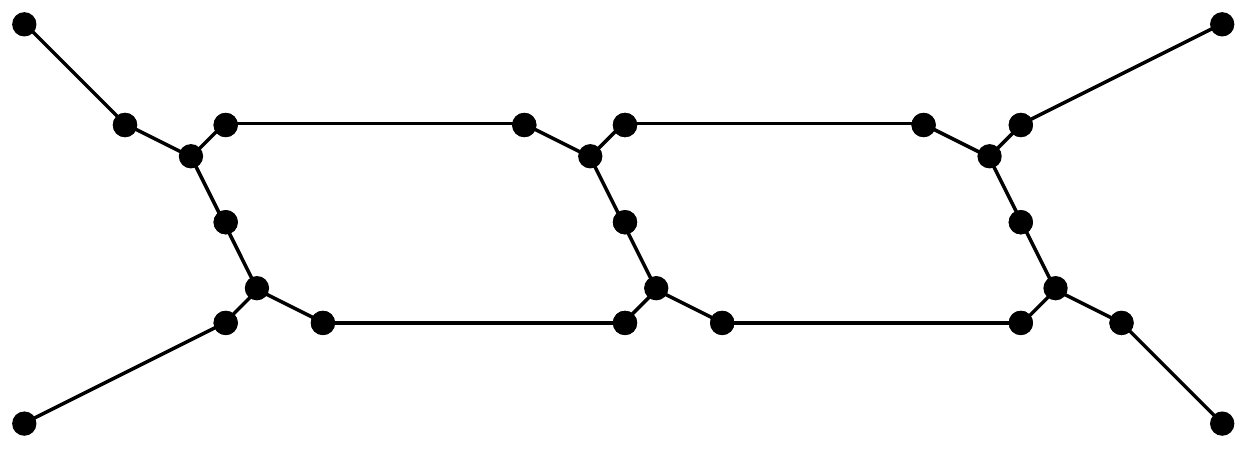}\qquad\qquad
	\begin{picture}(0,0) \thicklines 
	\put(-296,78){$F_1F_8q_4$}
	\put(-41,78){$F_5F_4q_3$}
	\put(-41,1){$F_5F_4q_2$}	
	\put(-296,1){$F_1F_8q_1$}		
	\end{picture}
	\begin{center}
		\caption{The Level-$1$ resistance network associated with $(V_{\alpha_3}, D_{\alpha_3})$.}
	\end{center}
\end{figure}
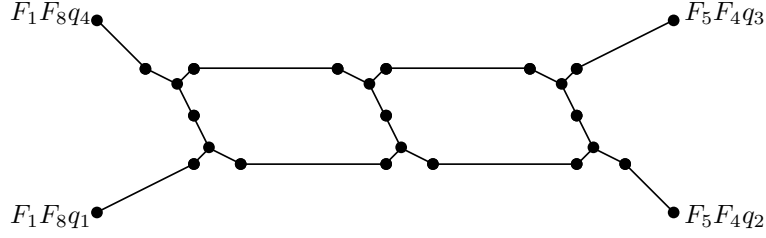

Finally, we look at the level-$1$ resistance network associated with $(V_{\alpha_3},D_{\alpha_3})$, shown in Figure 5.14. By using series connection operation, we    simplify the network into what Figure 5.15 presents.

\begin{figure}[h]
	\centering
	\includegraphics[width=8cm]{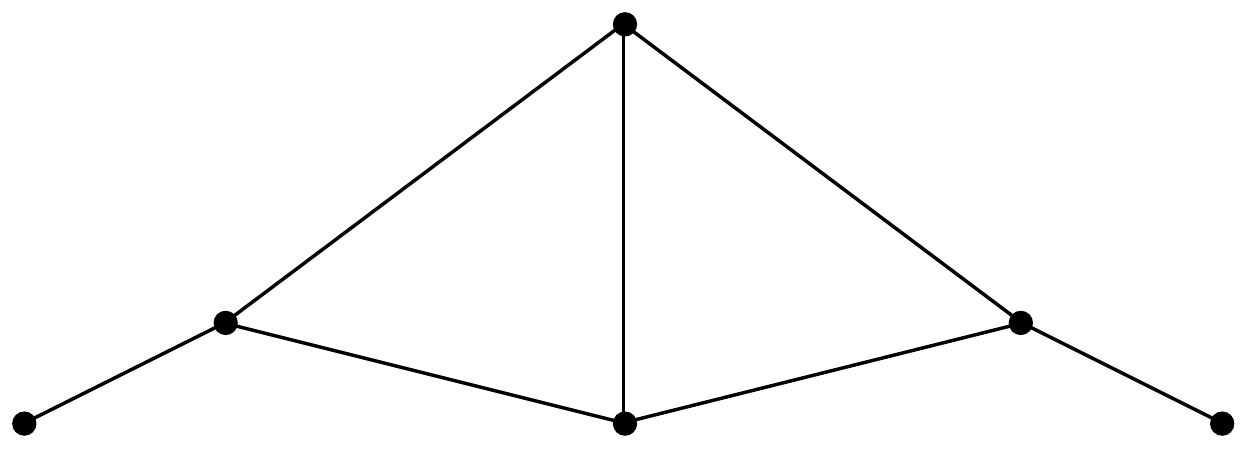}\qquad\qquad
	\begin{picture}(0,0) \thicklines 
	\put(-39,0){$F_5F_4q_2$}
	\put(-296,0){$F_1F_8q_1$}
	\put(-276,16){$r(b+d)$}
	\put(-210,20){$r(a+b+e)$}
	\put(-148,20){$r(a+b+e)$}
	\put(-255,57){$r(a+b+2c+e)$}
	\put(-125,57){$r(a+b+2c+e)$}
	\put(-152,40){$2rc$}
	\put(-63,15){$r(2a+c)$}
	\end{picture}
	\begin{center}
		\caption{A simplification of the resistance network in Figure 5.14.}
	\end{center}
\end{figure}
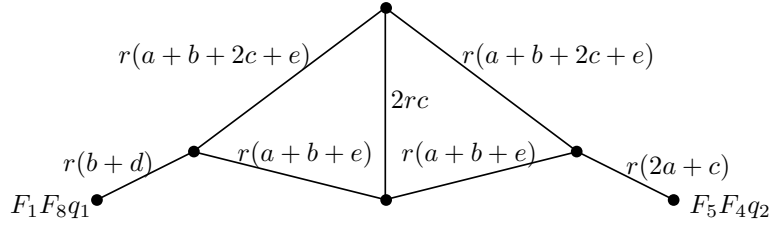

Then by using the symmetry of the above network, we can easily calculate the effective resistance between $F_1F_8q_1$ and $F_5F_4q_2$, which gives that 
\[r(2a+b+c+d+\frac{2}{(a+b+2c+e)^{-1}+(a+b+e)^{-1}})=e.\]
Substituting (5.8) and (5.9) into the above equation, we get
\begin{equation}
b=\frac{1}{16}(13+\sqrt{73})a,\quad c=2a, \quad d=\frac{3}{16}(13+\sqrt{73})a,\quad e=\frac{1}{16}(11+7\sqrt{73})a,\quad r=\frac{1}{4}.
\end{equation}

Thus we have

\textbf{Theorem 5.8.} \textit{For the $f.r.f.t.$ nested structure  of $\mathcal{V}^w$ shown in Example 3, there exists exactly one rotational symmetric homogeneous regular harmonic structures up to scalar constants.}

\subsection{Symmetrical overlapping gasket with closed bottom} We only consider the symmetric homogeneous regular harmonic structures for this example. Since the involved calculation is complicated comparing with the previous examples, instead of providing the exact result, we just give a  quick analysis about the dimension of the solutions, i.e, the number of free parameters to decide a symmetric homogeneous regular harmonic structure.

As before, we need to list the possible parameters to be determined, including the resistances in the networks $(V_{\alpha_i}, D_{\alpha_i}), i=1,2,$  and the renormalization factors. See Figure 5.16 for the  resistances, and Figure 5.17 for the renormalization factors.

\begin{figure}[h]
	\includegraphics[height=2.7cm]{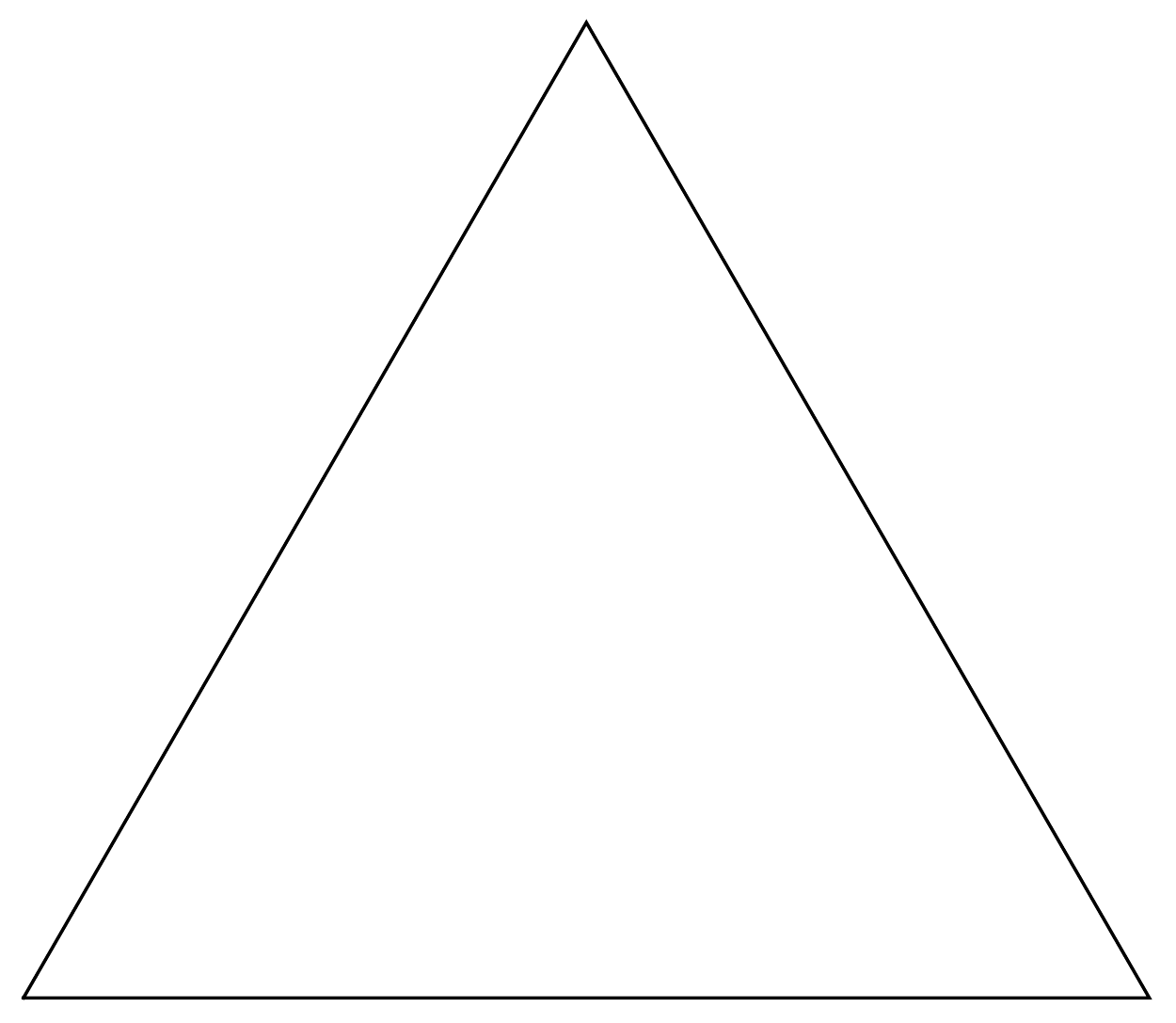}\qquad\qquad
	\includegraphics[height=2.7cm]{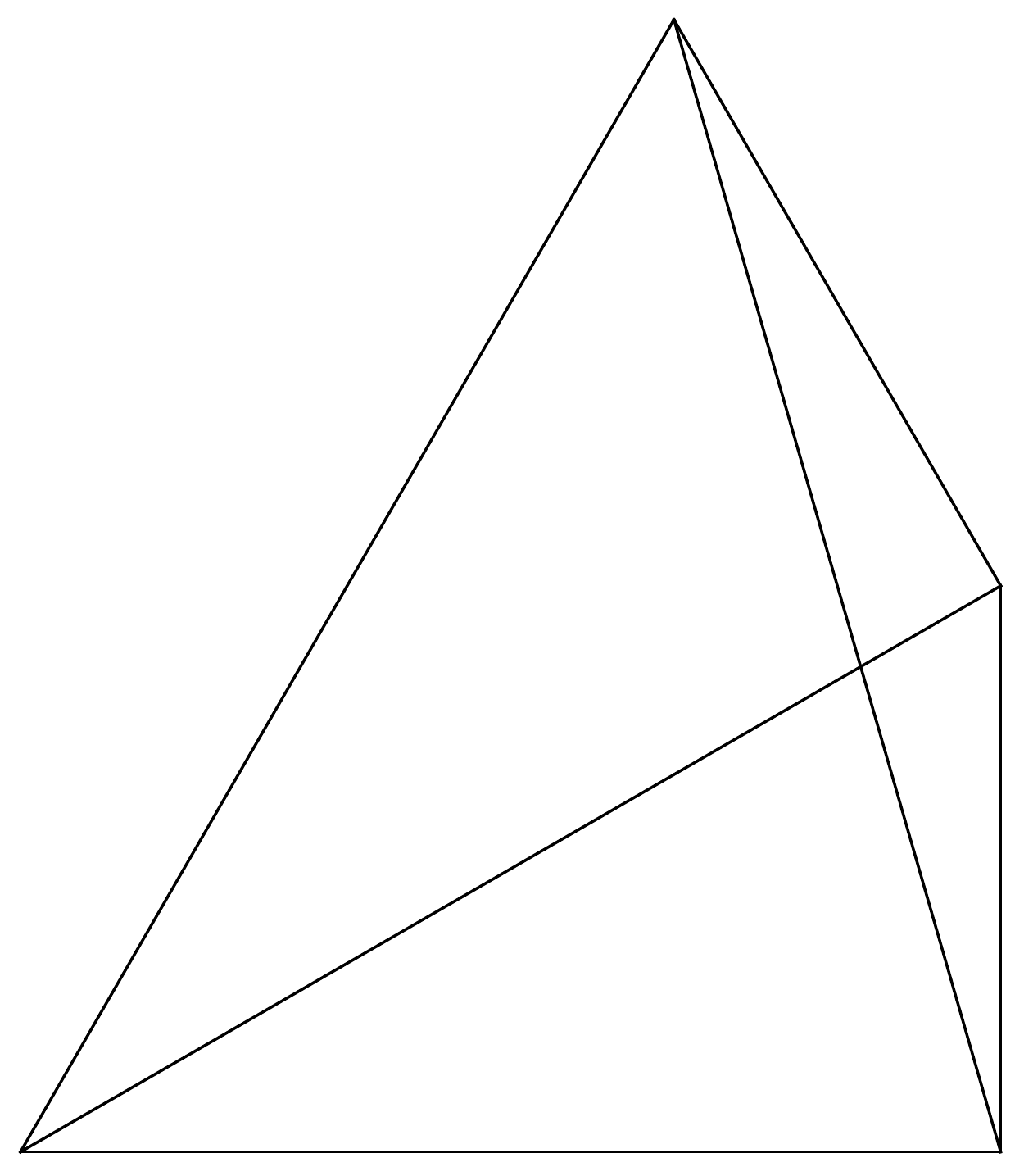}
	\begin{picture}(0,0)
	\put(-158,-7){$1$}
	\put(-132,40){$a$}
	\put(-182,40){$a$}
	
	\put(-38,-7){$b$}
    \put(-53,40){$c$}
	\put(-14,56){$d$}
	\put(-4,20){$e$}
	\put(-42,21){$f$}
	\put(-22,34){$g$}
	\end{picture}
	\caption{The resistance networks $(V_{\alpha_i}, D_{\alpha_i}), i=1,2$.}
\end{figure}

\begin{figure}[h]
	\includegraphics[height=3.4cm]{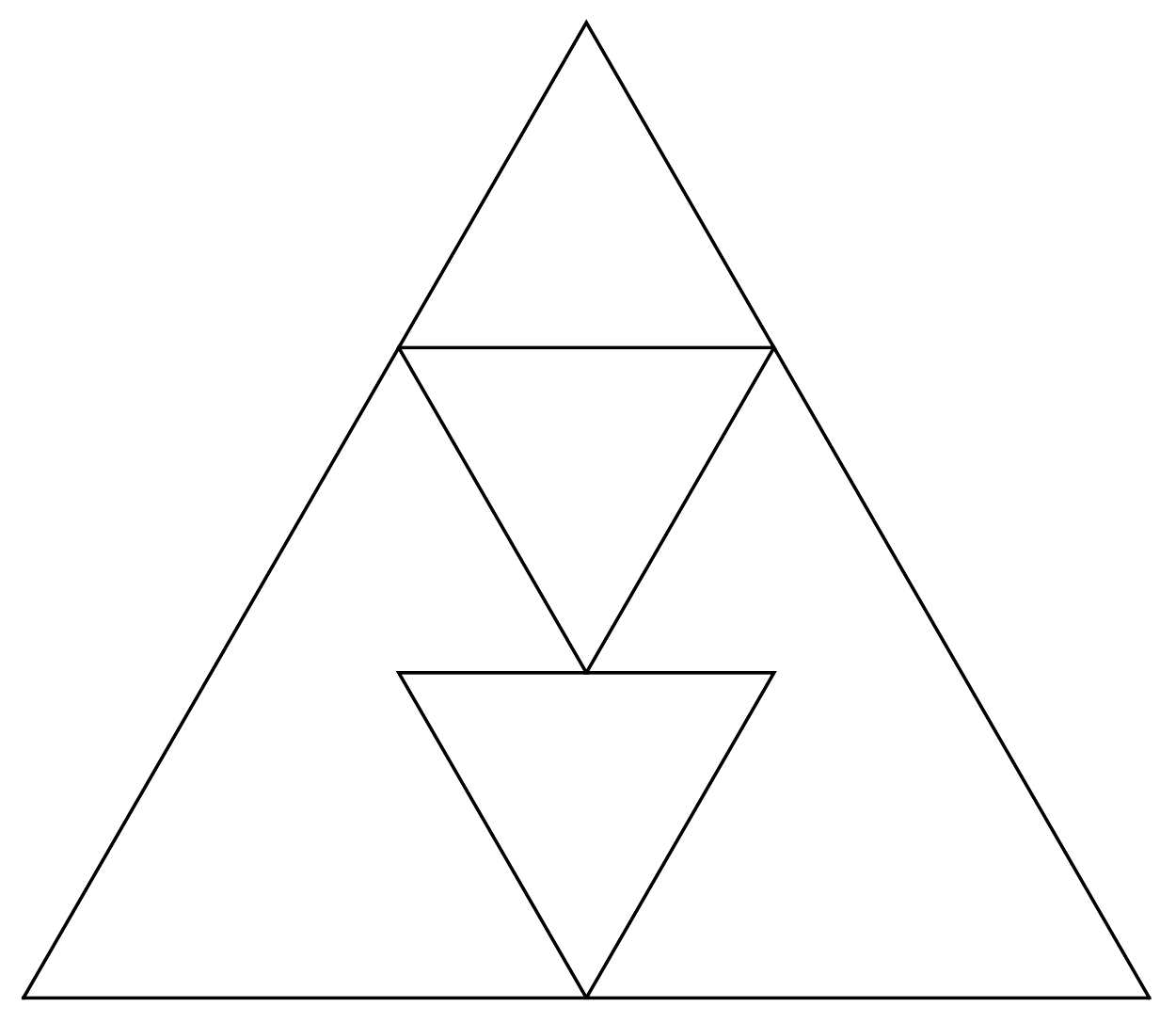}\qquad\qquad
	\includegraphics[height=3.4cm]{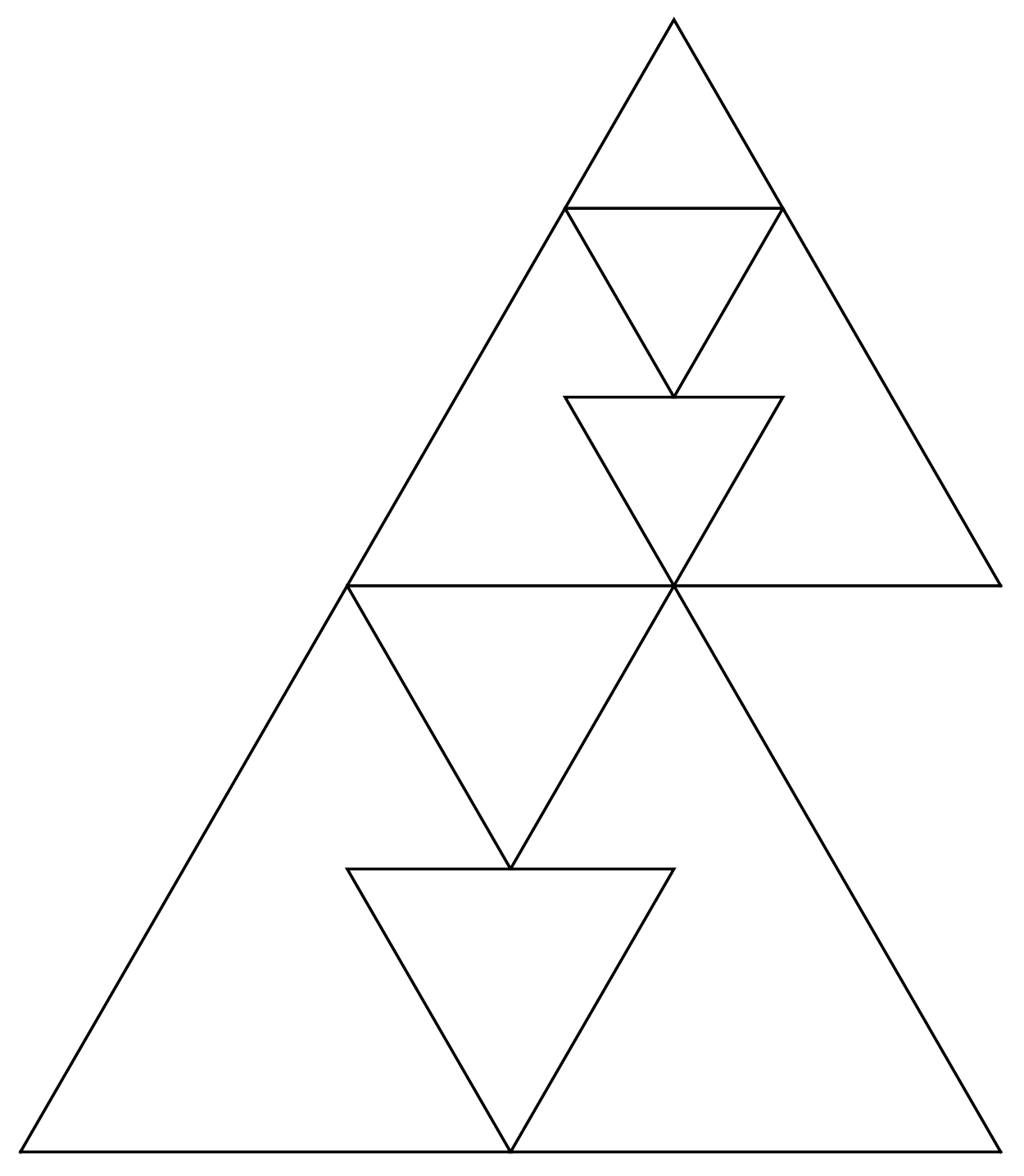}
	\begin{picture}(0,0)
	\put(-212,27){$1$}
	\put(-160,27){$1$}
    \put(-188,70){$r_1$}
    
    \put(-36,84){$r_2$}
    \put(-50,54){$r_3$}
    \put(-22,54){$r_3$}
    \put(-68,24){$r_4$}
    \put(-31,24){$r_4$}
	\end{picture}
	\caption{The renormalization factors.}
\end{figure}
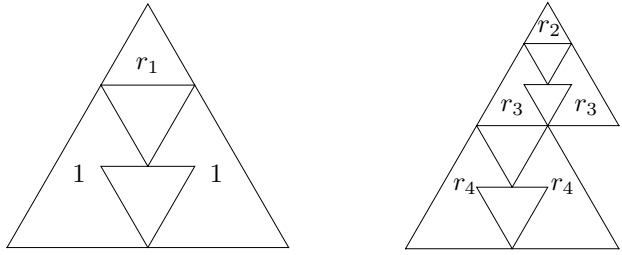

We have $7$ parameters for the resistances and $4$ parameters for the renormalization factors. The renormalization equations will give us $8$ equations among these parameters. Nevertheless, the homogeneity assumption will give us two more equations. In fact, by noticing that $F_4\mathcal{SG}^c$ is isometric to $F_1\mathcal{SG}^c$, and then comparing the pieces of these two copies, we have
\[\begin{cases}
r_1^2=r_2\cdot 1,\\
r_3\cdot 1=1\cdot r_1.
\end{cases}\]
So we have $10$ equations altogether. It is expectable that there exist the solutions of symmetric homogeneous regular harmonic structures with $1$ free parameter.

\section{Homogeneous regular harmonic structures}

In this section, we will focus on the homogeneous regular harmonic structures of an $f.r.f.t.$ self-similar set $K$. 

\textbf{Definition 6.1.} \textit{Let $(\{D_{\alpha_i}\}_{i=1}^M, \bm{r})$ be a harmonic structure of an $f.r.f.t.$ nested structure $\{K_\alpha,\Lambda_\alpha\}_{\alpha\in \Lambda}$.  We say $(\{D_{\alpha_i}\}_{i=1}^M, \bm{r})$ is} homogeneous \textit{ if for any two islands $K_\alpha$, $K_\beta$ with the same type and same size, $r_{\bm{e}(\vartheta,\alpha)}=r_{\bm{e}(\vartheta,\beta)}$.}

The resistance form $(\mathcal{E},\mathcal{F})$ generated by a homogeneous regular harmonic structure is \textit{translation invariant}, i.e., for any two islands $K_\alpha,K_\beta$ with the same type and same size, 
\[\mathcal{E}(u)=\mathcal{E}(u\circ\phi_{\alpha,\beta})\]
holds for any function $u\in\mathcal{F}$ supported in $K_\beta$.

In the case that the associated directed graph $G=(S,E)$ is \textit{strongly connected},  i.e., for any two states $\mathcal{T}_i,\mathcal{T}_j$ in $S$, there is a walk $\bm{e}$ such that $i(\bm{e})=\mathcal{T}_i, f(\bm{e})=\mathcal{T}_j$, we have the following proposition. The examples in Section 4 are all in this case.

\textbf{Proposition 6.2.} \textit{Let $(\{D_{\alpha_i}\}_{i=1}^M, \bm{r})$ be a homogeneous harmonic structure of an $f.r.f.t.$ nested structure $\{K_\alpha,\Lambda_\alpha\}_{\alpha\in \Lambda}$. Suppose the associated directed graph  $G=(S,E)$ is strongly connected. Then for any $1\leq i\leq M$, for any $\alpha,\beta\in \mathcal{T}_i$, the ratio $r_{\bm{e}(\vartheta,\alpha)}/r_{\bm{e}(\vartheta,\beta)}$ depends only on the ratio $diam(K_\alpha)/ diam(K_\beta)$, i.e., there exists a function $c(\cdot):\mathbb{R^+}\to \mathbb{R^+}$ such that $r_{\bm{e}(\vartheta,\alpha)}/r_{\bm{e}(\vartheta,\beta)}=c\big(diam(K_\alpha)/ diam(K_\beta)\big)$.} 

\textit{Proof.} Let $\tilde{\alpha},\tilde{\beta}\in \mathcal{T}_j$ be another pair of indices such that $\frac{diam(K_\alpha)}{diam(K_\beta)}=\frac{diam(K_{\tilde{\alpha}})}{diam(K_{\tilde{\beta}})}$. Since $G$ is strongly connected, we can always find a walk $\bm{e}$ such that $i(\bm{e})=\mathcal{T}_i$ and $f(\bm{e})=\mathcal{T}_1$. Connecting the walks $\bm{e}(\vartheta,\alpha)$ and $\bm{e}(\vartheta,\tilde{\beta})$ by $\bm{e}$, we  get a walk $\bm{e}_1=\bm{e}(\vartheta,\alpha)\bm{e}\bm{e}(\vartheta,\tilde{\beta})$, and similarly $\bm{e}_2=\bm{e}(\vartheta,\beta)\bm{e}\bm{e}(\vartheta,\tilde{\alpha})$. Note that both $\bm{e}_1$ and $\bm{e}_2$ are walks from $\mathcal{T}_i$ to $\mathcal{T}_j$. Obviously, $\psi_{\bm{e}_1}$ and $\psi_{\bm{e}_2}$ have the same similarity ratio, which gives that
\[r_{\bm{e}_1}=r_{\bm{e}_2},\]
since $(\{D_{\alpha_i}\}_{i=1}^M, \bm{r})$ is homogeneous. As a result, $r_{\bm{e}(\vartheta,\alpha)}r_{\bm{e}(\vartheta,\tilde{\beta})}=r_{\bm{e}(\vartheta,\beta)}r_{\bm{e}(\vartheta,\tilde{\alpha})}$. Thus the ratio $r_{\bm{e}(\vartheta,\alpha)}/r_{\bm{e}(\vartheta,\beta)}$ only depends on $diam(K_\alpha)/ diam(K_\beta)$.\hfill$\square$\\

Now let's look at the relation of the homogenous regular harmonic structures between different $f.r.f.t.$ nested structures. Let $K$ be a self-similar set, with two distinct $f.r.f.t.$ nested structures $\mathcal{S}:=\{K_\alpha,\Lambda_\alpha\}_{\alpha\in \Lambda}$  and $\mathcal{S}':=\{K_{\alpha'},\Lambda'_{\alpha'}\}_{\alpha'\in\Lambda'}$. We use  $G=(S,E)$ and $G'=(S',E')$ to denote their associated directed graphs respectively.

For an island $K_{\alpha'},\alpha'\in \Lambda'$, we can always find an at most countable set of indices $L_{\alpha'}\subset \Lambda$ such that
\begin{equation}
K_{\alpha'}\setminus V_{\alpha'}=\bigcup_{\alpha\in L_{\alpha'}}K_\alpha,
\end{equation}
with 
\begin{equation}
\#K_{\alpha}\cap K_{\beta}<\infty, \forall \alpha,\beta\in L_{\alpha'}.
\end{equation}
In fact, let $\pi$ be the canonical projection $\pi:E^\infty\to K$, then $\pi^{-1}(K_{\alpha'}\setminus V_{\alpha'})$ is an open set in $E^\infty$, so it is a countable disjoint union of cylinders $\pi^{-1}(K_{\alpha'}\setminus V_{\alpha'})=\bigcup_{\alpha\in L_{\alpha'}}E^\infty_\alpha$, where each cylinder is of the form $E^\infty_\alpha=\{\bm{\epsilon}\in E^\infty: [\bm{\epsilon}]_n=\bm{e}(\vartheta,\alpha)\}$ with $\alpha\in\Lambda$ and $|\alpha|=n$. This gives the desired set $L_{\alpha'}$. Obviously, for a given $\alpha'\in\Lambda'$, $L_{\alpha'}$ is not unique, but they are subdivisions of a unique $L_{\alpha'}$. 

\textbf{Definition 6.3.} \textit{We say $\mathcal{S'}$ can be} tiled \textit{by $\mathcal{S}$, and write $\mathcal{S}'\blacktriangleleft \mathcal{S}$, if for any $\alpha', \beta'\in \Lambda'$ with $\alpha'\sim\beta'$, there exist  $L_{\alpha'}, L_{\beta'}\subset \Lambda$ satisfying equations (6.1) and (6.2) such that there is a one to one correspondence $p_{\alpha',\beta'}:L_{\alpha'}\to L_{\beta'}$  satisfying}
\[\alpha\sim p_{\alpha',\beta'(\alpha)} \text{ \textit{and} }  \phi_{\alpha,p_{\alpha',\beta'}(\alpha)}=\phi_{\alpha',\beta'}, \quad\forall \alpha\in L_{\alpha'}.\]

{\textbf{Remark.} In simple cases, for an island $K_{\alpha'}$, we can find a finite subset of indices $\tilde{L}_{\alpha'}$ such that $K_{\alpha'}=\bigcup_{\alpha\in \tilde{L}_{\alpha'}}K_{\alpha}$,
with 
$\#K_{\alpha}\cap K_{\beta}<\infty,\forall \alpha,\beta\in \tilde{L}_{\alpha'}$, so that there is a simplified analogue of Definition 6.3, which requires that for any $\alpha'\sim\beta'$, we can find finite sets of indices $\tilde{L}_{\alpha'}$ and $\tilde{L}_{\beta'}$ such that there is a same one to one correspondence $p_{\alpha',\beta'}:\tilde{L}_{\alpha'}\to \tilde{L}_{\beta'}$ as in Definition 6.3. This simple case obviously implies $S'\blacktriangleleft S$.






\textbf{Example.} (a). For the two structures $\mathcal{S}$ and $\mathcal{S'}$ on $\mathcal{V}^o$ given by Figure 4.1 and Figure 4.4, we have both $\mathcal{S}\blacktriangleleft \mathcal{S'}$ and $\mathcal{S'}\blacktriangleleft \mathcal{S}$.  

(b). There are infinitely many pairs of $f.r.f.t.$ structures $\{\mathcal{S},\mathcal{S}'\}$ on $I=[0,1]$ such that $\mathcal{S}\not\blacktriangleleft \mathcal{S'}$ and $\mathcal{S}'\not\blacktriangleleft \mathcal{S}$.

(c). Consider the standard \textit{Sierpinski gasket} $\mathcal{SG}$. Let $q_1,q_2,q_3$ be the vertices of an equilateral triangle, and let $F_i: x\rightarrow \frac{1}{2}(x+q_i), i=1,2,3$. The Sierpinski gasket $\mathcal{SG}$ is the unique invariant set satisfying $\mathcal{SG}=\bigcup_{i=1}^3 F_i\mathcal{SG}$. It is a typical $p.c.f.$ self-similar set. The  $i.f.s.$ $\{F_i\}_{i=1}^3$ gives a natural $f.r.f.t.$ nested structure $\mathcal{S}=\{F_w \mathcal{SG},\Lambda_w\}_{w\in W_*}$ as discussed in Section 3. We can also generate $\mathcal{SG}$ by using the $i.f.s.$ consisting of all nine compositions $\{F_{ij}\}_{i,j=1}^3$, and this gives another $f.r.f.t.$ nested structure $\mathcal{S}'=\{F_{w'}\mathcal{SG},\Lambda_{w'}'\}_{w'\in {W}_*'}$. 

It is  easy to check that  $\mathcal{S}'\blacktriangleleft \mathcal{S}$. But the other direction is not true. In fact, the map $\phi_{1,11}=F_1$ is never a map of the form $\phi_{\alpha',\beta'}$ for some $\alpha',\beta'\in W_*'$.\\

\textbf{Theorem 6.4.} \textit{Let $K$ be a self-similar set, with two distinct $f.r.f.t.$ nested structures $\mathcal{S}=\{K_\alpha,\Lambda_\alpha\}_{\alpha\in \Lambda}$ and $\mathcal{S}'=\{K_{\alpha'},\Lambda_{\alpha'}'\}_{\alpha'\in \Lambda'}$. Suppose $(\{D_{\alpha_i}\}_{i=1}^M, \bm{r})$ is a homogeneous regular harmonic structure of $\mathcal{S}$, and $(\mathcal{E},\mathcal{F})$ is its induced resistance form.  Assume $S'\blacktriangleleft S$ and $G$ is strongly connected, where $G=(S,E)$ is the associated directed graph of $\mathcal{S}$. Then there is a homogeneous regular harmonic structure $(\{D_{{\alpha_j'}}\}_{j=1}^{M'}, \bm{r}')$ of $\mathcal{S}'$ inducing the same resistance form $(\mathcal{E},\mathcal{F})$.}

\textit{Proof.} For each function $u\in \mathcal{F}$, we denote its associated \textit{energy measure} by $\mu_{\mathcal{E},u}$, then for each $\alpha\in \mathcal{T}_i$ with $1\leq i\leq M$,
\[\mu_{\mathcal{E},u}(K_\alpha)=r^{-1}_{\bm{e}(\vartheta,\alpha)}\mathcal{E}_{\alpha_i}(u|_{K_\alpha}\circ \phi_{\alpha_i,\alpha}).\]  
The energy measure $\mu_{\mathcal{E},u}$ has no atom by a routine discussion, see for example [BST]. So for each $\alpha'\in\Lambda'$, we have
\[\mu_{\mathcal{E},u}(K_{\alpha'})=\sum_{\alpha\in L_{\alpha'}}\mu_{\mathcal{E},u}(K_\alpha),\]
where $L_{\alpha'}\subset \Lambda$ is a countable set of indices satisfying (6.1) and (6.2). 

For each $\alpha'\in\Lambda'$, let $\mathcal{F}_{\alpha'}:=\{u|_{K_{\alpha'}}: u\in\mathcal{F}\}$, and denote $\mathcal{E}_{\alpha'}(f):=\mu_{\mathcal{E},u}(K_{\alpha'})$ for each $f\in \mathcal{F}_{\alpha'}$ with $f=u|_{K_{\alpha'}}$. It is easy to check that the value of $\mathcal{E}_{\alpha'}(f)$ is independent of the choice of $u$. By using the polarization identity
$$\mathcal{E}_{\alpha'}(f,g):=\frac{1}{4}\big(\mathcal{E}_{\alpha'}(f+g)-\mathcal{E}_{\alpha'}(f-g)\big),\quad\forall f,g\in \mathcal{F}_{\alpha'},$$
we can get a bilinear form $(\mathcal{E}_{\alpha'},\mathcal{F}_{\alpha'})$. We claim that $(\mathcal{E}_{\alpha'},\mathcal{F}_{\alpha'})$ is a resistance form on $K_{\alpha'}$. Recall the definition of resistance forms, see [K5,K7]. We need to check the following conditions (RF1) through (RF5).

\textit{RF1.}  $\mathcal{F}_{\alpha'}$ is a linear space of functions on $K_{\alpha'}$ containing constants, and $\mathcal{E}_{K_{\alpha'}}(f)=0$ if and only if $f$ is constant on $K_{\alpha'}$.

\textit{RF2.} Write $f\sim g$ if $f-g$ is a constant. Then $(\mathcal{F}_{\alpha'}/\sim,\mathcal{E}_{\alpha'})$ is a Hilbert space.

In fact, by the standard theory of resistance forms (see [K7]), for any $f\in \mathcal{F}_{\alpha'}$, there exists a unique function $h(f)\in\mathcal{F}$ such that $h(f)|_{K_{\alpha'}}=f$ and 
$$\mathcal{E}\big(h(f)\big)=\min\{\mathcal{E}(v): v\in\mathcal{F}, v|_{K_{\alpha'}}=f\},$$ which is the harmonic extension of $f$. It is easy to see that
\begin{equation}
C\mathcal{E}\big(h(f)\big)\leq \mathcal{E}_{\alpha'}(f)\leq \mathcal{E}\big(h(f)\big),
\end{equation}
for some constant $0<C<1$ independent of $f$, since the boundary of $K_{\alpha'}$ is a finite set. If $\{f_n\}_{n\geq 1}$ is a Cauchy sequence in $(\mathcal{F}_{\alpha'}/\sim,\mathcal{E}_{\alpha'})$, then $\{h(f_n)\}_{n\geq 1}$, the sequence of  harmonic extension functions, is also a Cauchy sequence in $(\mathcal{F}/\sim, \mathcal{E})$. By the completeness of $(\mathcal{F}/\sim, \mathcal{E})$, $h(f_n)$ converges to some $u^*\in\mathcal{F}$ as $n$ goes to infinity, and this gives that $\lim_{n\to \infty}\mathcal{E}_{\alpha'}(f_n-f)=0$ for $f=u^*|_{K_{\alpha'}}$.

\textit{RF3.} For any finite subset $V\subset K_{\alpha'}$ and any function $f\in l(V)$, there exists $\tilde{f}\in \mathcal{F}_{\alpha'}$ such that $\tilde{f}|_{V}=f$. 

\textit{RF4.} For any $p,q\in K_{\alpha'}$, $\sup\{\frac{|f(p)-f(q)|^2}{\mathcal{E}_{\alpha'}(f)}: f\in \mathcal{F}_{\alpha'}, \mathcal{E}_{\alpha'}(f)>0\}$ is finite by the estimate (6.3).

\textit{RF5.} For any $f\in\mathcal{F}_{\alpha'}$, we have $\bar{f}:=\max\{\min\{f,1\},0\}\in \mathcal{F}_{\alpha'}$ and $\mathcal{E}_{\alpha'}(\bar{f})\leq \mathcal{E}_{\alpha'}(f)$.

Thus for each $\alpha'\in \Lambda'$, $(\mathcal{E}_{\alpha'},\mathcal{F}_{\alpha'})$ is indeed a resistance form on $K_{\alpha'}$. Moreover, we have the following claims on  $(\mathcal{E}_{\alpha'},\mathcal{F}_{\alpha'})$.
 
\textit{Claim 1. $\mathcal{E}_{\alpha'}(f)=\sum_{\beta\in \Lambda'_{\alpha'}}\mathcal{E}_{\beta'}(f|_{K_{\beta'}}),\quad \forall f\in \mathcal{F}_{\alpha'}.$}
 
This follows directly from 
\[\mathcal{E}_{\alpha'}(f)=\mu_{\mathcal{E},u}(K_{\alpha'})=\sum_{\beta'\in \Lambda'_{\alpha'}}\mu_{\mathcal{E},u}(K_{\beta'})=\sum_{\beta'\in \Lambda'_{\alpha'}} \mathcal{E}_{\beta'}(f|_{K_{\beta'}}),\]
with $u|_{K_{\alpha'}}=f.$\hfill$\square$

Let $L_{\alpha'}\subset\Lambda$ be an at most countable set of indices satisfying (6.1) and (6.2). Let $A$ be the set of accumulation points of $\bigcup_{\alpha\in L_{\alpha'}}V_\alpha$. We write $\mathcal{F}_{\alpha', A}:=\{f\in\mathcal{F}_{\alpha'}: f \text{ is constant in a neighborhood of } q, \forall q\in A\}.$ Obviously, for a function $f\in \mathcal{F}_{\alpha', A}$, $f$ is non-constant on only finitely many islands $K_\alpha$ with $\alpha\in L_{\alpha'}$.

{\textit{Claim 2. $\mathcal{F}_{\alpha',A}$ is dense in $(\mathcal{F}_{\alpha'},\mathcal{E}_{\alpha'})$ in energy, i.e., for any $f\in \mathcal{F}_{\alpha'}$ and any $\epsilon>0$, there exists $g\in\mathcal{F}_{\alpha',A}$ such that $\mathcal{E}_{\alpha'}(f-g)<\epsilon$.}

Let $f\in\mathcal{F}_{\alpha'}$ and $q\in A$. For $\epsilon>0$, we first prove that there is a function $\tilde{f}\in \mathcal{F}_{\alpha'}$ which is constant in a neighborhood of $q$, such that $\mathcal{E}_{\alpha'}(f-\tilde{f})<\epsilon$. In fact, from the definition of $\mathcal{F}_{\alpha'}$, there is a function $u\in\mathcal{F}$ such that $f=u|_{K_{\alpha'}}$. Choose a neighborhood $U_q$ of $q$, which is a finite union of islands $K_\beta$ with $\beta\in\Lambda$ such that 
\[\mu_{\mathcal{E},u}(U_q)<\frac{\epsilon}{4}.\]  
Note that $U_q\cap \overline{K\setminus U_q}<\infty$. Let $\bar{u}\in \mathcal{F}$ such that $\bar{u}|_{\overline{K\setminus U_q}}=u|_{\overline{K\setminus U_q}}$, $\bar{u}(q)=u(q)$, and $\bar{u}$ is harmonic in $U_q\setminus\{q\}$. Obviously, $\mu_{\mathcal{E},\bar{u}}(U_q)<\frac{\epsilon}{4}$. Then consider a decreasing nested sequence of neighborhoods $\{U_{q,n}\}_{n\geq 1}$ of $q$ contained in $U_q$, shrinking to $q$, assuming each of  $\{U_{q,n}\}_{n\geq 1}$ is also a finite union of islands $K_\beta$ with $\beta\in \Lambda$. For $n\geq 1$, define $\bar{u}_n\in\mathcal{F}$ by taking
\[\begin{cases}
\bar{u}_n|_{\overline{K\setminus U_q}}=u|_{\overline{K\setminus U_q}},\\
\bar{u}_n|_{{U}_{q,n}}=u(q),\\
\bar{u}_n \text{ is harmonic in } U_q\setminus {U}_{q,n}.
\end{cases}\]
Since $\bar{u}_n$ converges to $\bar{u}$ in energy, we could choose $n$ sufficient large such that
\[\mu_{\mathcal{E},\bar{u}_n}(U_q)<\frac{\epsilon}{4}.\] Denote this $\bar{u}_n$ by $\tilde{u}$ and write $\tilde{f}=\tilde{u}|_{K_{\alpha'}}$, then we have
\[\mathcal{E}_{\alpha'}(f-\tilde{f})\leq \mu_{\mathcal{E},u-\tilde{u}}(U_q)<\epsilon.\]

Noticing that $A$ consists of at most countably many points, we write  $A=\{q_1,q_2,\cdots\}$. Using the above argument, for $f\in\mathcal{F}_{\alpha'}$ and $\epsilon>0$, we could inductively construct a sequence of functions $\{f_n\}_{n\geq 1}$ in $\mathcal{F}_{\alpha'}$ with $f_0=f$, and for $n\geq 1$, $f_n$ is constant in a neighborhood of $q_k$, $\forall 1\leq k\leq n$,  and  $\mathcal{E}_{\alpha'}(f_{n-1}-f_n)<\frac{\epsilon}{4^n}$. Then $\{f_n\}_{n\geq 1}$ forms a Cauchy sequence in energy and the limit function $g$ of the sequence satisfies $g \in \mathcal{F}_{\alpha',A}$ and $\mathcal{E}_{\alpha'}(f-g)<\epsilon$. Thus we have proved Claim 2.\hfill$\square$

\textit{Claim 3. Suppose $\alpha'\in \mathcal{T}_j'$, $1\leq j\leq M'$, then $f\in \mathcal{F}_{\alpha'}$ if and only if $f\circ \phi_{\alpha'_j,\alpha'}\in \mathcal{F}_{\alpha'_j}$. In addition, for any $f\in \mathcal{F}_{\alpha'}$, 
\[
\mathcal{E}_{\alpha'}(f)=r^{-1}_{\bm{e}(\vartheta',\alpha')} \mathcal{E}_{\alpha'_j} (f\circ \phi_{\alpha'_j,\alpha'}),
\]
where $r_{\bm{e}(\vartheta',\alpha')}:=c\big(\frac{diam(K_{\alpha_j'})}{diam(K_{\alpha'})}\big)^{-1}$  with $c(\cdot)$ being the function introduced in Proposition 6.2. }

Since $S'\blacktriangleleft S$, we could choose $L_{\alpha'}$, $L_{\alpha_j'}\subset\Lambda$ satisfying (6.1) and (6.2) separately, and require that there is a one to one correspondence $p: L_{\alpha'}\rightarrow L_{\alpha_j'}$ such that
\begin{equation}
\phi_{\alpha,p(\alpha)}=\phi_{\alpha',\alpha_j'}, \quad\forall \alpha\in L_{\alpha'}.
\end{equation}
Let $A$ be the set of accumulation points of $\bigcup_{\alpha\in L_{\alpha'}}V_{\alpha}$ and denote by $\mathcal{F}_{\alpha',A}$ as in Claim 2. Write $p(A)$ the set of accumulation points of $\bigcup_{\alpha\in L_{\alpha'}}V_{p(\alpha)}$ for simplicity. Then for $f\in\mathcal{F}_{\alpha',A}$, by using Proposition 6.2 and (6.4), we have
\[\begin{aligned}
\mathcal{E}_{\alpha'}(f)
&=\sum_{i=1}^M\sum_{\alpha\in L_{\alpha'},\alpha\in\mathcal{T}_i}r^{-1}_{\bm{e}(\vartheta,\alpha)}\mathcal{E}_{\alpha_{i}}(f|_{K_{\alpha}}\circ \phi_{\alpha_{i},\alpha})\\
&=\sum_{i=1}^M\sum_{\alpha\in L_{\alpha'}, \alpha\in\mathcal{T}_i}r^{-1}_{\bm{e}(\vartheta,\alpha)}\mathcal{E}_{\alpha_{i}}\big(f|_{K_{\alpha}}\circ \phi_{p(\alpha),\alpha}\circ \phi_{\alpha_{i},p(\alpha)}\big)\\
&=c\big(\frac{diam(K_{\alpha'_j})}{diam(K_{\alpha'})}\big)\sum_{i=1}^M\sum_{\alpha\in L_{\alpha'},\alpha\in\mathcal{T}_i}r^{-1}_{\bm{e}\big(\vartheta,p(\alpha)\big)}\mathcal{E}_{\alpha_{i}}\big((f\circ \phi_{\alpha'_j,\alpha'})|_{K_{p(\alpha)}}\circ \phi_{\alpha_{i},p(\alpha)}\big)\\
&=r^{-1}_{\bm{e}(\vartheta',\alpha')} \mathcal{E}_{\alpha_j'} (f\circ \phi_{\alpha'_j,\alpha'}).
\end{aligned}\]
The last equality is valid since we can easily see that $f\circ \phi_{\alpha'_j,\alpha'}\in {\mathcal{F}}_{\alpha'_j,p(A)}$. Then Claim 3 holds by using Claim 2. \hfill$\square$

For each $1\leq j\leq M'$, let $D_{\alpha_j'}$ be the Laplacian induced by  the trace of $\mathcal{E}_{\alpha'_j}$ on $V_{\alpha'_j}$. Then, in view of  Claim 1 and Claim 3,  $(\{D_{\alpha'_j}\}_{j=1}^{M'},\bm{r}')$ is a homogenous regular harmonic structure of the $f.r.f.t.$ nested structure $\mathcal{S}'=\{K_{\alpha'},\Lambda_{\alpha'}'\}_{\alpha'\in \Lambda'}$ with  $r_{\bm{e}(\vartheta',\alpha')}:=c\big(\frac{diam(K_{\alpha_j'})}{diam(K_{\alpha'})}\big)^{-1}$ for any $\alpha'\in\mathcal{T}_{j}', 1\leq j\leq M'$.
\hfill$\square$

\textbf{Remark.} In the case that $\mathcal{S'}\not\blacktriangleleft \mathcal{S}$, the conclusion of Theorem 6.4 may fail to hold. For example, let $0<r<1$, consider the canonical $f.r.f.t.$ nested structure on the line segment $I=[0,1]$, denoted by $\mathcal{S}_r$,  generated by the $i.f.s.$, 
\[F_1: x\rightarrow rx,\quad F_2: x\rightarrow rx+1-r.\]
Suppose $r$ is not an algebraic number, then any harmonic structure on $\mathcal{S}_r$ is homogeneous, but there is only one of them inducing the same resistance form as that of  the  homogeneous regular harmonic structure on $\mathcal{S}_{1/2}$.\\

\section{Spectral asymptotics}
In this section, we briefly discuss the asymptotics of the eigenvalue counting functions associated with the Dirichlet forms on the $f.r.f.t.$ self-similar sets.  
Let $K$ be an $f.r.f.t.$ self-similar set.  Throughout this section, we assume that there exists a regular harmonic structure $(\{D_{\alpha_i}\}_{i=1}^M,\bm{r})$ with respect to an $f.r.f.t.$ nested structure  $\{K_\alpha,\Lambda_\alpha\}_{\alpha\in \Lambda}$ of $K$.  Let $G=(S,E)$ be the directed graph of the corresponding graph-directed construction. 

Let $\mu$ be a Borel probability measure  on $K$, $(\mathcal{E},\mathcal{F})$ be the self-similar Dirichlet form on $L_\mu^2(K)$ generated by $(\{D_{\alpha_i}\}_{i=1}^M,\bm{r})$ and $\mu$.  Let $\Delta_\mu$ be the associated \textit{Laplacian} of  $(\mathcal{E},\mathcal{F})$ on $K$. For $u\in\mathcal{F}$, the Laplacian $\Delta_\mu u$ of $u$ is a continuous function satisfying
\[\mathcal{E}(u,v)=-\int_K v\Delta_\mu u d\mu,\quad\forall v\in \mathcal{F}_0,\]
with $\mathcal{F}_0$ being the collection of functions in $\mathcal{F}$ vanishing at the vertices in $V_0$. Write $dom\Delta_\mu$  the domain of the Laplacian $\Delta_\mu$.

We are interested in the eigenfunctions and eigenvalues of the Laplacian $\Delta_\mu$. A function $u\in dom\Delta_\mu$ is said to be a \textit{Dirichlet eigenfunction} associated with an \textit{eigenvalue} $\lambda$ providing that it is a non-zero solution of 
\[\begin{cases}
-\Delta_\mu u=\lambda u,\\
u|_{V_0}=0.
\end{cases}\]
Equivalently, $u$ is a Dirichlet eigenfunction associated with an eigenvalue $\lambda$ of $\Delta_\mu$ if $u\in \mathcal{F}_0$ and 
\[\mathcal{E}(u,v)=\lambda\int_K uvd\mu,\quad \forall v\in \mathcal{F}_0.\]
Similarly, we say that $u$ is a \textit{Neumann eigenfunction} associated with an eigenvalue $\lambda$ if $u\in \mathcal{F}$ and 
\[\mathcal{E}(u,v)=\lambda\int_K uvd\mu,\quad \forall v\in \mathcal{F}.\]

It can be shown in a routine way that there is a discrete sequence of eigenvalues, call the \textit{spectrum} of the Laplacian $\Delta_\mu$, with $0\leq \lambda_1\leq \lambda_2\leq \cdots$, whose only accumulation point is $\infty$, for either the Dirichlet case or the Neumann case. So we can define the eigenvalue counting functions for the Dirichlet and Neumann eigenvalues as
\[N_\mathcal{D}(x)=\#\{\lambda\leq x: \lambda \text{ is a Dirichlet eigenvalue of }\Delta_\mu\},\] 
and
\[N_\mathcal{N}(x)=\#\{\lambda\leq x: \lambda \text{ is a Neumann eigenvalue of }\Delta_\mu\}.\]
We will consider the asymptotics of these eigenvalue counting functions as $x$ goes to the infinity.

We require the measure $\mu$ to \textit{match} the graph-directed construction of $K$. To be precise, we define a Borel probability measure $\mu$ on $K$  in the following way. For each edge $e\in E$, we assign a weight $\mu_e\in(0,1)$ such that $\sum_{i(e)=s} \mu_e=1$ for any state $s\in S$, and let $\mu_{\bm{e}}=\prod_{i=1}^{|\bm{e}|}\mu_{e_i}$ for any walk $\bm{e}$ of finite length. For any $\beta\in\Lambda$, define $\mu(K_\beta)$ as 
\[\mu(K_\beta)=\mu_{\bm{e}(\vartheta,\beta)}, \]
where $\bm{e}(\vartheta,\beta)$ is the same as defined in Section 5. In the case that the corresponding $f.r.g.d.$ fractal family of $K$ satisfies the open set condition in the sense of Mauldin and Williams [MW], the normalized Hausdorff measure $\nu$ satisfies the above construction. Indeed, this is true for the four examples in Section 4. 

For $e\in E$, we write $\mu_e$ the weight corresponding to the measure $\mu$. For $\delta>0$, we introduce a $\# S\times\#S$-matrix $M_\delta$ with entries given by $M_\delta(i,j)=\sum_{\{e: i(e)=i, f(e)=j\}}(r_e\mu_e)^\delta$, and let $\Psi(\delta)$ denote the spectral radius of the matrix $M_\delta$. Obviously, there exists a unique positive real number $\delta$ such that $\Psi(\delta)=1$.

In [HN], Hambly and Nyberg have made a thorough analysis for the spectral asymptotics for Laplacians on $f.r.g.d.$ fractal families with respect to graph-directed measures. Since by Theorem 3.4, an $f.r.f.t.$ self-similar set is naturally an $f.r.g.d.$ fractal, the spectral asymptotics for $\Delta_\mu$ on $K$ then follows directly. In the following, we only state the result for the spectral asymptotic ratio in the case that the graph $G=(S,E)$ is strongly connected, or equivalently, the matrix $M_\delta$ is irreducible.

\textbf{Theorem 7.1.} \textit{Let $K$ be a $f.r.f.t.$ self-similar set, whose associated directed graph is strongly connected. Let $\delta$ be the positive real number such that $\Psi(\delta)=1$. Then
\[\begin{aligned}
0<\liminf_{x\rightarrow\infty} N(x)x^{-\delta}\leq \limsup_{x\rightarrow\infty} N(x)x^{-\delta}<\infty,
\end{aligned}\]
where $N(x)$ stands for either $N_{\mathcal{D}}(x)$ or $N_{\mathcal{N}}(x)$.}

We remark that the limit $\lim_{x\rightarrow\infty} N(x)x^{-\delta}$ may not exist. The leading-order term in the asymptotic expansion of $N(x)$ is either a constant or a periodic function. In the case that the directed graph is not strongly connected, things will be more complicated. See more precise results in [HN], where a multidimensional renewal theorem was established to solve the problem. The proof is inspired by the idea in [KL] dealing with the $p.c.f.$ self-similar sets. 

Due to the well-known result of Hutchinson [H], for any given probability weight vector $\bm{\rho}=(\rho_1,\rho_2,\cdots, \rho_N)$ with $0<\rho_i<1$ and $\sum_{i=1}^{ N}\rho_i=1$, there always exists a probability measure on $K$, called a \textit{self-similar measure} associated with the $i.f.s.$ of $K$, denoted by $\nu_{\bm{\rho}}$, satisfying that
\[\begin{aligned}
\nu_{\bm{\rho}}(A)=\sum_{i=1}^N\rho_i\cdot\nu_{\bm{\rho}}\big(F_i^{-1}(A)\big),\end{aligned}\]
for any Borel set $A$ in $K$. However, in general, it is not expectable that such measures match the graph-directed construction of $K$. Thus the method in [HN] for the spectral asymptotics of the Laplacians is not applicable for the self-similar measures. The following question arises naturally.

\textbf{Question 7.2.} \textit{How about the asymptotic behavior of the eigenvalue counting functions for Laplacians on $f.r.f.t.$ self-similar sets with respected to the self-similar measures?}

We leave the question for further study.

\section{Appendix}
In this appendix, we give two examples, organized as following. 

Example 1. A self similar set with the finite neighboring type property, not satisfying the finite chain length property.

Example 2. A self similar set with the finite chain length property, not satisfying the finite neighboring type property.\\

\textbf{Example 1.(Golden ratio Sierpinski gasket)} 
Let $\{q_i\}_{i=1}^3$ be the three vertices of an equilateral triangle, and $\{F_i\}_{i=1}^3$ be the three contractive similitudes,
\[
\begin{aligned}
F_1: x\rightarrow &\rho^2(x-q_1)+q_1,\\
F_2: x\rightarrow \rho(x-q_2)&+q_2,\quad F_3: x\rightarrow \rho(x-q_3)+q_3,
\end{aligned}
\] 
with $\rho=\frac{\sqrt{5}-1}{2}$. The \textit{golden ratio Sierpinski gasket} $\mathcal{SG}^g$ is the invariant set of the \textit{i.f.s.} $\{F_i\}_{i=1}^3$, i.e., $\mathcal{SG}^g=\bigcup_{i=1}^3 F_i(\mathcal{SG}^g)$, see Figure \ref{goldensg}. It is a slight variant of Example 5.4 in [NW].

\begin{figure}
	\includegraphics[width=5cm]{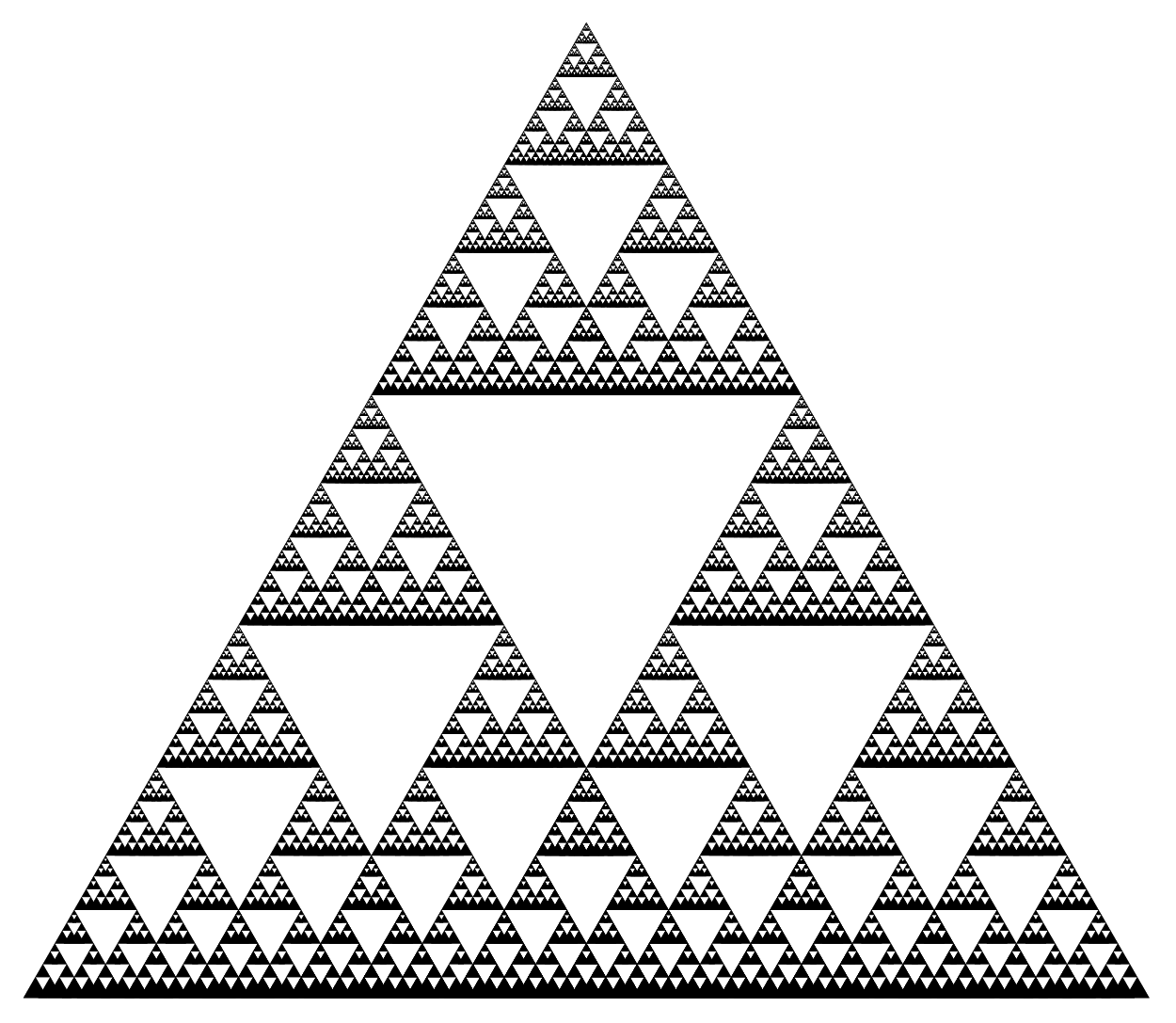}
	\caption{The golden ratio Sierpinski gasket $\mathcal{SG}^g.$}\label{goldensg}
\end{figure} 

Obviously, $\mathcal{SG}^g$ satisfies the finite neighboring type property, see a detailed discussion in [NW]. However, $\mathcal{SG}^g$ does not satisfy the finite chain length property. In fact, consider the collection of copies $\{F_w\mathcal{SG}^g|w\in \{2,3\}^n\}$, $n\geq 1$, located along the bottom line of the fractal. By ordering the words in lexicographical order, i.e., letting $w^{(1)}=22\cdots 2,w^{(2)}=22\cdots 23$, $\cdots$, $w^{(2^n)}=33\cdots3$, and removing the completely overlapping ones, we can find that the collection $\gamma_n=(F_{w^{(1)}}\mathcal{SG}^g, F_{w^{(2)}}\mathcal{SG}^g, \cdots,F_{w^{(2^n)}}\mathcal{SG}^g)$ provides a $\delta$-overlapping chain for any $0<\delta<1$. Since $n$ can be arbitrarily large, $\mathcal{SG}^g$ does not satisfy the finite chain length property.\\

\textbf{Example 2.($\lambda$-gaskets with irrational moving sliders)} Let $\lambda\in (0,1)$. Define the following $i.f.s. \{F_i\}_{i=0}^4$ in $\mathbb{R}^2$,
\[\begin{aligned}
F_0: x\rightarrow \frac{1}{3}x,\quad F_1: x\rightarrow \frac{1}{3}x&+(\frac{1}{3},0),\quad F_2: x\rightarrow \frac{1}{3}x+(\frac{2}{3},0),\\
F_3: x\rightarrow \frac{1}{3}x+(\frac{1}{3},\frac{\sqrt{3}}{3}),&\quad F_4: x\rightarrow \frac{1}{3}x+(\frac{1}{6}+\frac{\lambda}{3},\frac{\sqrt{3}}{6}).
\end{aligned}\]
Let $K_\lambda$ be the invariant set of this \textit{i.f.s.} See Figure \ref{klambda} for an illustration. 
\begin{figure}[h]
\includegraphics[width=5cm]{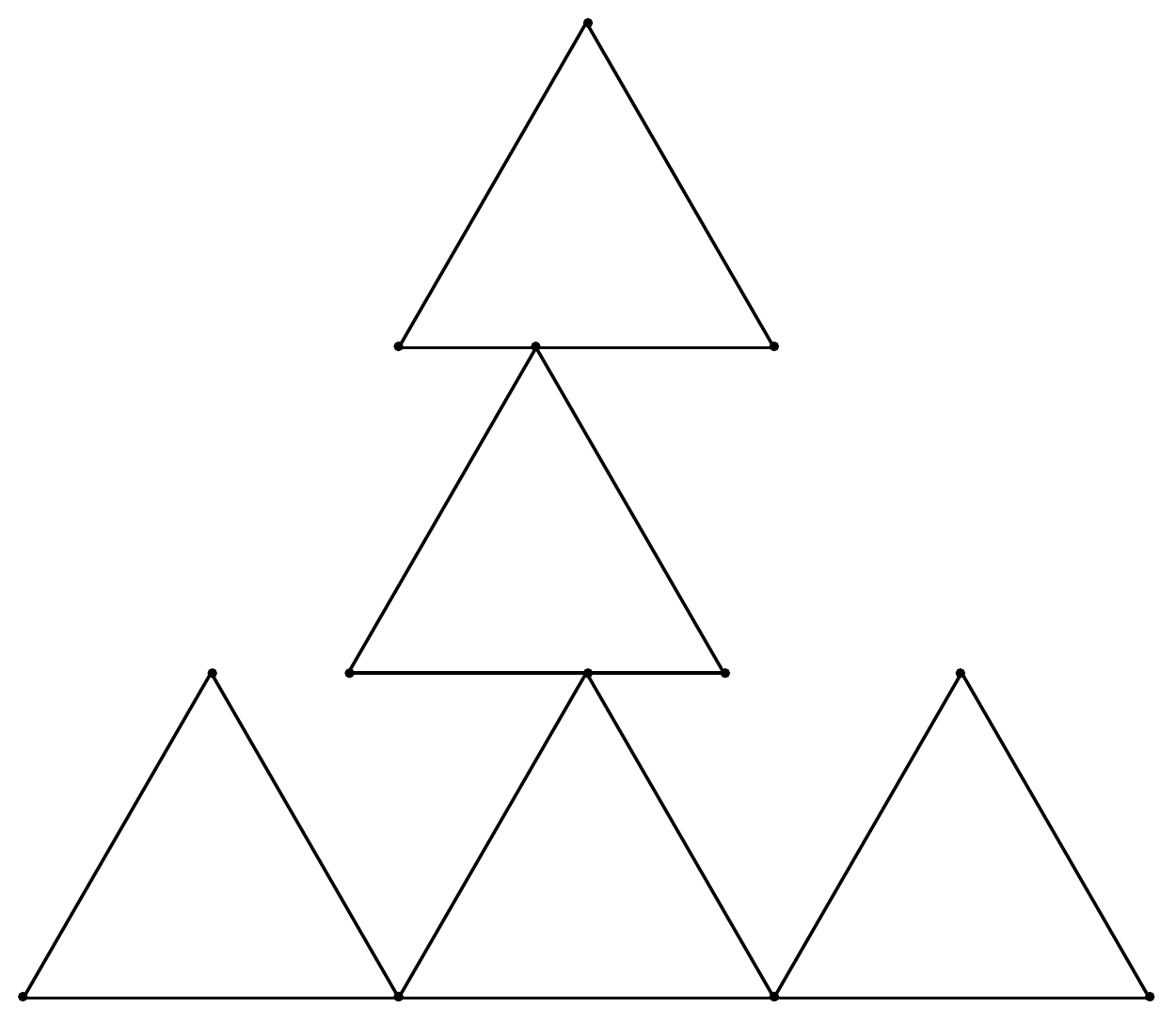}
\begin{picture}(0,0)
\put(-95,73){$\frac{\lambda}{3}$}
\put(-126,13){$F_0$}
\put(-80,13){$F_1$}
\put(-34,13){$F_2$}
\put(-80,93){$F_3$}
\put(-86,53){$F_4$}
\end{picture}
\caption{An illustration for $K_\lambda$.}
\label{klambda}
\end{figure}
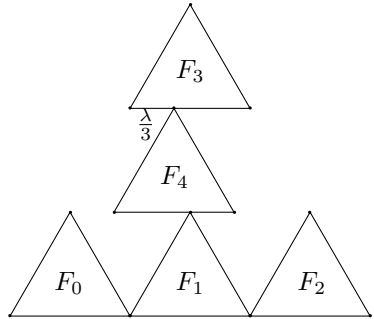

Obviously, $K_\lambda$ satisfies the finite chain length property. In fact, given any two copies $F_wK_\lambda$ and $F_uK_\lambda$, either they intersect each other by at most one point, or one contains the other. So any overlapping chain has length at most $1$.

On the other hand, $K_\lambda$ does not satisfy the finite neighboring property when $\lambda$ is an irrational number. In fact, write the ternary expansion of $\lambda$,
\[\lambda=\sum_{i=1}^{\infty} l_i3^{-i}, \]
with $l_i\in\{0,1,2\},\forall i\geq 1$. By shifting the coefficients in this expansion, we get a sequence of irrational numbers 
\[\lambda_k=\sum_{i=1}^{\infty} l_{i+k}3^{-i}, k\geq 1.\]
For the sake of uniformity,  write $\lambda_0=\lambda$. 
It is not hard to see that for any $k\geq 0$, $F_4F^k_3K_\lambda\cap F_3F_{[l]_k}K_\lambda\neq \emptyset$, where $[l]_k=l_1l_2\cdots l_k\in \{0,1,2\}^k\subset W_*$. Moreover, a calculation yields that 
\[ (F_3F_{[l]_k})^{-1}\circ F_4F^k_3: x\rightarrow x-(\frac{1}{2},\frac{\sqrt{3}}{2})+(\lambda_k,0).\]
Notice that $\lambda_k\neq \lambda_{k'}$ when $k\neq k'$ since $\lambda$ is irrational. This shows that $K_\lambda$ does not satisfy the finite neighboring property. 

This example was introduced in [W](Example 3) for a different purpose. Some adjustment is made in our setting.

\end{document}